\documentclass{amsart}
\usepackage{amssymb,amscd}
\usepackage[all]{xy}
\newcommand\datver[1]{\def\datverp%
 {\par\boxed{\boxed{\text{#1; Run: \today}}}}}

\newcommand\boxb[1]{\square_b}
\newcommand\Ahat{\widehat{A}}
\newcommand\AS{\operatorname{AS}}

\newcommand\ff{\operatorname{ff}}

\numberwithin{equation}{section}

\newcommand\paperbody%
        {\numberwithin{theorem}{section}
        \numberwithin{definition}{section}

	}

\newcommand\fammap{\phi}
\newcommand\famfib{Z}
\newcommand\famtot{M}
\newcommand\fambas{B}
\newcommand\sfamrel{\efamtot/\famtot}
\newcommand\famrel{\famtot/\fambas}
\newcommand\tfammap{\phi'}
\newcommand\sfamfib{F}
\newcommand\bfamtot{M'}
\newcommand\tfamrel{\bfamtot/\famtot}
\newcommand\bfamrel{\bfamtot/\fambas}
\newcommand\bboumap{\Phi'}
\newcommand\bboubas{D'}
\newcommand\efammap{\tphi}
\newcommand\efamfib{Z}
\newcommand\efamtot{\tM}
\newcommand\efambas{B}
\newcommand\efamrel{\efamtot/\efambas}
\newcommand\efamext{(\efammap/\fammap)^{!}}
\newcommand\mfammap{\pi}

\newcommand\mfamtot{\fambas\times\bbB^{p+1}}

\newcommand\Mfammap{\pi}

\newcommand\Mfamtot{\fambas\times\bbB^{p+1}\times\bbS^q}

\newcommand\boufib{X}
\newcommand\boufibmap{\psi}
\newcommand\boufibbas{Y}
\newcommand\boumap{\Phi}
\newcommand\bourel{\pa\famtot/\boubas}
\newcommand\boubasrel{\boubas/\fambas}

\newcommand\boubas{D}

\newcommand\eboumap{\tPhi}

\newcommand\eboubas{\tD}

\newcommand\Mboumap{\pi'}

\newcommand\Mboubas{\efambas\times\bbS^p}

\newcommand\Kfn[3]{\operatorname{KK}^{#1}_{#2}(#3)}
\newcommand\qr{ \quan_{r}}
\newcommand\qb{ \quan_{-\infty}}
\newcommand\qcn{ \quan_{c}}
\newcommand\icn{ i_{c}}
\newcommand\rcn{ \iota_{c}}
\newcommand\ecn{ e_{c}}
\newcommand\cn{ M_{c}}
\newcommand\cnbm{\Phi_{c}}
\newcommand\cnfm{\phi_{c}}

\newcommand\mc{\mathcal{C}_{\iota}}
\newcommand\cnf{\pi_{\Phi}}
\newcommand\ccn{\varepsilon_{c}}
\newcommand\dcn{\varepsilon_{d}}

\newtheorem{theorem}{Theorem}

\newtheorem{proposition}[theorem]{Proposition}
\newtheorem{corollary}[theorem]{Corollary}
\newtheorem{lemma}[theorem]{Lemma}

\theoremstyle{definition}
\newtheorem{definition}{Definition}

\newtheorem{remark}[definition]{Remark}

\newcommand\sE{\mathbb{E}}
\newcommand\sF{\mathbb{F}}
\newcommand\sG{\mathbb{G}}
\newcommand\sH{\mathbb{H}}

\newcommand\bbB{\mathbb B}
\newcommand\bbC{\mathbb C}
\newcommand\bbD{\mathbb D}
\newcommand\bbE{\mathbb E}

\newcommand\bbL{\mathbb L}

\newcommand\bbR{\mathbb R}
\newcommand\bbS{\mathbb S}

\newcommand\bbZ{\mathbb Z}

\newcommand\CI{{\mathcal{C}}^{\infty}}
\newcommand\COS{\mathcal{S}\mathcal{C}}

\newcommand\B{\mathcal{B}}

\newcommand\F{\mathcal{F}}
\newcommand\K{K^{0}_{c}}
\newcommand\Kv{\mathcal{K}}

\newcommand\Hv{\mathcal{H}}

\newcommand\Pt{\operatorname{p}}
\newcommand\Ps[2]{\Psi^{#1}(#2)}
\newcommand\fsP[3]{\Psi_{#1-\sus}^{#2}(#3)}
\newcommand\ptP[3]{\Psi_{#1-\Pt}^{#2}(#3)}
\newcommand\ptPS[3]{\mathcal{S}_{#1-\Pt}^{#2}(#3)}
\newcommand\sptP[3]{\Psi_{\sus,#1-\Pt}^{#2}(#3)}
\newcommand\sfptP[4]{\Psi_{#1-\sus,#2-\Pt}^{#3}(#4)}

\newcommand\fcP[3]{\Psi_{\fC{#1}}^{#2}(#3)}
\newcommand\fC[1]{\operatorname{#1-\cusp}}
\newcommand\aD[1]{\operatorname{#1-ad}}
\newcommand\ad{\operatorname{ad}}
\newcommand\adP[3]{\Psi_{\aD{#1}}^{#2}(#3)}
\newcommand\sadP[4]{\Psi_{#1-\sus,\aD{#2}}^{#3}(#4)}
\newcommand\syM[3]{S^{#2}_{#1}(#3)}
\newcommand\adfcP[4]{\Psi_{\aD{#2},\fC{#1}}^{#3}(#4)}
\newcommand\adcuP[3]{\Psi_{\aD{#1},\cusp}^{#2}(#3)}
\newcommand\cuP[2]{\Psi_{\cusp}^{#1}(#2)}
\newcommand\ptfcP[4]{\Psi_{#1-\Pt,\fC{#2}}^{#3}(#4)}

\newcommand\ptcuP[3]{\Psi_{#1-\Pt,\cusp}^{#2}(#3)}
\newcommand\ptcuPS[3]{\mathcal{S}_{#1-\Pt,\cusp}^{#2}(#3)}

\newcommand\fcs{{}^{\fC{\boumap}}S}

\newcommand\fsn[1]{\Psi_{\boumap s(#1)}}
\newcommand\jS[2]{A_{#1}^{#2}}

\newcommand\IJS[2]{J_{#1}^{#2}}
\newcommand\ptjS[3]{A_{#1-\Pt}^{#2}(#3)}
\newcommand\ptcujS[2]{A_{#1-\Pt,\cusp}^{#2}}
\newcommand\IUS[2]{G_{#1}^{#2}}

\newcommand\Gsfc{G^{-\infty}_{\boumap s}}
\newcommand\Gsfcn[1]{G^{-\infty}_{\boumap s(#1)}}
\newcommand\Gsbn[1]{G^{0}_{\boumap s(#1)}}

\newcommand\Gs{G^{-\infty}_{\sus}}
\newcommand\Gsp{G^{-\infty}_{\esus}}

\newcommand\quan{\operatorname{quan}}
\newcommand\quanp{\operatorname{quan}_{a}}
\newcommand\ind{\operatorname{ind}}
\newcommand\indAS{\operatorname{ind_{\AS}}}
\newcommand\inda{\operatorname{ind_{\text{a}}}}
\newcommand\indt{\operatorname{ind_{\text{t}}}}
\newcommand\inv{\operatorname{inv}}
\newcommand\Thom{\operatorname{Thom}}
\newcommand\Id{\operatorname{Id}}
\newcommand\Ld{\operatorname{L}^{2}}

\newcommand\pa{\partial}
\newcommand\pt{\operatorname{pt}}

\newcommand\KK{\operatorname{KK}}
\newcommand\Kf[2]{\operatorname{KK}^{0}_{#1}(#2)}
\newcommand\Kfo[2]{\operatorname{KK}^{1}_{#1}(#2)}
\newcommand\Ef[2]{\bbE_{#1}(#2)}
\newcommand\Kc{\operatorname{K_{\text{c}}}}
\newcommand\Kcn[1]{\operatorname{K^{#1}_{\text{c}}}}
\newcommand\Kt{\operatorname{K}}
\newcommand\tKt{\operatorname{\widetilde{K}}}
\newcommand\Kto{\operatorname{K^1}}
\newcommand\Kco{\operatorname{K^1_{\text{c}}}}
\newcommand\cuK[1]{\operatorname{K_{\cusp}}(#1)}

\newcommand\scK{\operatorname{K_{\scat}}}
\newcommand\fcK[1]{\operatorname{K_{\operatorname{#1-\cusp}}}}
\newcommand\fcS[2]{\mathcal{S}_{\operatorname{#1-\cusp}}^{#2}}
\newcommand\fcKn[2]{\operatorname{K^{#1}_{\operatorname{#2-\cusp}}}}
\newcommand\fcKo[1]{\operatorname{K^1_{\operatorname{#1-\cusp}}}}
\newcommand\dc{\eth_{\cusp}}
\newcommand\cuS{{}^{\cusp}S}

\newcommand\cuT{{}^{\cusp}T}
\newcommand\scT{{}^{\scat}T}
\newcommand\fcT[1]{{}^{\operatorname{#1-\cusp}}T}
\newcommand\Cfc{\mathcal{C}_{\boumap}}
\newcommand\Ccu{\mathcal{C}_{\cusp}}

\newcommand\scat{\operatorname{sc}}

\newcommand\bo{\operatorname{b}}

\newcommand\cusp{\operatorname{cu}}

\newcommand\coF{{}^{\mathcal{C}}\kern-2pt\Lambda}

\newcommand\cFTs{{}^{\Phi}\overline{T}\kern-1pt{}^*}

\newcommand\esus{\operatorname{e-sus}}
\newcommand\sus{\operatorname{sus}}

\newcommand\Ch{\operatorname{Ch}}

\newcommand\tphi{\tilde{\phi}}

\newcommand\tM{\widetilde M}
\newcommand\tD{\widetilde D}

\newcommand\tPhi{\widetilde \Phi}
\newcommand\com[1]{\overline{#1}}

\newcommand\ie{i\@.e\@. }

\newcommand\cA{\mathcal{A}}
\newcommand\cB{\mathcal{A}'}

\newcommand\cS{\mathcal S}
\newcommand\CIc{{\mathcal{C}}^{\infty}_c}

\newcommand\CO{\mathcal{C}}
\newcommand\COo{\mathcal{C}_{0}}

\newcommand\Diag{\operatorname{Diag}}

\newcommand\ci{${\mathcal{C}}^\infty$}

\newcommand\dCI{\dot{\mathcal{C}}^{\infty}}

\newcommand\nul{\operatorname{null}}

\newcommand\longhookrightarrow{\hookrightarrow}

\newcommand\Mand{\text{ and }}

\newcommand\Mat{\text{ at }}

\newcommand\Mif{\text{ if }}

\newcommand\Min{\text{ in }}

\newcommand\Mwhere{\text{ where }}

\datver{1.0J; Revised: 04-03-2006}

\begin{document}
\title[K-theory of cusp operators]
{Index in K-theory for Families\\ of fibred cusp operators}

\author{Richard Melrose AND Fr\'ed\'eric Rochon}
\address{Department of Mathematics, Massachusetts Institute of Technology}
\email{rbm@math.mit.edu}
\address{Department of Mathematics, State University of New York, Stony Brook}
\email{rochon@math.sunysb.edu}

\begin{abstract} A families index theorem in K-theory is given for the
setting of Atiyah, Patodi and Singer of a family of Dirac operators with
spectral boundary condition. This result is deduced from such a K-theory
index theorem for the calculus of cusp, or more generally fibred cusp,
pseudodifferential operators on the fibres (with boundary) of a fibration;
a version of Poincar\'e duality is also shown in this setting, identifying
the stable Fredholm families with elements of a bivariant K-group.
\end{abstract}
\thanks\datverp
\maketitle

\tableofcontents

\section*{Introduction}

The index theorem of Atiyah, Patodi and Singer gives a formula for the
index of a Dirac operator on a compact manifold with boundary with boundary
condition given by projection onto the range of the positive part of the
boundary Dirac operator. Versions of this result for families, with the
formula being in cohomology for the Chern character of the (virtual) index
bundle, were given by Bismut-Cheeger \cite{Bismut-Cheeger1,
Bismut-Cheeger2} and the first author and Piazza
\cite{MR99a:58144,MR99a:58145}. Here we formulate an index theorem in
K-theory in the wider context of the algebras of pseudodifferential of
fibred-cusp type, so generalizing the K-theory formulation of Atiyah and
Singer in the boundaryless case. Our result specializes to give an index
theorem in K-theory for the families of Dirac operators in
the earlier contexts cited above. In a subsequent paper we will show how to
derive a formula for the Chern character of the index class, reducing to
the known formula in the Dirac case. There is also a relation, discussed
below, with the `direct' proof of the theorem of Atiyah, Patodi and Singer
given by Dai and Zhang in \cite{Dai-Zhang1}.

We consider a smooth fibration of compact manifolds where the fibres are
manifolds with boundary, with the boundary possibly carrying a finer
fibration. With the addition of a normal trivialization of the boundary
along the fibres we call this a fibration with fibred cusp structure. For
such objects we formulate a notion of K-theory, denoted
$\fcK{\boumap}(\fammap)$ (where $\fammap$ is the overall fibration and
$\boumap$ the boundary fibration) as the stable homotopy classes of the
Fredholm families of corresponding `fibred-cusp' pseudodifferential
operators (introduced by Mazzeo and the first author in \cite{Mazzeo-Melrose4}) on the fibres. 
In the boundaryless
case this reduces to the compactly supported K-theory of the fibre
cotangent bundle as in \cite{MR43:5554}. The analytic index arises
from the realization as Fredholm operators either through perturbation to
make the null spaces have constant rank or through Kasparov's bivariant
K-theory. In close analogy with the boundaryless case we define a
topological index map by using an embedding of the fibration in the product
of the base and a ball and then we show the equality of analytic and
topological index homomorphisms

\begin{theorem}\label{faficu.18} For a fibration with fibred cusp structure, the
analytic and topological index maps, to K-theory of the base,
\begin{equation}
\xymatrix@1{\fcK{\boumap}(\fammap)
\ar@<.5ex>[r]^(.55){\inda}\ar@<-.5ex>[r]_(.55){\indt}&\Kt(\fambas)}
\label{faficu.19}\end{equation}
are equal.
\end{theorem}

We also give an analogue of Atiyah's Poincar\'e duality in K-theory.

\begin{theorem}\label{faficu.88} For a fibration with fibred cusp
structure, there is a natural `quantization' isomorphism 
\begin{equation}
\quan:\fcK{\boumap}(\fammap)\longrightarrow
\Kf{\fambas}{\CO_{\boumap}(\famtot),\CO(\fambas)}
\label{faficu.62}\end{equation}
where $\CO_{\boumap}(\famtot)$ is the $C^*$ algebra of those continuous functions
on the total space of the fibration which are constant on the fibres of
$\boumap$ on the boundary.  
\end{theorem}
A brief review of $KK$-theory and a definition of  $\Kf{\fambas}{\CO_{\boumap}(\famtot),\CO(\fambas)}$
can be found in section~\ref{KK.0} below.  The `cusp' case of these results (meaning
$\boumap=\pa\fammap=\fammap\big|_{\pa\famtot}$ and in which case we write
the K-group as $\cuK{\fammap})$ applies to give a families index theorem in
K-theory for problems of Atiyah-Patodi-Singer type.

\begin{theorem}\label{faficu.153} Let $\eth$ be a family of Dirac operators
  associated to a unitary Clifford module for a family of metrics of product
  type near the boundary of the fibres of a fibration and suppose $P$ is a
  spectral section for the boundary Dirac operator, then $(\eth, P)$
  defines a class $[(\eth,P)]\in\cuK{\fammap}$ and
\begin{equation*}
         \ind (\eth,P)= \ind_{a}([\eth,P)]).                
\end{equation*}
\end{theorem}
This is a refinement at a K-theoretic level of the index theorem of
\cite{MR99a:58144} for families Dirac operators with Atiyah-Patodi-Singer 
boundary conditions.
Note that \emph{only} the cusp case is needed for this application, so the
  proof requires only a relatively small part of the discussion in the
  body of the paper, in particular
  Sections~\ref{stes.0}--\ref{Adiabaticpassage} are not required for this.

Since we generalize it here, let us briefly recall the families index
in the boundaryless case. The index theorem in K-theory of Atiyah
and Singer takes the form of the equality of an analytically defined and a
topologically defined index for a family of elliptic pseudodifferential
operators $P\in\Psi^m(\famrel;\sE)$ on the fibres of a fibration
\begin{equation}
      \xymatrix{\famfib\ar@{-}[r]&\famtot\ar[d]^{\fammap}\\&\fambas}
\label{faficu.6}\end{equation}
where $\sE=(E_+,E_-)$ is a $\bbZ_2$-graded complex vector bundle over the
total space, $\famtot,$ of the fibration and $P$ maps sections of $E_+$ to
sections of $E_-.$ The analytic index is the element in $\Kt(\fambas)$
which is the formal difference of vector bundles
\begin{equation}
\inda(P)=[\nul(P+A)\ominus\nul(P^*+A^*)]\in \Kt(\fambas)
\label{faficu.10}\end{equation}
for a perturbation $A\in\Psi^{-\infty}(\famrel;\sE)$ such that the null spaces
have constant rank (and for any choice of data leading to the
adjoints). Such a perturbation always exists and $\inda(P)\in \Kt(\fambas)$ is
independent of choices. The vanishing of the index is equivalent to the
existence of such a perturbation with $P+A$ invertible. The symbol of the
family $P$ defines an element in the (compactly supported) K-theory of the
fibre cotangent bundle, $[\sigma(P)]\in \Kc(T^*(\famrel))$ and the analytic
index of $P$ depends only on $[\sigma (P)].$ Since all K-classes with
compact support on $T^*(\famrel)$ arise in this way, the analytic index
gives a map
\begin{equation}
\inda:\Kc(T^*(\famrel))\longrightarrow \Kt(\fambas).
\label{faficu.7}\end{equation}

Alternatively, the analytic index may be defined via the bivariant K-theory
of Kasparov, consisting of equivalence classes of almost-unitary Fredholm
modules. Thus, the choice of a selfadjoint family
$Q\in\Psi^{-m}(\famrel;E_+),$ such that $Q^2P^*P-\Id$ is smoothing, fixes a class
\begin{equation}
\left[\begin{pmatrix}0&QP^*\\PQ&0\end{pmatrix}\right]\in
\Kf{\fambas}{\CO(\famtot),\CO(\fambas)}
\label{faficu.29}\end{equation}
which also depends only on the class of the symbol, \ie defines a
homomorphism which combined with the natural push-forward map 
again gives the analytic index
\begin{equation}
\xymatrix{\Kc(T^*(\famrel))\ar[rr]^{\inda}\ar[rd]
&&\Kt(\fambas)=\Kf{\fambas}{\CO(\fambas),\CO(\fambas)}\\
&\Kf{\fambas}{\CO(\famtot),\CO(\fambas)}\ar[ru]&}.
\label{faficu.30}\end{equation}
In fact the map on the left is an isomorphism which is a realization of
Atiyah's map from elliptic pseudodifferential operators to K-homology, \ie
is a form of Poincar\'e duality.

The topological index is defined as a Gysin map, using Bott
periodicity. By a standard generalization of Whitney's embedding theorem
, the fibration \eqref{faficu.6} may be embedded as a
subfibration of a real vector bundle $V$ over $\fambas$ (indeed the bundle may be
taken to be a product); the K-theory of $T^*(V/\fambas)$ is then canonically
isomorphic to the K-theory of the base. The composite map 
\begin{equation}
\xymatrix{&\Kc(T^*(V/\fambas))\ar[rd]^{\Thom}&\\
\Kc(T^*(\famrel))\ar[rr]^{\indt}\ar[ru]&&\Kt(\fambas)}.
\label{faficu.8}\end{equation}
is the `topological index' and is again independent of choices. The index
theorem of Atiyah and Singer \cite{Atiyah-Singer1}, \cite{MR43:5554}
in K-theory is the equality of these two maps;
the one obtained by `quantizing' symbols by use of pseudodifferential
operators, the other by `trivializing' the symbols using embeddings. 

In this paper the corresponding problem is considered for fibred-cusp
pseudodifferential operators. The full case is discussed below, initially
the discussion is restricted to the cusp algebra; this indeed is the
special case which is most closely related to the index theorem of Atiyah,
Patodi and Singer. This relationship itself is made precise below following
the discussion of the pseudodifferential index theorem.

The algebra of cusp pseudodifferential operators $\cuP m{\famfib}$ on a
compact manifold with boundary has properties similar to those of the usual
algebra on a compact manifold without boundary $\famfib;$ it is discussed
extensively in \cite{fipomb}.  The definition of cusp operators depends on
the choice of a boundary defining function $x\in \CI(\famfib)$, that is, a nonegative function 
which is zero on 
$\pa\famfib$, positive everywhere else and such that $dx$ is non-zero on
$\pa\famfib$. Given such a defining function, consider the Lie algebra
of cusp vector fields. These are arbitrary smooth vector fields in the
interior which, near the boundary are of the form 
\begin{equation*}
ax^{2}\pa_{x}+\sum_{j=1}^{m} a_{j} \pa_{z_{j}},\ a, a_{j}\in \CI(\famfib)
\label{faficu.227}\end{equation*}
where $(x,z)$  are coordinates near $\pa \famfib.$ The cusp differential
operatorsform the universal enveloping algebra of this Lie algebra (as a
$\CI(\famfib)$-module). By microlocalization this algebra is extended to the
cusp pseudodifferential operators. A typical example of cusp differential
operator is obtained by considering the Laplacian associated to a
Riemannian metric $g$ which in a collar neighborhood of the boundary $\pa
\famfib\times [0,1)_{x}\subset \famfib$ takes the form 
\begin{equation*}
g=\frac{dx^{2}}{x^{4}}+h,\ h\in \CI(\pa\famfib\times [0,1)_{x}; T^{*}(\pa\famfib)
                 \otimes T^{*}(\pa\famfib)).
\label{faficu.228}\end{equation*}

The algebra of cusp pseudodifferential operators is closely related to the
b-pseudodifferential algebra but has the virtue of admitting a \ci\
functional calculus. The main difference, compared to the boundaryless case
and as far as Fredholm properties are concerned, is that as well as a
symbol map in the usual sense, taking values in functions homogeneous of
degree $m$ on the cusp cotangent bundle, $\sigma_m:\cuP
m{\famfib}\longrightarrow \syM{\cusp}m{\famfib},$ there is a
non-commutative `indicial' (or normal) homomorphism taking values in
suspended families of pseudodifferential operators on the boundary
\begin{equation}
N:\cuP m{\famfib}\longrightarrow \Psi^m_{\sus}(\pa \famfib);
\label{faficu.9}\end{equation}
the suspended algebra consists of pseudodifferential operators with a
symbolic parameter representing the dual to the normal bundle to the
boundary. As already noted, the cusp algebra is not quite naturally
associated to a compact manifold with boundary but is fixed by the choice of a
trivialization of the normal bundle to the boundary. This may be thought of
as a residue of the product-type structure in the theorem of Atiyah, Patodi
and Singer. A cusp pseudodifferential operator is Fredholm on the
(weighted) cusp Sobolev spaces if and only if \emph{both} the symbol and
the normal operator are invertible; we describe such an operator
as `fully elliptic.' The data consisting of (compatible) pairs $\sigma
_m(P)$ and $N(P)$ constitutes the \emph{joint symbol}, so an operator is
fully elliptic when its joint symbol is invertible.

Following the model of the theorem of Atiyah and Singer described above we
consider a fibration as in \eqref{faficu.6} where now the model fibre, $\famfib,$
is a compact manifold with boundary. We define analogues of the objects
described above in the boundaryless case including a `symbolic' K-group
$\cuK{\fammap}$ as the set of stable homotopy classes of compatible and
invertible joint symbols. The analytic and topological indexes are
homomorphisms from this group into the topological K-group of the base.

The definition of the analytic index is a straightforward extension of the
boundaryless case. That is, for a fully elliptic family of cusp
pseudodifferential operators $P\in\cuP m{\famrel;\sE}$ on the fibres
of the fibration there exists a perturbation by a family of smoothing
operators supported in the interior, such that the null bundle has constant
rank and then \eqref{faficu.7} again fixes an element $\inda(P)\in
\Kt(\fambas)$ which is independent of the perturbation. In fact it only
depends on the stable homotopy class of the joint symbol, within invertible
symbols, and so defines the top map in \eqref{faficu.19} in this case.

The symbol, as opposed to the joint symbol, homomorphism leads to a short
exact sequence
\begin{equation}
\xymatrix@1{\Kt(\fambas)\ar[r]&
\cuK{\fammap}\ar[r]^(.28){\sigma}&
\Kt^0(\overline{\cuT^*(\famrel)};\cuS^*(\famrel))}=\Kc(T^*(\famrel))
\label{faficu.23}\end{equation}
where $\overline{\cuT^*(\famrel)}$ is the radial compactification of the fibrewise cusp
cotangent bundle $\cuT^*(\famrel)$(which is isomorphic to $T^*(\famrel))$ .  The image group can also be interpreted as the compactly supported
`absolute' K-theory of
$\cuT^*(\famrel)$-- the more
intricate notation here emphasizes that it is `relative' to fibre infinity
but `absolute' as regards the boundary of $\famtot.$ The first map in
\eqref{faficu.23} represents the inclusion of the fully elliptic operators
of the form $\Id+\cuP{-\infty}{\famrel;E}.$ The exactness of
\eqref{faficu.23} follows from proposition 2.2 and theorem 5.2  in \cite{fipomb}. In fact it is
shown there that any family of elliptic operators
$P\in\cuP m{\famrel;\sE},$ so only assuming the invertibility of the
symbol family, $\sigma (P),$ can be perturbed by an element of
$\cuP{-\infty}{\famrel;\sE}$ to be invertible (hence of course
fully elliptic). This leads to a splitting of the sequence
\eqref{faficu.23} which we can write
\begin{equation}
\xymatrix@1{\Kt(\fambas)\ar@<.5ex>[r]^{i_*}&
\cuK{\fammap}\ar@<.5ex>[r]^(.4){\sigma}\ar@<.5ex>[l]^{\inda}&
\Kc(T^*(\famrel))\ar@<.5ex>[l]^(.6){\inv}.}
\label{faficu.24}\end{equation}
In this sense the index element of $\Kt(\fambas)$ is a `difference class'
measuring the twisting of the given fully elliptic family relative to the
invertible perturbation; since the index map from homotopy classes of fully
elliptic families with a fixed elliptic symbol to $\Kt(\fambas)$ is an
isomorphism, the K-theory of the base can be represented by the fully
elliptic perturbations of any one elliptic family. Were it easy to
determine the map $\inv,$ \ie to find an invertible family corresponding to
a given symbol, the index problem would be much simpler!

As in the boundaryless case, this analytic index map can also be defined via
bivariant K-theory. Let $\CO_{\cusp}(\famtot)\subset\CO(\famtot)$ be
the $C^*$ subalgebra of the continuous functions on $\famtot$ consisting of those
which are constant on the boundary of each fibre of $\fammap\,;$ thus there
is a short exact sequence of $C^*$ algebras
\begin{equation}
\xymatrix@1{\COo(\famtot)\ar[r]&
\CO_{\cusp}(\famtot)\ar[r]^{R}&\CO(\fambas)}
\label{faficu.31}\end{equation}
where $\COo(\famtot)$ is the $C^*$ algebra of continuous functions on $\famtot$
vanishing on the boundary and $R$ is restriction to the boundary. The
quantization of an invertible joint symbol to a Fredholm cusp
pseudodifferential operator (and choice of parametrix) then gives an
alternative definition of the analytic index as in \eqref{faficu.29},
\eqref{faficu.30}:
\begin{equation}
\xymatrix{\cuK{\fammap}\ar[rr]^{\inda}\ar[dr]_{\quan}&&
\Kf{\fambas}{\CO(\fambas),\CO(\fambas)}=\Kt(\fambas).\\
&\Kf{\fambas}{\CO_{\cusp}(\famtot),\CO(\fambas)}\ar[ur]&}
\label{faficu.32}\end{equation}
The quantization map here is an isomorphism which we can interpret as
Poincar\'e duality. In fact we show that there is a commutative diagram
\begin{equation}
\xymatrix{\Kt(\fambas)\ar[r]\ar@{<->}[d]&
\cuK{\fammap}\ar[r]^(.4){\sigma}\ar@{<->}[d]&
\Kc(T^*(\famrel))\ar@{<->}[d]
\\
\Kf{\fambas}{\CO(\fambas),\CO(\fambas)}\ar[r]&
\Kf{\fambas}{\CO_{\cusp}(\famtot),\CO(\fambas)}\ar[r]& 
\Kf{\fambas}{\COo(\famtot),\CO(\fambas)}.}
\label{faficu.33}\end{equation}
Here the left isomorphism is the natural identification of the KK group
with K-theory and the right isomorphism is an absolute version of Atiyah's
isomorphism as discussed in \cite{Melrose-Piazza1}; it follows that the
central map is also an isomorphism.

To complete the analogy with the Atiyah-Singer theorem we proceed to define
a `topological' index using embeddings of the fibration and then prove the
equality of the two index maps. In this case, because of the
non-commutative structure of the definition of elements of $\cuK{\fammap},$
there is more of an analytic flavour to the definition of this `topological
index' than in the boundaryless case.

The main step in defining $\indt$ is to show that for an embedding of
$\fammap:\famtot\longrightarrow \fambas$ as a subfibration of
$\efammap:\efamtot\longrightarrow \fambas$ there is a natural lifting map
\begin{equation}
\efamext:\cuK{\fammap}\longrightarrow \cuK{\efammap}
\label{faficu.15}\end{equation}
corresponding to tensoring with the `Bott element' for the normal bundle of
the smaller fibration. Our construction is closely related to that of
Atiyah and Singer in \cite{Atiyah-Singer1} in the boundaryless case, but is of
necessity more intricate since we need to preserve the invertibility
of the indicial family, not just symbolic ellipticity. For this reason we
replace the construction in \cite{Atiyah-Singer1} by a slightly different one
involving pseudodifferential operators `of product type.' In view of the
identification above, this can also be considered as a lifting construction
for KK theory, giving a map, which we believe but do not show, is an
           explicit realization of the map dual to the
restriction $R:\CO_{\cusp}(\efamtot)\longrightarrow
\CO_{\cusp}(\famtot)$
\begin{equation}
R^*:\Kf{\fambas}{\CO_{\cusp}(\famtot),\CO(\fambas)}
\longrightarrow \Kf{\efambas}{\CO_{\cusp}(\efamtot),\CO(\efambas)}.
\label{faficu.34}\end{equation}
It is important for the subsequent computation of the Chern character that
this map is given by an explicit smooth construction with the corresponding
$\cuK{\fammap}.$ For the proof of the index theorem it is essential that
\eqref{faficu.15} be consistent both with the index and with the symbolic
lifting construction of \cite{Atiyah-Singer1}, which is to say that it
leads to two commutative diagrams (one for the left-directed and one for the
right-directed arrows)
\begin{equation}
\xymatrix{&
\cuK{\efammap}\ar@<.5ex>[r]^(.4){\sigma}\ar@<.5ex>[dl]^{\inda}&
\Kc(T^*(\efamrel))\ar@<.5ex>[l]^(.6){\inv}\\
\Kt^0(\fambas)\ar@<.5ex>[ur]\ar@<.5ex>[dr]\\
&
\cuK{\fammap}\ar@<.5ex>[r]^(.4){\sigma}\ar[uu]_{\efamext}
\ar@<.5ex>[ul]^{\inda}&
\Kc(T^*(\famrel))\ar[uu]_{\efamext}\ar@<.5ex>[l]^(.6){\inv}}
\label{faficu.35}\end{equation}

There is always an embedding of the fibration of compact manifolds with
boundary as a subfibration of the product
$\pi_N:\fambas\times\bbB^N\longrightarrow \fambas$ of the base with a ball
$\bbB^{N}$
of sufficiently large dimension $N$ and \eqref{faficu.15} allows
$\cuK{\fammap}$ to be mapped into the group for such a product
\begin{equation}
\cuK{\fammap}\longrightarrow\cuK{\pi_N},\ N\text{ large.}
\label{faficu.16}\end{equation}

The K-groups for these products may be computed directly (see section~\ref{Products})
and for odd
dimensions there is a `Thom isomorphism' 
\begin{equation*}
\Thom:\cuK{\pi_{2N+1}}\longleftrightarrow \Kt(\fambas).
\label{faficu.39}\end{equation*}

Hence  there is a well-defined `topological index map'
\begin{equation}
\indt=\Thom\circ(\pi_{2N+1}/\phi)^!:\cuK{\phi}\longrightarrow \Kt(\fambas)
\label{faficu.17}\end{equation}
and this completes the constructions of the ingredients in the statement of
Theorem~\ref{faficu.18} in the cusp case.

The more general case of \emph{fibred cusp} operators is similar to, but a
little more complicated than, the cusp operators discussed above. These
algebras of operators on a compact manifold correspond to the choice of a
fibration of the boundary and the behaviour of the boundary becomes more
`commutative' as the fibres become smaller -- that is, the cusp case is the
most non-commutative. Thus $\famfib,$ the typical fibre, is a compact
manifold with boundary, $\pa\famfib,$ which carries a distinguished fibration
\begin{equation}
      \xymatrix{&&\famfib\\
\boufib\ar@{-}[r]&\pa\famfib\ar@{^{(}->}[ur]\ar[d]^{\boufibmap}\\&\boufibbas}
\label{faficu.1}\end{equation}
where $\boufibbas$ and $\boufib$ are compact manifolds without
boundary. Associated with this structure, and a choice of trivialization,
along the fibres of $\boufibmap,$ of the normal bundle to the boundary, is
an algebra of `fibred cusp pseudodifferential operators'
$\fcP{\boufibmap}{*}{\famfib}$ introduced in \cite{Mazzeo-Melrose4}. The
cusp case corresponds to the extreme case of 
the one-fibre fibration $\boufibbas=$\{\text{pt}\}. The other extreme case,
choosing the point fibration of the boundary, $\boufibbas=\pa\famfib,$
$\boumap=\Id,$ is the `scattering case' in which the index theorem is
reducible directly to the usual Atiyah-Singer setting by a doubling
construction (briefly discussed in \cite{MelroseGST} and below); this case is
the simplest in most senses. It is discussed separately below, since a
special case of the index theorem for scattering operators is used to
handle the families index for perturbations of the identity in the general
case. It is perhaps helpful to think of the fibred cusp operators as
associated to the topological space $\famfib/\boufibmap$ in which the
boundary is smashed to the base $\boufibbas.$ Thus, for the cusp calculus,
the boundary is smashed to a point whereas for the scattering calculus it is
left unchanged.

Here we consider locally trivial families of such structures. Thus, suppose that
$\famtot$ is a manifold with boundary which admits a fibration over a base
$\fambas$ with typical fibre $\famfib$ which can be identified with the
manifold in \eqref{faficu.1}. We suppose that the boundary of $\famtot$ has
a finer fibration than over the base $\fambas$ giving a commutative diagram
\begin{equation}
\xymatrix{&&\famfib\ar@{-}[r]&\famtot\ar[d]^{\fammap}\\
\boufib \ar@{-}[r]&\pa\famfib\ar@{^(->}[ur]\ar[d]^{\boufibmap}
\ar@{-}[r]&\pa\famtot\ar[r]^{\pa\fammap}\ar[d]^{\boumap}
\ar@{^(->}[ur]&\fambas\\
&\boufibbas\ar@{-}[r]&\ar[ur]\boubas
}
\label{faficu.2}\end{equation}
in which $\boufib$ is the typical fibre of a fibration of $\pa\famtot$ over
$\boubas$ and $\boufibbas$ is the typical fibre of a fibration of $\boubas$ over
$\fambas.$ We call such a pair of fibrations, together with a choice of
normal trivialization, a fibration with fibred cusp structure. Associated
to this geometry is an algebra of fibred cusp pseudodifferential operators
acting on the fibres of $\famtot$ over $\fambas,$ which we shall denote
$\fcP{\boumap}*{\famrel}.$ An element is a family parameterized by
$\fambas$ with each operator related to the fibre (in $\famtot)$ above that
point with its boundary smashed to the fibre of $\boubas$ above that point
of $\fambas$ by the fibre of $\boumap.$

In this more general setting we obtain similar results to those described
above for the cusp calculus. First we define an abelian group,
$\fcK{\boumap}(\fammap),$ with elements which are the stable homotopy
classes of joint elliptic symbols, and a corresponding odd K-group,
$\fcKo{\boumap}(\fammap).$ The definition of analytic index maps, taking
values respectively in $\Kt(\fambas)$ and $\Kto(\fambas)$ is essentially as
above. We also give an analogue of \eqref{faficu.23} but now as a 6-term
exact sequence
\begin{equation}
\xymatrix{\Kc(T^*(\boubasrel))\ar[r]^(0.6){i_{0}}&
\fcK{\boumap}(\fammap)\ar[r]^(.4){\sigma_{0}}&
\Kc(T^*(\famrel))\ar[d]^{I_{0}}\\
\Kco(T^*(\famrel))\ar[u]^{I_{1}}&
\fcKo{\boumap}(\fammap)\ar[l]^(.40){\sigma_{1}}&
\Kco(T^*(\boubasrel))\ar[l]^(.55){i_{1}}}
\label{faficu.36}\end{equation}
where $\sigma_{0},$ $\sigma_{1}$ are maps related to the symbol, $I_{0},$ 
$I_{1}$ are
forms of the index map of Atiyah-Singer and $i_{0},$ $i_{1}$ correspond to
inclusion of perturbations of the identity by operators of order $-\infty.$ 
We show below that this is isomorphic to the corresponding sequence in
KK-theory arising from the short exact sequence of $C^*$ algebras
(replacing \eqref{faficu.31})
\begin{equation}
\begin{gathered}
\COo(\famtot)\longrightarrow\CO_{\boumap}(\famtot)\longrightarrow \CO(\boubas),\\
\CO_{\boumap}(\famtot)=\left\{f\in\CO(\famtot);
f\big|_{\pa\famtot}=\boumap^*g,\ g\in\CO(\boubas)\right\},
\end{gathered}
\label{faficu.37}\end{equation}
namely 
\begin{equation}
\xymatrix{\Kf{\fambas}{\CO(\boubas),\CO(\fambas)}\ar[r]&
\Kf{\fambas}{\CO_{\boumap}(\famtot),\CO(\fambas)}\ar[r]&
\Kf{\fambas}{\COo(\famtot),\CO(\fambas)}\ar[d]\\
\Kfo{\fambas}{\COo(\famtot),\CO(\fambas)}\ar[u]&
\Kfo{\fambas}{\CO_{\Phi}(\famtot),\CO(\fambas)}\ar[l]&
\Kfo{\fambas}{\CO(\boubas),\CO(\fambas)}\ar[l].}
\label{faficu.38}\end{equation}
where the four isomorphisms between the corresponding spaces on the left
and right in \eqref{faficu.38} and \eqref{faficu.36} are the Poincar\'e
duality maps of Atiyah, as realized in Kasparov's KK theory. The remaining
two isomorphisms are given by quantization maps in the fibred-cusp calculus
amounting as before to Poincar\'e duality; this is Theorem~\ref{faficu.88} in the
general case. Note that the fact that the symbol map in \eqref{faficu.36}
is not (in general) surjective shows that in some sense the `universal
case' is that of the cusp calculus since only through it can \emph{every}
elliptic symbol be quantized to a Fredholm family.

This universality appears explicitly in the form of a natural map in which
the finer fibration of the boundary is `forgotten' 
\begin{equation}
\xymatrix{\fcK{\boumap}(\fammap)\ar[rr]^{q_{\ad}}\ar[dr]_{\inda}&&
\cuK{\fammap}\ar[dl]^{\inda}\\
&\Kt(\fambas)}
\label{faficu.220}\end{equation}
through which the index factors; it is defined through an adiabatic
limit. We then define the topological index as the composite with the
topological index in the cusp case. It is also possible to proceed more
directly, through an appropriate embedding of the fibration.

Since we make substantial use below of various classes (`calculi') of
pseudodifferential operators we have tried to use a uniform notation. In
the families index theorem of Atiyah and Singer the quantization map giving
the analytic index is in terms of families of pseudodifferential operators
on the fibres of a fibration. We use what is the standard notation for
these, $\Ps m{\famrel;\sE},$ except that $\sE=(E_+,E_-)$ is a
$\bbZ_2$-graded bundle and the operators act from $E_+$ to $E_-.$ We will 
sometime need to consider tensor products of the form
\[
             \sE\otimes_{+}F= (E_{+}\otimes F, E_{-})
\]
where $F$ is a bundle which is not $\bbZ_{2}$-graded.
We also use $\bbZ_{2}$-graded bundle notation for the operators in the
quantization map \eqref{faficu.62}, on a fibration with fibred cusp
structure 
\begin{equation*}
\fcP{\boumap}m{\famrel;\sE}
\label{faficu.196}\end{equation*}
where the suffix alone indicates the `uniformity type' at the
boundary. Several more such classes arise below. First the indicial
operator for families as in \eqref{faficu.9} takes values in the (singly)
suspended algebra for which we use the suffix `sus.' These are the
pseudodifferential operators on $\bbR\times\famfib,$ or for families on the
fibres of $\bbR\times\pa\famtot,$ which are translation-invariant and
rapidly decreasing at infinity in the real factor -- taking the Fourier
transform therefore gives a parameterized family but the parameter enters
as a symbolic variable. So more generally we use notation such as
\begin{equation*}
\fsP{\boufibmap}m{\boufib},\ \ptP{\boufibmap}{m',m}{\boufib},\
\adP{\boufibmap}m{\boufib} 
\label{faficu.197}\end{equation*}
to denote, respectively, pseudodifferential operators which are suspended
with respect to the cotangent variables of a fibration, pseudodifferential
operators which are of product type with respect to a fibration and
pseudodifferential operators which depend adiabatically on a parameter (and
pass from operators on the total space to fibre operators in the limit in
the parameter). The general approach to pseudodifferential operators, by
defining them in terms of classical conormal distributions on some blown up
version of the product space, allows these types to be combined where
required. Thus we use adiabatic families of fibred cusp operators to pass
from one boundary fibration to another in the definition of
\eqref{faficu.220} and product-type cusp pseudodifferential operators in
the lifting construction -- these have product-type suspended operators as
normal operators with corresponding notation.

We do not put the $0$'s at the end of short exact sequences.

In \S\ref{Analyticindex} the analytic index, in K-theory, of fully elliptic
fibred cusp pseudodifferential operators is described as a map from the
K-group associated to invertible full symbols. In the special case of the
scattering structure, the index theorem is derived from the Atiyah-Singer
index theorem in \S\ref{Scattering} and used to discuss families which are
perturbations of the identity for general fibred cusp structures in
\S\ref{Minusinfinity}. The quantization map, to KK-theory, is introduced
in \S\ref{KK.0} and shown to be an isomorphism for the cusp structure in
\S\ref{aco.0}. The 6-term symbol sequences in K-theory and KK-theory are
described in \S\ref{stes.0} and related in \S\ref{GeneralPoincare},
resulting in the proof of Theorem~\ref{faficu.88}. The map
\eqref{faficu.220} is constructed in \S\ref{Adiabaticpassage} and a result on the
extension of fibred cusp structures is contained in \S\ref{efc.0} and the
related multiplicativity and lifting properties are examined in
\S\S\ref{Multiplicativity},\ref{Lifting}. In \S\ref{Products} the model
cases are analyzed and used to define the topological index map and to prove
Theorem~\ref{faficu.18} in \S\ref{ti.0}. The application of the cusp case
to families of Atiyah-Patodi-Singer type is made in \S\ref{APS} and
Theorem~\ref{faficu.153} is proved there. In the appendices the various classes
of pseudodifferential operators appearing in the body of the paper are
described, including product-type, fibred cusp and adiabatic algebras and
their combinations.

\paperbody

\section{Analytic index}\label{Analyticindex}

As already briefly described above, the general setting of this paper is a
compact manifold with boundary, $\famtot,$ with a fibration
$\fammap:\famtot\longrightarrow \fambas$ over a compact manifold usually
without boundary; we also assume, without loss of generality, that the base
is connected. In fact it is convenient at various points to allow the base
(and hence also the total space) to be a manifold with corners, to allow
especially products with intervals. We still treat $\famtot$ as a manifold
with boundary, in that `the' boundary is then the union of the boundary
faces of the fibres. Thus, $\fammap$ is a smooth surjective map with
surjective differential at every point. It follows that the fibres,
$\fammap^{-1}(b)$ for $b\in\fambas,$ are compact manifolds with boundary,
all diffeomorphic to a fixed manifold $\famfib$ for which we use the
notation \eqref{faficu.6}. It is often notationally convenient to assume
that the fibres are also connected, but it is by no means necessary and we
believe that the paper is written so that this assumption is not used. Such
a fibration is always locally trivial.

In addition we assume that the boundary of the total space of the
fibration, $\pa\famtot,$ carries a second, finer, fibration, giving a
commutative diagram of fibrations 
\begin{equation}
\xymatrix{\pa\famtot\ar[rd]_{\boumap}\ar[rr]^{\pa\fammap}&&\fambas.\\
&\boubas\ar[ur]_{\pa\fammap/\boumap}}
\label{faficu.66}\end{equation}
Thus the boundary of each fibre carries a fibration and the overall fibre, with
this fibration of its boundary, is always diffeomorphic to the model fibre
with its model boundary fibration, as in \eqref{faficu.1}. The maps fit
together as in \eqref{faficu.2} and it is always the case that there is a local
trivialization of the overall fibration $\fammap$ in which the fibration of
the boundary of the fibres is reduced to this normal form; in this sense
the structure is locally trivial.

To associate with $\fammap$ and $\boumap$ a class of pseudodifferential
operators on the fibres of $\fammap,$ of the type introduced in
\cite{Mazzeo-Melrose4}, we need one more piece of information. Namely we need to
choose a trivialization of the normal bundle to the boundary of the fibres
of $\fammap$ along the fibres of $\boumap.$ This simply amounts to the
choice of a boundary defining function $x\in\CI(\famtot)$ (so $x\ge0,$
$\pa\famtot=\{x=0\}$ and $dx\not=0$ on $\pa\famtot).$ Such a choice gives a
trivialization of the normal bundle to $\famtot$ everywhere (simply choose
the inward pointing normal vector $v_p$ at $p\in\pa\famtot$ to satisfy
$v_px=1).$ Two such choices $x,$ $x',$ are equivalent, in the sense that
they give the same trivialization along the fibres of $\boumap$ if and only
if 
\begin{equation}
x'=ax+bx^2,\ a,b\in\CI(\famtot),\ a\big|_{\pa\famtot}=\boumap^*a',\
a'\in\CI(\boubas).
\label{faficu.67}\end{equation}
Equivalent choices turn out to lead to the same algebra of pseudodifferential
operators. Since even inequivalent choices are homotopic and lead to
isomorphic structures the dependence on this choice of normal
trivialization will not be emphasized; it should be fixed throughout but
none of the results depend on which choice is made.

\begin{definition}\label{faficu.198} A \emph{fibration with fibred cusp
structure} is a fibration of compact manifolds \eqref{faficu.6} with fibres
modelled on a fixed compact manifold with boundary, a finer fibration as in
\eqref{faficu.66} of the boundary of the total space and a choice of
trivialization of the normal bundle to the boundary along the fibres as
discussed above.
\end{definition}

Notice that the extreme cases in which $\boumap=\pa\fammap$ is simply the
restriction of the fibration to the boundary of $\famtot$ is a particularly
interesting case, the `cusp' case, which includes the setting of the index
theorem of Atiyah, Patodi and Singer. The other extreme, in which $\boumap$
is the identity, is the essentially commutative `scattering' case.  In
terms of \cite{Atiyah-Bott}, this corresponds to operators for which the
obstruction in K-theory to the existence of local elliptic boundary
conditions vanishes (in contrast with the global nature of the
Atiyah-Patodi-Singer boundary condition).  The general fibred cusp case is
intermediate between these two extremes.

In general, even if we were assuming that $\famtot$ is connected, it would
be artificial to assume that the boundary $\pa\famtot$ is connected, since
many interesting examples involve a disconnected boundary. We therefore
avoid any such assumption. When the boundary is not connected, the cusp
case is still interpreted as $\boumap=\pa\fammap.$ As discussed in
\cite{fipomb} and \cite{rochon}, this means one should allow terms of order
$-\infty$ acting between different components of the boundary and
correspondingly in the definition of the indicial family.

Under these conditions we may associate to a fibration with fibred cusp
structure an algebra of `fibred-cusp' pseudodifferential operators,
$\fcP{\boumap}{\bbZ}{\famrel}.$ The reader is referred to
\cite{Mazzeo-Melrose4} for the original definition of the algebra and to
the discussion in \cite{rochon}. A generalization to product-type
operators, used in the lifting construction below, is given in
Appendix~\ref{PTFCO}. For our purposes here the main interest lies in the
boundedness, compactness, Fredholm and related symbolic properties of these
operators. Thus if $\sE=(E_+,E_-)$ is a $\bbZ_2$-graded complex vector bundle
we will denote by $\fcP{\boumap}{\bbZ}{\famrel;\sE}$ the space of these
operators acting from sections of $E_+$ to sections of $E_-;$ they always
give continuous linear operators on weighted spaces of smooth sections
\begin{equation}
P\in\fcP{\boumap}{\bbZ}{\famrel;\sE}:x^s\CI(\famtot;E_+)\longrightarrow x^s\CI(\famtot;E_-),\
\forall\ s\in\bbR.
\label{faficu.68}\end{equation}
Both the symbol and the normal operator can be defined by appropriate
`oscillatory testing' of these maps. The symbol map gives a short exact sequence 
\begin{equation}
\xymatrix@1{\fcP{\boumap}{m-1}{\famrel;\sE}\ar[r]&
\fcP{\boumap}{m}{\famrel;\sE}\ar[r]^(.5){\sigma}&
\fcS{\boumap}{m}(\famrel; \sE) }
\label{faficu.69}\end{equation}
with 
\begin{equation}
\fcS{\boumap}{m}(\famrel;\sE)=\CI(\fcs^*(\famrel);\hom(\sE)\otimes N_{-m})
\label{symbol}\end{equation}
where $\fcs^*(\famrel)$ is the sphere bundle of the fibred cusp
cotangent bundle of the fibres (isomorphic, but not naturally so, to the
`usual' bundle $T^*(\famrel))$ and $N_{-m}=N^{-m}$ is the bundle with
sections which are the functions homogeneous of degree $m$ on the fibres
(then $N$ can be identified with the normal bundle to the boundary of the
radial compactification) and $\hom(\sE)=\hom(E_+,E_-).$ The normal (or
indicial) operator gives a short exact sequence 
\begin{equation}
\xymatrix@1{x\fcP{\boumap}{m}{\famrel;\sE}\ar[r]&
\fcP{\boumap}{m}{\famrel;\sE}\ar[r]^{N}&
\fsP{\boumap}m{\bourel;\sE}}.
\label{faficu.70}\end{equation}
Here the image is the space of suspended pseudodifferential operators on
the fibres of the fibration $\boumap$ with symbolic parameters in 
\begin{equation}
\fcT{\boumap}^*_{\pa\famtot}(\famrel)/T^*(\pa\famrel)
\label{faficu.71}\end{equation}
which is to say the duals to the fibre variables of $\boumap$ together with
a dual to the normal variable. The symbol and normal operators are
connected by the identity
\begin{equation}
\sigma\big|_{\pa\famtot}=\sigma\circ N\text{ on }\fcP{\boumap}{m}{\famrel;\sE}.
\label{faficu.72}\end{equation}
This is the only constraint, so if we set 
\begin{multline}
\jS{\boumap}m(\fammap;\sE)=\big\{(\sigma ,N)\in\\
   \fcS{\boumap}{m}(\famrel;\sE)\oplus\fsP{\boumap}m{\bourel;\sE};
\sigma \big|_{\pa\famtot}=\sigma _m(N)\big\}
\label{faficu.74}\end{multline}
the joint symbol sequence becomes 
\begin{equation}
\xymatrix@1{x\fcP{\boumap}{m-1}{\famrel;\sE}\ar[r]&
\fcP{\boumap}{m}{\famrel;\sE}\ar[r]^(.58){(j,N)}&
\jS{\boumap}m(\fammap;\sE)}.
\label{faficu.73}\end{equation}

Such a family extends by continuity to an operator between the natural
continuous families of weighted fibred-cusp Sobolev spaces
\begin{equation}
P:\CO(\fambas;x^sH^{M}_{\fC{\Phi}}(\famrel;E_+))\longrightarrow
\CO(\fambas;x^sH^{M-m}_{\fC{\Phi}}(\famrel;E_-)),\ M,s\in\bbR.
\label{faficu.40}\end{equation}
It is a Fredholm family on these spaces if and only if its image in 
$\jS{\boumap}m(\fammap;\sE)$ is invertible, in which case we say that the
family is fully elliptic.

\begin{lemma}\label{faficu.41} If $P\in\fcP{\boumap}{\bbZ}{\famrel;\sE}$ is a
  fully elliptic family then there exists a smoothing perturbation in
  $x^\infty\fcP{\boumap}{-\infty}{\famrel;\sE}$ such that $P+A$ has null
  space in \eqref{faficu.40} forming a trivial smooth vector bundle, over
  $\fambas,$ contained in $\dCI(\famtot;E_+).$
\end{lemma}

\begin{proof} The algebra is invariant under conjugation by $x^s$ and such
  conjugation affects neither the symbol nor the normal operator; thus we
  can take $s=0.$ For simplicity of notation we will also suppose that $M=m$ in
  \eqref{faficu.40}; this is all that is used below and the general case
  follows easily. The full ellipticity and properties of the calculus allow
  us to construct a parametrix, which is to say a family
  $Q\in\fcP{\boumap}{-m}{\famrel;\sE^-},$ where $\sE^-=(E_-,E_+),$ such that
\begin{equation}
PQ=\Id_+-R_+,\ QP=\Id_--R_-,\ R_\pm\in
  x^{\infty}\fcP{\boumap}{-\infty}{\famrel;E_{\pm}}.
\label{faficu.42}\end{equation}
As a bundle of Hilbert spaces,
$\CO(\fambas;H^m_{\fC{\boumap}}(\famrel;E_+))$ is necessarily
trivial. Thus there is a sequence of orthogonal projections $\Pi_N,$
tending strongly to the identity, corresponding to projection onto the
first $N$ terms of an orthonormal basis in some trivialization. In
particular for a compact operator such as $R_+,$ $R_+\Pi_N\longrightarrow
R_+$ in the topology of norm continuous families of bounded operators.

Furthermore, $\Pi_N$ may be approximated by smooth families of projections
of the same rank. To see this, fixing $N,$ we may certainly find a sequence of
smoothing operators with supports disjoint from the boundary such that
$W_j\longrightarrow \Pi_N$ in the norm topology; it suffices to
use a partition of unity and work locally. Then $W_j'=W_j\Pi_N W_j$ is a
continuous (in $\fambas)$ family of smoothing operators, supported in the
interior, which approximates $\Pi_N$ and for large $j$ has rank $N.$ Taking
$j$ sufficiently large the range of $W_j'$ is a trivial subbundle which has
a basis $e_{l,j}\in\CO(\fambas;\CI(\famrel;E_+))$ with supports always in
the interior. These sections can themselves be uniformly approximated by
sections of the same form but smooth over $\fambas.$ Replacing $W_j'$ by
the orthogonal projection onto the span of these sections, for $j$ large,
we have succeeded in approximating $\Pi_N$ by smooth families of
projections with kernels supported in the interior. Thus we may suppose
that the $\Pi_N$ are smooth families of smoothing operators with kernels
vanishing to infinite order at both boundaries of the product.

Returning to \eqref{faficu.42} we may compose on the right with $\Id-\Pi_N$
and, using the fact that it is also a projection, deduce that 
\begin{equation}
QP(\Id-\Pi_N)=(\Id-B)(\Id-\Pi_N),\ B=R_+(\Id-\Pi_N).
\label{faficu.43}\end{equation}
For large $N$ it follows that $B$ has uniformly small norm. The inverse of
$\Id-B$ is then also of the form $\Id+B'$ with $B'$ a smoothing operator with
kernel vanishing to infinite order at the boundary, and so is an element of
$x^\infty\fcP{\boumap}{-\infty}{\famrel;E_+}.$ Replacing $Q$ by the new
parametrix $(\Id+B')Q$ we have replaced the first identity in
\eqref{faficu.42} by
\begin{equation}
Q(P+A)=\Id-\Pi_N,\ A=-P\Pi_N\in x^\infty\fcP{\boumap}{-\infty}{\famrel;\sE}.
\label{faficu.44}\end{equation}
This perturbs $P$ as desired.
\end{proof}

Now once the null space of $P+A$ is arranged to be a smooth bundle it
follows that its range has a complement of the same type 
(but in general not trivial), namely the null
bundle of $P^*+A^*$ for some choice of inner products and smooth
density.

\begin{proposition}\label{faficu.45} The element 
\begin{equation*}
\nul(P+A)\ominus\nul(P^*+A^*)\in\Kt(\fambas)
\label{faficu.75}\end{equation*}
defined by any fully elliptic element of $\fcP{\boumap}m{\famrel;E}$ using
Lemma~\ref{faficu.41} is independent of the choice of perturbation $A,$ 
is constant under smooth homotopy and is additive under direct sums.
\end{proposition}

\begin{proof} First suppose $P$ has been perturbed so that
  \eqref{faficu.44} holds for some parametrix $Q$ and some family of smooth
  projections $\Pi_N.$ Changing to another family $\Pi'_k$ we can choose
  $k$ so large that $B'=\Pi_N(\Id-\Pi'_k),$ has small norm and then 
\begin{equation}
(\Id-B')^{-1}Q(P+A')=\Id-\Pi'_k,\ A'=-P\Pi'_k+A(\Id-\Pi'_k).
\label{faficu.76}\end{equation}
The null space of the new family is the range of $\Pi'_k$ into which the
  range of $\Pi_N$ is mapped isomorphically by $\Id-\Pi'_k.$ The range of
  $P+A'$ is then the direct sum of an isomorphic image of the complement of
  this image of $\Pi_N$ plus the previous range. Thus the element in
  $\Kt(\fambas)$ is unchanged.

To see the independence of the choice of stabilizing perturbation, suppose
$A$ and $A'$ are two such perturbations. Consider the family depending on
an additional parameter $P+(\cos\theta A+\sin\theta A'),$ $\theta \in\bbS.$
The circle can be included in the base of the product fibration and then
the argument above can be applied to stabilize the new family. This shows
that the pairs of bundles resulting from different stabilizations are
homotopic and so define the same element in $\Kt(\fambas).$ Indeed the same
argument applies to a homotopy of the operator itself, through fully elliptic
operators.

It is then immediate that the index class is additive under direct sums.
\end{proof}

\begin{remark}\label{faficu.221} Various `stabilization' constructions like
  this are used below. For instance, if $P_t$ is a family of totally elliptic
  operators depending smoothly on an additional parameter $t\in[0,1]$ and
  it is invertible for $t=0$ then there is a finite rank smoothing
  perturbation, vanishing at $t=0,$ which makes the family invertible for
  all $t\in[0,1].$
\end{remark}

To formulate the index as a map we now consider the K-group which arises from
the symbol algebra.

\begin{definition}\label{faficu.77} For a fibration with fibred cusp
structure, $\fcK{\boumap}(\fammap)$ denotes the set of
equivalence classes of the collection of the invertible elements of
the $\jS{\boumap}0(\fammap,\sE)$ (with inverse in $\jS{\boumap}0(\fammap,\sE^-)$,
$\sE^-= (E^{-},E^{+})$)
where the equivalence relation is a finite chain consisting of the following 
\begin{gather}
\begin{gathered}
(\sigma, N)\in\jS{\boumap}0(\fammap;\sE)\sim(\sigma _1,N_1)
\in\jS{\boumap}0(\fammap;\sF)\\
\text{ if there exist bundle isomorphisms over
}\famtot,\ a_\pm:E_\pm\longrightarrow F_\pm\\
\text{ such that }\sigma
=a_-\circ\sigma _1\circ a_+\Mand
N=(a_-\big|_{\pa\famtot})\circ N _1\circ (a_+\big|_{\pa\famtot})
\end{gathered}
\label{faficu.199}
\\
\begin{gathered}
(\sigma, N)\in\jS{\boumap}0(\fammap;\sE)\sim(\sigma _2,N_2)
\in\jS{\boumap}0(\fammap;\sE)\\
\text{ if there exists a homotopy }(\sigma(t),N(t)),\
t\in[0,1],\ (\sigma (t)^{-1},N(t)^{-1})\in\jS{\boumap}0(\fammap;\sE^-)\\
\text{such that }\sigma(0)=\sigma,\ \sigma(1)=\sigma_2,\ N(0)=N,\ N(1)=N_2,\
\end{gathered} 
\label{faficu.200}
\\
\begin{gathered}
(\sigma, N)\in\jS{\boumap}0(\fammap;\sE)\sim(\sigma _3,N_3)
\in\jS{\boumap}0(\fammap;\sE\oplus F)\\
\Mif F\text{ is ungraded},\ \sigma _3=\sigma \oplus\Id_F\Mand N_3=N\oplus\Id_F.
\end{gathered}
\label{faficu.78}\end{gather}
\end{definition}

There is also a corresponding odd K-group.  Given the diagram 
\eqref{faficu.2}, let $I=[0,1]$ be the unit interval
and consider the suspended version 
\begin{equation}
\xymatrix{&&\famfib\ar@{-}[r]&\famtot\times I\ar[d]^{s\fammap}\\
\boufib \ar@{-}[r]&\pa\famfib\ar@{^(->}[ur]\ar[d]^{\boufibmap}
\ar@{-}[r]&\pa\famtot\times I\ar[r]^{\pa s\fammap}\ar[d]^{s\boumap}
\ar@{^(->}[ur]&\fambas\times I\\
&\boufibbas\ar@{-}[r]&\ar[ur]\boubas\times I }
\label{sK.0}\end{equation} 
where $s\fammap=\fammap\times \Id$ and $s\boumap=\boumap\times\Id.$
Let $E$ be an ungraded complex vector bundle. To the fibration \eqref{sK.0}
corresponds the space of joint symbols $\jS{s\boumap}m(s\fammap;E).$

\begin{definition}\label{faficu.77a} For a fibration with fibred cusp
structure, $\fcKo{\boumap}(\fammap)$ denotes the set of equivalence classes
of the collection of the invertible elements of the
$\jS{s\boumap}0(s\fammap,E)$ (with $E$ ungraded and with inverse in
$\jS{s\boumap}0(s\fammap,E)$) which are the identity when restricted to
$B\times\{0,1\}$ under, the equivalence relation corresponding to a finite
chain as in \eqref{faficu.199}, \eqref{faficu.200} and \eqref{faficu.78}
with bundle transformations and homotopies required to be the identity on
$B\times\{0,1\}.$
\end{definition}

\begin{proposition}\label{faficu.79} For any fibration with fibred cusp
structure, both $\fcK{\boumap}(\fammap)$ and $\fcKo{\boumap}(\fammap)$ are
abelian groups  under direct sum (or equivalently stabilized product) which
are naturally independent of the choice of boundary trivialization and the index
construction of Lemma~\ref{faficu.41} defines group homomorphisms
\begin{equation}
\inda:\fcK{\boumap}(\fammap)\longrightarrow \Kt(\fambas),\
\inda:\fcKo{\boumap}(\fammap)\longrightarrow \Kto(\fambas). 
\label{faficu.80}\end{equation}
\end{proposition}

\begin{proof} The abelian group structure follows as in the boundaryless
  case. Since changing the boundary defining function only affects the
  calculus through a change of the trivialization of the normal bundle and
  all such trivializations are homotopic the resulting abelian groups are
  independent of this choice and Proposition~\ref{faficu.45} shows that the
  analytic index map, \eqref{faficu.80}, is then well defined and additive.
\end{proof}

\section{The scattering case}\label{Scattering}

The scattering case, in which the fibres of the boundary fibration are
reduced to points is effectively `commutative' compared to the others. In
particular it is very close to the setting of the original Atiyah-Singer
index theorem and we show here that it is reducible to it. The resulting
identification of the analytic and topological indexes allows us to derive,
in the next section, an index theorem in the setting of a general boundary
fibration but for perturbations of the identity of order $-\infty.$ 

\begin{lemma}\label{faficu.128} For the scattering structure on any compact
fibration, \eqref{faficu.6},
\begin{equation}
\scK(\fammap)=\fcK{\Id}(\fammap)
\equiv\Kc(T^*(\famrel); T^*_{\pa\famtot}(\famrel)),
\label{faficu.63}\end{equation}
is identified with the compactly supported K-theory of the fibre cotangent
bundle of the interior of $\famtot.$
\end{lemma}

\begin{proof} In this case an element of $\jS{\Id}0(\fammap;\sE)$ is a pair
$(\sigma ,b)$ each of which is a bundle isomorphism, taking values in the
lift of $\hom(\sE).$ The `symbolic part' $\sigma$ is defined on the sphere
bundle at infinity of the (radial compactification of the) appropriately rescaled
fibre cotangent bundle $\scT^*(\famrel)$ and the boundary part is smooth on
the radial compactification of the restriction,
$\scT_{\pa\famtot}^*(\famrel)$ of this bundle to the boundary. Since they
are compatible at the intersection, which is to say the corner of
$\com{\scT^*}(\famrel),$ together this gives \emph{precisely} a section of
$\hom(\sE)$ lifted to the boundary of $\com{\scT^*}(\famrel).$ This is the
data needed for the standard definition of a compactly supported K-class in
the interior of a manifold with boundary (to which such a manifold with
corners is homeomorphic) and all classes arise this way. This gives
\eqref{faficu.63}.
\end{proof}

Poincar\'e duality reduces to one of the cases discussed (for
$\fambas=\{\text{pt}\})$ in \cite{Melrose-Piazza1}
\begin{equation}
\Kc(T^*(\famrel); T^*_{\pa\famtot}(\famrel))\longleftrightarrow
\Kf{\fambas}{\CO(\famtot),\CO(\fambas)}.
\label{faficu.64}\end{equation}

The index theorem in the scattering case has a `simple' formulation and
proof in the sense that it reduces directly to the Atiyah-Singer
theorem through the following observation from \cite{MelroseGST}.

\begin{lemma}\label{faficu.129} In the scattering case $\scK(\fammap)$ is
  generated by the equivalence classes of elements of the subset 
\begin{equation}
\{(\sigma ,b)\in\jS{\Id}0(\fammap;\sE);\text{ near }\pa\famtot,\
E_+=E_-=\bbC^N,\ \sigma =b=\Id\}\subset\jS{\Id}0(\fammap;\sE).
\label{faficu.130}\end{equation}
\end{lemma}

\begin{proof} First we show that any invertible joint symbol is homotopic
to an element in which, near the boundary, both $\sigma$ and $b$ are the
lifts of some bundle isomorphism from $E_+$ to $E_-.$ Given an element
$(s,b)\in\jS{\Id}0(\fammap;\sE)$ the bundle isomorphism can be taken to be
any extension off $\pa\famtot$ of the section $b$ restricted to the zero
section (identified with $\pa\famtot)$ of $\scT^*_{\pa\famtot}(\famrel).$
First perturb $b$ to be constant on the (linear) fibres of
$\scT^*_{\pa\famtot}(\famrel)$ near the zero section. An extension of the
radial expansion of the vector bundle allows it to be deformed, with
$\sigma$ to keep the consistency condition, to be fibre constant in this
sense in a neighbourhood of the boundary.

Now, given that $\sigma$ and $b$ are identified with a bundle isomorphism
near the boundary, this isomorphism can be used to modify $E_-$ to be equal
to $E_+$ close to the boundary so that both symbols become the identity on
$E_+$ very close to the boundary without changing the equivalence
class. Then $E_+$ can be complemented to a trivial bundle.
\end{proof}

Once the K-group is identified with the set of equivalence classes as in
\eqref{faficu.130}, the quantization map can also be arranged to yield
pseudodifferential operators, in the ordinary sense, which are equal to the
identity in a neighborhood of the boundary. Such operators can be extended to the
double of $\famtot,$ across the boundary, to be the identity on the
additional copy of $\famtot$ and the Atiyah-Singer index theorem then applies.

\begin{theorem}\label{faficu.201} For the scattering structure on a
  fibration the analytic index for fully elliptic scattering operators on
  the fibres factors through the Atiyah-Singer index map for the double
  $2\famtot=\famtot\cup\famtot^-$ 
\begin{multline}
\inda:\scK(\fammap)\overset{\sim}\longrightarrow
  \Kc(T^*(\famrel);T^*_{\pa\famtot}(\famrel))=\\
\Kc(T^*(2\famrel);T^*(\famtot^-/\fambas))\longrightarrow 
\Kc(T^*(2\famrel))\overset{\indAS}\longrightarrow \Kt(\fambas).
\label{faficu.202}\end{multline}
\end{theorem}

\begin{proposition}\label{faficu.131} Suppose that $E$ is a complex vector
bundle over the total space of a fibration with scattering structure and
$b\in\cS(\scT^*_{\pa\famtot}(\famrel);\hom(E))$ is such that $\Id+b$ is
invertible then any family $\Id+B,$ $B\in\Psi_{\scat}^{-\infty}
(\famrel;E)$ with
$N(\Id+B)=\Id+b$ has analytic index equal to the image  
\begin{equation}
[\Id+b]\in\Kco(\bbR\times T^*(\pa\famrel))=\Kc(T^*(\pa\famrel))\overset{\indAS}
\longrightarrow \Kt(\fambas)
\label{faficu.132}\end{equation}
under the Atiyah-Singer index map for the boundary.
\end{proposition}

\begin{proof} First we can complement $E$ to be trivial, stabilizing the
symbol by the identity. By a small perturbation we can also arrange that the
boundary symbol $b$ is compactly supported on
$\scT^*_{\pa\famtot}(\famrel)=\bbR\times T^*(\pa\famrel)$ and is equal to a
bundle map near the zero section. Thus, the varying part of $b$ can be
confined to a compact subset of $(0,\infty)\times W,$ where $W$ is
the boundary of the radial compactification of the vector bundle
$\bbR_s\times T^*(\pa\famrel)$ and the first variable is the radial
variable. In the deformation in the proof of Lemma~\ref{faficu.129} above,
across the corner at infinity and into the interior, the radial variable
becomes the normal variable to the boundary, $x,$ in a product
decomposition near the boundary with the variation now in $(0,1)$ and $W$
is identified with the boundary of the radial compactification of
$T^*(\famrel)\big|_{x=\frac12}.$ Thus, $\Id+b$ has been identified with a
symbol in the conventional sense on the sphere bundle of $T^*((0,1)_x\times
\pa\famrel)$ reducing to the identity near $x=1$ and some bundle
isomorphism near $x=0.$ By Bott periodicity of the Atiyah-Singer index map
this reduces to the identifications in \eqref{faficu.132}.
\end{proof}

\begin{remark}\label{faficu.203} From this result a families index theorem
  in K-theory in the setting of Callias' index theorem follows. See also
  the discussion by Anghel \cite{MR1233861} and Kucerovsky \cite{MR1860509}
for single operators.
\end{remark}

\section{Perturbations of the identity}\label{Minusinfinity}

The discussion in the previous section of the index in the scattering case
allows us to compute the index of Fredholm perturbations of
the identity by fibred cusp operators of order $-\infty,$ \ie to give
analogues of Proposition~\ref{faficu.131} for all fibred cusp
structures. This was done in \cite{rochon} for the numerical index.

\begin{proposition}\label{faficu.133} Suppose that $E$ is a complex vector
bundle over the total space of a fibration with fibred cusp structure and
$b\in\fsP{\boumap}{-\infty}{\bourel;E}$ is such that
$\Id+b$ is invertible, then any family $\Id+B,$
$B\in\fcP{\boumap}{-\infty}{\famrel;E}$ with $N(\Id+ B)=\Id+b$ has analytic 
index equal to the image
\begin{equation}
[\Id+b]\in\Kco(\bbR\times T^*(\boubasrel))=\Kc(T^*(\boubasrel))\overset{\indAS}
\longrightarrow \Kt(\fambas) 
\label{faficu.134}\end{equation}
under the Atiyah-Singer index map for the fibration of $\boubas$ over $\fambas.$
\end{proposition}

\begin{proof} First we may use an `excision' construction
  to replace $\famtot$ by the simpler manifold $\pa\famtot\times[0,1]_x.$
  Indeed, taking a product decomposition of $\famtot$ near the boundary and
  identifying it with a neighbourhood of $x=0$ in the product, the
  quantization of $b$ may be localized to vanish outside this
  neighbourhood, \ie to have kernel vanishing outside the product of this
  neighbourhood with itself, and so can be identified with an operator on
  the model product with the same index. Thus it suffices to consider the
  product case $\famtot=\pa\famtot\times[0,1],$ with a bundle lifted from
  the boundary and with $b$ trivial at the $x=1$ boundary.

Thus the boundary fibration extends to the whole space as a fibration over
$[0,1]\times\boubas$ and the strategy is to reduce the problem to the
scattering calculus on this space. Consider a family of smoothing
projections $\Pi_N$ as in Section~\ref{Analyticindex} for the bundle $E$
over the fibration of the boundary given by $\boumap$ extended to act
trivially in the variable $x\in[0,1].$ Then $\Pi_Nb\Pi_N \longrightarrow b$
as $N\to\infty$ uniformly on $\bbR\times T^*(\boubasrel)$ in view of the rapid
decay. Thus we may replace $b$ by $\Pi_Nb\Pi_N$ for sufficiently large $N$
and hence assume that it acts on some finite rank subbundle of the smooth
sections of $\CI(\pa\famtot/\boubas)$ (namely the range of $\Pi_N)$ pulled
back to $\bbR\times T^*(\boubasrel).$ Quantizing $b$ to an operator
$B\in\fcP{\boumap}{-\infty}{\famtot;E}$ we may arrange that $\Pi_NB=B=B\Pi_N$
by replacing $B$ by $\Pi_NB\Pi_N.$ Note that $\Pi_N$ is not a fibred cusp
pseudodifferential operator, because it kernel is singular on the fibre
diagonal, however its composite with an element of
$\fcP{\boumap}{-\infty}{\famtot;E}$ is in the same space; this follows from
an examination of the kernels, see Appendix~\ref{Appendix-PT}. Now in fact the
same local analysis shows that
$\Pi_NB\Pi_N\in\Psi_{\scat}^{-\infty}([0,1]\times\boubasrel;W_N)$ has boundary
symbol $b,$ where $W_N$ is the range of $\Pi_N.$ Thus the result follows
from Proposition~\ref{faficu.131}.
\end{proof}

\section{Analytic classes in $\KK$ theory}\label{KK.0}

Baum, Douglas and Taylor, in \cite{Baum-Douglas-Taylor}, associate a
$\KK$-class with a Dirac operator on a manifold with boundary with
Atiyah-Patodi-Singer boundary condition. In this section, we extend, and
refine, their construction to the situation of families of fully elliptic
fibred cusp operators. This can also be thought as an adaptation to the
context of fibred cusp operators of a similar discussion for
b-pseudodifferential operators in \cite{MR99a:58144}, except that even in
that special case, here we associate to a Fredholm b-pseudodifferential
operator on $X$ a class in the K-homology of $X/\pa X,$ the manifold with
the boundary smashed to a point, rather than the absolute (or relative)
K-homology (see the final remark of \cite{Baum-Douglas-Taylor}). Although
small, this is an important difference in that it is at the heart of our
formulation of a families index theorem in K-theory, including the
Atiyah-Patodi-Singer case.

For a quick review of $\KK$-theory and the relation with elliptic operators
see \cite{Baum-Douglas-Taylor} and \cite{Melrose-Piazza1} and the books 
\cite{Blackadar1} and \cite{Higson-Roe} as well as the papers 
\cite{kasparov1}, \cite{kasparov2}, and \cite{kasparov} where KK-theory
was initially developed. To describe families of elliptic operators via
KK-theory we recall the extra feature introduced in \cite{kasparov}.

\begin{definition} If $X$ is a compact manifold and $\cA$ a $\bbZ_2$-graded
$C^{*}$ algebra then a $\CO(X)$-algebra structure on $\cA$ is a graded unital
homomorphism
\begin{equation*}
    r: \CO(X)\longrightarrow Z(\mathcal{M}(\cA))
\end{equation*}
where $Z(\mathcal{M}(\cA))$ is the center of $\mathcal{M}(\cA),$ the multiplier
algebra of $\cA,$ and where $\CO(X)$ is trivially graded; in particular
this gives $\cA$ a $\CO(X)$-module structure. 
\label{kk.1}\end{definition}

\begin{definition}
Let $(\cA,r_{\cA})$ and $(\cB,r_{\cB})$ be graded $\CO(X)$-algebras where
$X$ is a  compact manifold. Then $\Ef{X}{\cA,\cB}$ is the set of all triples
$(E,\phi,F)$ where $E$ is countably generated graded Hilbert module over $\cB,$
$\phi$ is a graded $*$-homomorphism from $\cA$ to $\mathcal{B}(E)$ and $F$ is
an operator in $\mathcal{B}(E)$ of degree 1, such that for all 
$a\in\cA,$ $b\in\cB,$ $e\in E$ and $f\in \CO(X),$
 \begin{equation}
\begin{aligned}
&\text{(i) }[\phi(a), F]\in \Kv(E),\\
&\text{(ii) } \phi(a)(F^{2}-\Id) \in \Kv(E),\\
&\text{(iii) }\phi(a)(F-F^{*}) \in \Kv(E),\\
&\text{(iv) }\phi(a\cdot r_{\cA}(f))(e\cdot b)= 
\phi(a)(e\cdot (r_{\cB}(f)\cdot
b)).
\end{aligned}
\label{kk.3}\end{equation}
Here $\Kv(E),$ defined for instance in \cite{Blackadar1}, is the analog
of compact operators for Hilbert modules. The elements of $\Ef{X}{\cA,\cB}$
are called Kasparov $\CO(X)$-modules for $(\cA,\cB).$  We denote by
$\bbD_{X}(\cA,\cB)$ the set of triples in $\Ef{X}{\cA,\cB}$ for which
$[F,\phi(a)],$ $(F-F^{*})\phi(a)$ and $(F^{2}-1)\phi(a)$ vanish for all
$a\in\cA.$  The elements of $\bbD_{X}(\cA,\cB)$ are called degenerate
Kasparov modules.
\label{kk.2}\end{definition}

Condition (iv) is the extra feature needed to deal with families of 
elliptic operators. It requires equivariance for the 
$\CO(X)$-module structure of $\cA$ and $\cB.$  One recovers the definition
of standard Kasparov modules by dropping (iv).

An element $(E,\phi,F)\in \Ef{X}{\cA,\CO([0,1];\cB)}$ generates a family 
\begin{equation*}
    \{ (E_{t},\phi_{t},F_{t})\in\Ef{X}{\cA,\cB}; t\in [0,1]\} 
\end{equation*}
obtained by evaluation at each $t\in [0,1].$ This family and the triple
itself will be called a homotopy between $(E_{0},\phi_{0},F_{0})$ and
$(E_{1},\phi_{1},F_{1})$ and these modules are then said to be homotopic
with the relation written
\begin{equation*}
(E_{0},\phi_{0},F_{0})\sim_{h}(E_{1},\phi_{1},F_{1}).
\end{equation*}

See for example \cite{Blackadar1} for a proof that all degenerate Kasparov
modules are homotopic to the trivial Kasparov $\CO(X)$-module.

\begin{definition}
We denote by $\Kf{X}{\cA,\cB}$ the set of equivalence
classes of $\Ef{X}{\cA,\cB}$ under the equivalence relation $\sim_{h}$ and
similarly define $\Kfo{X}{\cA,\cB}$ by
\begin{equation*}
         \Kfo{X}{\cA,\cB}=\Kf{X}{\cA,\mathcal{S}\cB} 
\end{equation*}
where 
\[
\mathcal{S}\cB= \{ f\in \CO([0,1]; \cB); \; f(0)=f(1)=0 \}
\] is the the suspension of the $C^{*}$
algebra $\cB.$
\label{kk.6}\end{definition}
\noindent In the paper of Kasparov \cite{kasparov}, the notation $\mathcal{R}\KK(X;\cA,\cB)$ is
used, we prefer a more compact notation.

As discussed in \cite{Blackadar1}, there are various other useful
equivalence relations which give $\Kf{X}{\cA,\cB}.$  The set of equivalent
classes $\Kf{X}{\cA,\cB}$ is
an abelian group with addition given by direct sum 
\begin{equation*}
         [(E_{0},\phi_{0},F_{0})]+ [(E_{1},\phi_{1},F_{1})]=
  [(E_{0}\oplus E_{1},\phi_{0}\oplus \phi_{1}, F_{0}\oplus F_{1})]. 
\end{equation*}

We are now in a position to define the $\KK$-classes associated to fully
elliptic fibred cusp operators. In our situation, $X=\fambas$ is the base
of the fibration, while $\cA$ will be variously $\COo(\famtot),$
$\Cfc(\famtot),$ $\CO(D),$ etc., and $\cB=\CO(\fambas).$  In particular, we
will only consider commutative $C^{*}$ algebras which are trivially graded
and the $\CO(\fambas)$-algebra structure will always be the obvious one.
  
Let $P\in\fcP{\boumap}{m}{\famrel;\sE}$ be a family of fully elliptic fibred
cusp operators of order $m,$ where $\sE$ is a $\bbZ_2$-graded complex vector
bundle on $\famtot.$ Introducing a graded inner product on $\sE$ and family of
metrics on $\famrel$ it follows that $P^{*}P$ is also a family of fully
elliptic operators of order $2m$ with strictly positive symbol and indicial
family. By standard pseudodifferential constructions there is an
approximate inverse square-root $Q\in\fcP{\boumap}{-m}{\famrel;E_+}$
which is invertible and positive definite such that
\begin{equation}
Q^{2}\circ P^{*}P-\Id \in x^{\infty}\fcP{\boumap}{-\infty}{\famrel;E_+},
\label{KK.5}\end{equation}
where $x$ is a boundary defining function for $\pa\famtot.$ Consider the
family of operators $A=PQ\in\fcP{\boumap}{0}{\famrel;\sE}.$ It is fully
elliptic and by construction it is almost unitary in the sense that
\begin{equation}
A^{*}A-\Id \in 
 x^{\infty}\fcP{\boumap}{-\infty}{\famrel;E_+}\text{ and }
AA^{*}-\Id \in 
x^{\infty}\fcP{\boumap}{-\infty}{\famrel;E_-}
\label{KK.1}\end{equation}
are families of compact operators. Let $\Cfc(\famtot)$ denote the space of 
continuous functions on $\famtot$ which are constant along the fibres of
$\boumap$, that is 
\begin{equation}
  \Cfc(\famtot)= \{f\in\CO(\famtot);f\big|_{\pa\famtot}=
g\circ\Phi\text{ for some } g\in\CO(\boubas)\}.
\label{faficu.46}\end{equation}
The injection $\CO(\fambas)\hookrightarrow\Cfc(\famtot)$ gives
$\Cfc(\famtot)$ a $\CO(\fambas)$-algebra structure and
$\CO(\fambas)$ itself has the $\CO(\fambas)$-algebra structure given by the
identity map $\CO(\fambas)\longrightarrow \CO(\fambas).$ Also let
$\Cfc^{\infty}(\famtot)\subset\Cfc(\famtot)$ denote the subspace of those 
functions which are smooth, thus 
\begin{equation*}
  \Cfc^\infty(\famtot)= \{f\in\CI(\famtot);f\big|_{\pa\famtot}=
g\circ\Phi\text{ for some } g\in\CI(\boubas)\}.
\label{faficu.57}\end{equation*}
Let us denote by $\mu$ the action of the $C^*$ algebra $\Cfc(\famtot)$
through multiplication
\begin{equation*}
\mu(f)\in\mathcal{B}(\Hv),\
\Hv=\Ld(\famrel;\sE)=\Ld(\famrel;E_+)\oplus\Ld(\famrel;E_-),\ f\in\Cfc(\famtot).
\label{faficu.47}\end{equation*}

\begin{lemma}\label{Kasp-mod}The triple $(\Hv,\mu,\F)$ where
\begin{equation*}
   \F= \begin{pmatrix}
                      0 & A^{*} \\
                      A & 0 
               \end{pmatrix}\in \B(\Hv)
\label{faficu.48}\end{equation*}
gives rise to a well-defined Kasparov module in
$\Ef{\fambas}{\Cfc(\famtot),\CO(\fambas)}$ and hence a class 
\begin{equation}
[P]=[(\Hv,\mu,\F)]\in \Kf{\fambas}{\Cfc(\famtot),\CO(\fambas)}.
\label{faficu.81}\end{equation}
\end{lemma}

\begin{proof}
By Kuiper's theorem, $\Hv$ is a countably generated $\bbZ_2$-graded Hilbert
module  over $\CO(\fambas).$ To show that $(\Hv,\mu,\F)$ is a module in
the sense of Kasparov we need to check the following properties 
\begin{equation}
\begin{aligned}
&\text{(i) }[\mu(f), \F]\in \Kv(\Hv),\\
&\text{(ii) } \mu(f)(\F^{2}-\Id) \in \Kv(\Hv),\\
&\text{(iii) }\mu(f)(\F-\F^{*}) \in \Kv(\Hv),\\
&\text{(iv) }\mu(b_{1}f)(h\cdot b_{2})= \mu(f)(h\cdot (b_{2}b_{1})),
\end{aligned}
\quad
\forall\ f\in \Cfc(\famtot),\ b_{1},b_{2}\in\CO(\fambas)\text{ and }h\in \Hv.
\label{faficu.50}\end{equation}
Property (iv) is immediate. Property (iii) follows directly from the fact
that $\F^{*}=\F.$ Property (ii) is a consequence of \eqref{KK.1}.
To check property (i), we may restrict to $f\in\Cfc^{\infty}(\famtot)$
since these smooth functions are dense in $\Cfc(\famtot)$ and the map 
\begin{equation*}
           [\mu(\cdot),\F]: \Cfc(\famtot)\longrightarrow \B(\Hv)
\label{faficu.49}\end{equation*}
is continuous.

For any $f\in\Cfc^{\infty}(\famtot),$
$\mu(f)\in\fcP{\boumap}{0}{\famrel;\sE}$ has symbol $f$ and indicial family
which can be identified with $f\big|_{\pa\famtot}$ which is to say a
constant multiple of the identity on each fibre, so commuting with the
normal operator of any other element. Thus

\begin{equation*}
           [\mu(f),\F]\in x\fcP{\boumap}{-1}{\famrel;\sE}
\label{faficu.51}\end{equation*}
is a family of compact operators.
\end{proof}

A fully elliptic operator also defines a class in the group
$\KK(\Cfc(\famtot),\CO(\fambas));$ to get Poincar\'{e} duality, we need to
take into account the $\CO(\fambas)$-algebra structure.

\begin{lemma}\label{KK.4}
The class $[P]\in \Kf{\fambas}{\Cfc(\famtot),\CO(\fambas)}$ associated to a fully
elliptic family of fibred cusp pseudodifferential operators by
Lemma~\ref{Kasp-mod} does not depend on the choice of $Q$ in \eqref{KK.5}
and in fact only depends on the homotopy class of $P$ in the space of fully
elliptic operators.
\end{lemma}

\begin{proof} Any two choices of a family of positive definite approximate
square-roots differ by a family in
$x^\infty\fcP{\boumap}{-\infty}{\famrel;E_+}$ so the resulting families
$\F$ differ by compact families and hence define the same element in
$\Kf{\fambas}{\Cfc(\famtot),\CO(\fambas)}.$

To prove the second part of the lemma, let $p_{t}\in
\fcP{\boumap}{m}{[0,1]\times\famrel;\sE}$ be a smooth curve of families of fully
elliptic operators, where $t\in[0,1].$  Then there exists a smooth curve
$Q_{t}\in\fcP{\boumap}{-m}{[0,1]\times\famrel;E_+}$ of invertible
approximate inverse square-roots such that 
\begin{equation*}
 Q_{t}^{2}\circ P^{*}_{t}P_{t}-\Id\in
x^{\infty}\fcP{\boumap}{-\infty}{[0,1]\times\famrel;\sE}.
\label{faficu.53}\end{equation*}
Hence, $(\Hv, \mu, \F_{t})\in \bbE(\Cfc(\famtot),\CO(\fambas))$ and if
$A_{t}= P_{t}Q_{t}$ then
\begin{equation*}
\F_{t}=\begin{pmatrix}
                      0 & A_{t}^{*} \\
                      A_{t} & 0 
               \end{pmatrix} \in \B(\Hv),
\label{faficu.54}\end{equation*}
defines an operator homotopy between the modules $(\Hv, \mu, \F_{0})$
and $(\Hv, \mu, \F_{1}).$ This implies that $[P_{0}]=[P_{1}]$ in
$\Kf{\fambas}{\Cfc(\famtot),\CO(\fambas)}.$ 
\end{proof}

This $\KK$-class also behaves in the expected manner under direct sums, so
if $P\in \fcP{\boumap}{m}{\famrel;\sE}$ and
$R\in\fcP{\boumap}m{\famrel;\sF}$ are families of fully elliptic operators,
then
\begin{equation*}
     [P\oplus R] = [P]+ [R] \text{ in }\Kf{\fambas}{\Cfc(\famtot),\CO(\fambas)}.
\label{faficu.55}\end{equation*}

It follows that this construction defines a `quantization' homomorphism of
abelian groups
\begin{equation}
\quan:\fcK{\boumap}(\fammap)\longrightarrow
\Kf{\fambas}{\Cfc(\famtot),\CO(\fambas)}.
\label{faficu.58}\end{equation}
The analytical index of the family $P$ factors through this map. Let
$c_{\boumap}:\CO(\fambas)\longrightarrow\Cfc(\famtot)$ be the inclusion of
constant functions along the fibres of
$\fammap:\famtot\longrightarrow\fambas.$ Then, at the level of 
$\KK$-theory, $c_{\boumap}$ defines a contravariant functor
\begin{equation*}
 c_{\boumap}^{*}:\Kf{\fambas}{\Cfc(\famtot),\CO(\fambas)}\longrightarrow 
	     \Kf{\fambas}{\CO(\fambas),\CO(\fambas)}.
\label{faficu.56}\end{equation*}

\begin{lemma}\label{KK.7} Under the standard identification  
\begin{equation*}
\Kf{\fambas}{\CO(\fambas),\CO(\fambas)}\cong
\KK(\bbC,\CO(\fambas))\cong\Kt^{0}(\fambas),
\label{faficu.59}\end{equation*}
there is a commutative diagram
\begin{equation*}
\xymatrix{\fcK{\boumap}(\fammap)\ar[r]^(0.37){\quan}\ar[dr]_{\inda}&
\Kf{\fambas}{\Cfc(\famtot),\CO(\fambas)}\ar[d]^{c_{\boumap}^*}\\
&\Kt(\fambas)}
\label{faficu.60}\end{equation*}
\end{lemma}

\begin{proof}
This follows from the discussion in \cite{Blackadar1}, more precisely 
proposition 17.5.5, corollary 12.2.3 and paragraph 8.3.2. It is also a
simple consequence of the stabilization of the null space of $P,$ by
perturbation, of Lemma~\ref{faficu.41}.
\end{proof}

It is also possible to define a quantization map for
$\fcKo{\boumap}(\fammap).$ Given a joint symbol $(\sigma,N)\in
\jS{s\boumap}0(s\fammap;E_{+})$ representing a class in
$\fcKo{\boumap}(\fammap),$ let $P_{s}\in \fcP{s\boumap}{0}{\famtot\times I/
\fambas\times I ;E}$ be a family of fully elliptic fibred cusp operators with
joint symbol $(\sigma, N)$ such that
$P_{s}\big|_{\fambas\times\{0,1\}}\equiv \Id.$  Here, recall that 
$s\boumap: \pa\famtot\times I \to \boubas\times I$ is the boundary fibration
of \eqref{sK.0}.  Let $ Q_{s}\in
\fcP{s\boumap}{0}{\famtot\times I/\fambas\times I;E_{+}}$ be an approximate
positive definite inverse square root
\begin{equation*}
   Q_{s}^{2}\circ P_{s}^{*}P_{s}-\Id \in x^{\infty}\fcP{s\boumap}{-\infty}
{\famtot\times I /\fambas\times I, E_{+}}
\end{equation*} 
such that $Q\big|_{\fambas\times\{0,1\}}\equiv \Id,$
where $x$ is the boundary defining function for $\famtot.$
Consider the family of operators $A=PQ\in\fcP{s\boumap}{0}{\famtot\times I/
\fambas\times I;E_{+}}.$ It is fully elliptic,
$A\big|_{\fambas\times\{0,1\}}\equiv\Id,$  and by construction it almost
unitary in the sense that
\begin{equation*}
   A^{*}A-\Id \in x^{\infty}\fcP{s\boumap}{-\infty}
{M\times I/B\times I, E_{+}},\
AA^{*}-\Id \in x^{\infty}\fcP{s\boumap}{-\infty}
{M\times I/B\times I, E_{-}}
\end{equation*}
are families of compact operators which vanish on $\fambas\times\{0,1\}.$
Notice that there is a natural inclusion $\Cfc(\famtot)\subset
\mathcal{C}_{s\Phi}(\famtot\times I).$ Let us denote by $\mu_{s}$ the
action of the $C^{*}$ algebra $\Cfc(\famtot)$ through multiplication 
\begin{multline}
  \mu_{s}(f)\in\mathcal{B}(\mathcal{H}_{s}),\
\mathcal{H}_{s}= \Ld(\famtot\times I/B\times I;\sE)\\
=
\Ld(\famtot\times I/B\times I;E_{+})\oplus 
\Ld(\famtot\times I/B\times I;E_{-})
\end{multline}
for $f\in \Cfc(\famtot),$ where in this odd context $E_{+}=E_{-}.$
Let $\COS(\fambas)$ be the the $C^{*}$ algebra of continuous functions
\begin{equation}
 \COS(\fambas)=\{f\in\CO(\fambas\times I);f\big|_{\fambas\times
   \{0,1\}}\equiv 0\}.
\label{sK.2}\end{equation}

\begin{lemma}
The triple $(\mathcal{H}_{s},\mu_{s},\mathcal{F}_{s})$ where
\begin{equation*}
          \F_{s}= \begin{pmatrix}
                      0 & A^{*} \\
                      A & 0 
               \end{pmatrix} \in \B(\Hv_{s}),
\end{equation*}
gives rise to a well-defined Kasparov module $\Ef{\fambas}{\Cfc(\famtot),
  \COS(\fambas)}$ and a class 
\begin{equation*}
           [(\sigma,N)]=[P_{s}]\in \Kf{\fambas}{\Cfc(\famtot),\COS(\fambas)}
\end{equation*}
which only depends on the class $[(\sigma,N)]\in\fcKo{\boumap}(\fammap)$
\label{sK.3}\end{lemma}

\begin{corollary}
Under the standard identification 
\begin{equation*}
\Kf{\fambas}{\Cfc(\famtot),\COS(\fambas)}\cong
\Kfo{\fambas}{\Cfc(\famtot),\CO(\fambas)},
\end{equation*}
Lemma~\ref{sK.3} gives us a well defined quantization map
\begin{equation}
\quan:\fcKo{\boumap}(\fammap)\longrightarrow
\Kfo{\fambas}{\Cfc(\famtot),\CO(\fambas)}. 
\label{sK.5}\end{equation}
\label{sK.4}\end{corollary}

\section{Poincar\'e duality for cusp operators}\label{aco.0}

We shall prove Theorem~\ref{faficu.88} in the case of the cusp
structure. To do so, we first discuss the short exact sequence \eqref{faficu.23}.

Consider the collection of invertible joint symbols of cusp
pseudodifferential operators, 
\begin{multline}
\IJS{\cusp}0(\fammap,\sE)=\big\{(\sigma ,N)\in
\syM{\cusp}0{\famrel);\sE}\oplus\fsP{\boumap}0{\pa\famrel;\sE};\\
\sigma \big|_{\pa\famtot}=\sigma _0(N),\ (\sigma ^{-1},N^{-1})
\in\syM{\cusp}0{\famrel);\sE^-}\oplus\fsP{\boumap}0{\pa\famrel;\sE^-}\big\}
\label{faficu.85}\end{multline}
on $\bbZ_2$-graded bundles where
$\syM{\cusp}0{\famrel;\sE}=\CI(\cuS^*(\famrel);\hom(\sE)).$ This naturally
maps by restriction to the collection of invertible symbols 
\begin{equation}
\begin{gathered}
\IUS{\cusp}0(\fammap,\sE)=\left\{\sigma\in\syM{\cusp}0{\famrel;\sE};
\sigma^{-1}\in\syM{\cusp}0{\famrel;\sE^-}\right\},
\\
\sigma:\IJS{\cusp}0(\fammap,\sE)\longrightarrow\IUS{\cusp}0(\fammap,\sE).
\end{gathered}
\label{faficu.87}\end{equation}
In fact this map is surjective. That is, for every family of elliptic
symbols there does exist a family of invertible normal operators which is
compatible with it. This is an aspect of the cobordism invariance of the index
and is shown in this form in \cite{fipomb}. 

The clutching construction associates to each invertible symbol an element
of the compactly supported $\Kt$ theory of $\cuT^*(\famrel)$ which is
`absolute' with respect to the boundary of $\famtot.$ As in the
boundaryless case, the resulting element is stable under homotopy, bundle
isomorphisms of $E_\pm$ and stabilization, so \eqref{faficu.87} descends to
a surjective map
\begin{equation}
\cuK{\fammap}\longrightarrow \Kc(T^*(\famrel))
\label{faficu.84}\end{equation}
where we use the fact that the cusp and standard cotangent bundles are
isomorphic, with the isomorphism natural up to homotopy. Thus
\eqref{faficu.84} is a surjective group homomorphism.

If $E_+=E_-$ we can consider those invertible suspended families of
pseudodifferential operators which are smoothing perturbations of the
identity on a fixed bundle over $\pa\famtot$
\begin{equation}
\Gs(\fammap;E)=\left\{N\in\Id_E+\Psi_{\sus}^{-\infty}(\pa\famrel;E);
N^{-1}\in\Id_E+\Psi_{\sus}^{-\infty}(\pa\famrel;E)\right\}.
\label{faficu.86}\end{equation}
Let $\Gs(\fammap;*)$ be the union over $E$ and let $\Gsp(\fammap;*)$ be the
subset corresponding to bundles $E$ which bound a bundle over $\famtot.$
Since we may complement a bundle to be trivial, and hence extendible, the
stable homotopy classes of elements of $\Gsp(\fammap;*)$ and
$\Gs(\fammap;*)$ are the same.

\begin{lemma}\label{faficu.89} Passing to the set of stable homotopy
classes, with equivalence also under bundle isomorphisms, the inclusion and
restriction maps give the split short exact sequence of Abelian groups
\eqref{faficu.24}:
\begin{equation}
\xymatrix{\Gsp(\fammap;*)\ar[r]^i\ar[d]^{\sim}&
\IJS{\cusp}0(\fammap,*)\ar[r]^\sigma\ar[d]^{\sim}\ar[ld]_{\inda}&
\IUS{\cusp}0(\fammap,*)\ar[d]^{\sim}\\
\Kt(\fambas)\ar@<.5ex>[r]^{i}&
\cuK{\fammap}\ar@<.5ex>[r]^(0.4){\sigma}\ar@<.5ex>[l]^{\inda}&
\Kc(T^*(\famrel))\ar@<.5ex>[l]^(0.6){\inv}.}
\label{faficu.90}\end{equation}
\end{lemma}

\begin{proof} That the set of stable homotopy classes, also allowing smooth
identification of bundles, of $\Gsp(\fammap,*),$ and hence also
$\Gs(\fammap,*),$ is canonically identified with $\Kt(\fambas)$ is a
standard result (see for instance \cite{HHH}) when the fibration is
trivial, $\famtot=\famfib\times\fambas.$ It remains true in the general
case, this can be shown using the families of
projections $\Pi_N$ of Section~\ref{Analyticindex}, see also \cite{rochon}.

Since direct sums behave consistently, the resulting maps are group
homomorphisms and form a complex, since $\Gsp(\fammap,*)$ clearly maps to
the identity in $\IUS{\cusp}0(\fammap;*).$

We have already noted the surjectivity of the second map. To see exactness
in the middle, suppose that a compatible pair $(\sigma,N)$ induces a
trivial class in $\Kc(T^*(\famrel)).$ Since stabilization and the action of
bundle isomorphisms is the same on the full and symbolic data, we may
suppose that $\sigma$ is homotopic to the identity through elliptic
symbols. In particular, $E_+=E_-.$ Adding the homotopy variable as an
additional base variable, the surjectivity of the symbol map allows the
homotopy of symbols to be lifted to an homotopy of joint symbols, see
Remark~\ref{faficu.221}. Thus $(\sigma ,N)$ may be deformed by homotopy to
$(\Id,\Id+A)$ where necessarily $A\in\Psi_{\sus}^{-1}(\pa\famrel;E).$ By a
further small perturbation this is homotopic to $\Id+A\in\Gsp(\fammap,*),$
$A\in\Psi_{\sus}^{-\infty}(\pa\famrel;E),$ showing exactness at 
$\cuK{\fammap}.$ 
Injectivity of the first map follows from the fact that the index map 
provides a
right inverse for it, which is a consequence of the families index of
Proposition~\ref{faficu.133}.   Note that the existence of an invertible 
family with
a given elliptic symbol defines the map $\inv$ which shows that the
sequence splits.
\end{proof}

There is also a corresponding exact sequence at the level of KK-theory.

\begin{lemma}\label{aco.8}
The short exact sequence of $C^{*}$-algebras \eqref{faficu.31} leads to a
split short exact sequence
\begin{equation}
\Kf{\fambas}{\CO(\fambas),\CO(\fambas)}\overset{\iota^{*}}\longrightarrow
 \Kf{\fambas}{\Ccu(\famtot),\CO(\fambas)}\overset{s^{*}}
\longrightarrow\Kf{\fambas}{\COo(\famtot),\CO(\fambas)}.
\label{faficu.93}\end{equation}
\end{lemma}
\begin{proof}
From standard results in $\KK$-theory (cf. theorem 19.5.7 in 
\cite{Blackadar1}),
the short exact sequence \eqref{faficu.31} leads to a six-term exact sequence
\begin{equation}
\xymatrix{
\Kf{\fambas}{\CO(\fambas),\CO(\fambas)}\ar[r]^{\iota^{*}}&
\Kf{\fambas}{\Ccu(\famtot),\CO(\fambas)}\ar[r]^{s^{*}}&
\Kf{\fambas}{\COo(\famtot), \CO(\fambas)}\ar[d]^{\delta}\\
\Kfo{\fambas}{\COo(\famtot),\CO(\fambas)}\ar[u]^{\delta}&
\Kfo{\fambas}{\Ccu(\famtot),\CO(\fambas)}\ar[l]_{s^{*}}&
\Kfo{\fambas}{\CO(\fambas),\CO(\fambas)}\ar[l]_{\iota^{*}}}
\label{aco.10}\end{equation}
where both boundary homomorphisms, $\delta,$ are obtained by multiplying
by a specific element
$\delta_{\iota}\in\Kfo{\fambas}{\CO(\fambas),\COo(\famtot)}.$ More
precisely, under the identification of
$\Kfo{\fambas}{\CO(\fambas),\COo(\famtot)}$ with
$\Kf{\fambas}{S\CO(\fambas),\COo(\famtot)},$ $\delta_{\iota}=i^{*}u$ where
$u\in\Kf{\fambas}{C_{\iota},\COo(\famtot)}$ and $i:
S\CO(\fambas)\longrightarrow C_{\iota}$ is the natural inclusion.  Here,
$C_{\iota}$ is the mapping cone  
\begin{equation*}
\begin{gathered}
C_{\iota}=\{(x,f)\in\Ccu(\famtot)\oplus \COo([0,1)\times \fambas);\iota(x)=f(0)\},\
S\CO(\fambas)= \COo((0,1)\times \fambas)
\end{gathered}
\label{faficu.91}\end{equation*}
and $i(f)=(0,f)\in C_{\iota}$ for $f\in S\CO(\fambas).$  In this situation
there is an injective map 
\begin{equation*}
j:\COo([0,1)\times\fambas)\ni f\longmapsto (\phi^{*}(f(0)),f)\in C_{\iota}
\label{faficu.92}\end{equation*}
so we may interpret $i$ as a map from $S\CO(\fambas)$ to
$\COo([0,1)\times \fambas).$ Since $\COo([0,1)\times \fambas)$ is a
contractible $C^{*}$-algebra, the class $[i]$ of $i$ in
$\Kf{\fambas}{S\CO(\fambas), C_{\iota}}$ is zero. In particular, this means
that $i^{*}u=0$ since $i^{*}u$ can be interpreted as the Kasparov product
of $[i]$ with $u$ (see for example 18.4.2a in \cite{Blackadar1}). Thus
$\delta=0$ and we get the short exact sequence \eqref{faficu.93}.

To see that this short exact sequence splits, consider the natural
injective $C^{*}$-homomorphism
$i_{\fammap}:\CO(\fambas)\longrightarrow\Ccu(\famtot).$ This satisfies 
$\iota i_{\fammap}=\Id,$ so $i_{\fammap}^{*}$ 
is a left inverse for $\iota^{*}.$
\end{proof}

There is a correspondence between the short exact sequences of
Lemmas~\ref{faficu.89} and \ref{aco.8}.  Consider the diagram 
\begin{equation}
\xymatrix{
\Kt(\fambas)\ar[r]^{i_{*}}\ar[d]_{\qb}&
\cuK{\phi}\ar[r]^{\sigma_{*}}\ar[d]_{\quan}&
\Kc(T^{*}(\famrel))\ar[d]_{\quanp}\\
\Kf{\fambas}{\CO(\fambas),\CO(\fambas)}\ar[r]^{\iota^{*}}&
\Kf{\fambas}{\Ccu(\famtot),\CO(\fambas)}\ar[r]^{s^{*}}&
\Kf{\fambas}{\COo(\famtot),\CO(\fambas)}}
\label{aco.11}\end{equation}
where the top row is the short exact sequence of Lemma~\ref{faficu.89} and the 
bottom row is the short exact sequence of Lemma~\ref{aco.8}. The central
map is the quantization map of Lemma~\ref{KK.7}. Similarly, the map
$\quanp$ is the quantization map as discussed in
\cite{Melrose-Piazza1}.

The map $\qb$ on the left is essentially the same quantization map as
in Lemma~\ref{KK.7}. Thus if $N\in\Gsp(\fammap,E_{+})$ there exists $P=\Id+L,$
$L\in\cuP{-\infty}{\famrel;E_{+}}$ with $N(P)=N.$ Choose an invertible,
positive, approximate inverse square-root,
$Q\in\Id+\cuP{-\infty}{\famrel;E_{+}},$
\begin{equation*}
Q^{2}\circ P^{*}P -\Id \in x^{\infty}\fcP{\boumap}{-\infty}{\famrel;E_{+}}
\label{faficu.94}\end{equation*}
and consider $A=PQ.$ Let $\Hv_{\fambas}$ be the $C(\fambas)$-Hilbert module 
\begin{equation}
      \Hv_{\fambas}=\Hv_{\fambas}^{+}\oplus\Hv_{\fambas}^{-}=
\Ld(\famtot;E)\oplus \Ld(\famtot;E).
\label{aco.12}\end{equation}
Then $\qb(N)$ is the $\KK$-class associated to the Kasparov module
$(\Hv_{\fambas},\mu_{\fambas},\F_{\fambas})$ in
$\Ef{\fambas}{\CO(\fambas),\CO(\fambas)}$ where
\begin{equation}
    \F_{\fambas}=\begin{pmatrix}
                 0 & A^{*} \\
                 A &  0 
                \end{pmatrix}
\label{aco.15}\end{equation}
and
$\mu_{\fambas}:\CO(\fambas)\longrightarrow\mathcal{\fambas}(\Hv_{\fambas})$
is just given by multiplication using the $\CO(\fambas)$-module structure.

\begin{proposition}
The diagram~\eqref{aco.11} is commutative and $\qb ,$ $\quan$
and $\quanp$ are isomorphisms.
\label{aco.13}\end{proposition}

\begin{proof} The fact that $\quanp$ is an isomorphism is established in
\cite{Melrose-Piazza1}, where the Poincar\'{e} duality of \cite{kasparov} is
generalized. From the definition of $\quan$ and $\quanp,$ it is
straightforward to check that the right square of the diagram is
commutative.

That $\qb$ is an isomorphism follows from the identification
\begin{equation*}
\Kf{\fambas}{\CO(\fambas),\CO(\fambas)}\cong
\KK^{0}(\bbC,\CO(\fambas))\cong \Kt(\fambas)
\label{faficu.95}\end{equation*}
and that under this $\qb$ gives the index map.

The commutativity of the left square is not quite obvious since it does not
commute in terms of Kasparov modules. Indeed, if $P\in\Id+K$ is as
discussed above then $\iota^{*}\qb(P)$ is represented by the Kasparov
module $(\Hv_{\fambas}, \iota^{*}\mu_{\fambas},\F_{\fambas})$ with
$\Hv_{\fambas}$ and $\F_{\fambas}$ are as in \eqref{aco.12} and
\eqref{aco.15} whereas $\quan i_{*}(P)$ is represented by the Kasparov module
$(\Hv_{\fambas},\mu,\F_{\fambas})$ with $\mu: \Ccu(\famtot)\longrightarrow 
\mathcal{\fambas}(\Hv_{\fambas})$ given by multiplication by $\Ccu(\famtot).$
Since
\begin{equation*}
     \iota^{*}\qb(P)= \iota^{*}(i_{\fammap})^{*}\quan i_{*}(P)
\label{faficu.97}\end{equation*}
as Kasparov modules, the commutativity of the second square of the diagram,
$s^{*}\quan i_{*}(P)=0$ shows that
\begin{equation*}
     \iota^{*}\qb(P)=\iota^{*}(i_{\fammap})^{*}\quan i_{*}(P)=\quan i_{*}(P)\Min
  \Kf{\fambas}{\Ccu(\famtot),\CO(\fambas)}
\label{faficu.96}\end{equation*}
and so these two modules give the same element in the K-group.

This implies that $\quan$ is also an isomorphism.
\end{proof}

\section{The 6-term exact sequence}\label{stes.0}

As noted above, it is always possible to perturb an elliptic family of cusp
operators by a cusp operator of order $-\infty$ so that it becomes
invertible and this leads to the short exact sequence \eqref{faficu.24}. In
the general fibred cusp case there is an obstruction in K-theory to the
existence of such a perturbation and this results in the 6-term exact
sequence \eqref{faficu.36}.

Consider a fibration with fibred cusp structure, as in \eqref{faficu.2}, so
that the algebra $\fcP{\boumap}{*}{\famrel}$ of families of fibred cusp
operators is well-defined. Let $r_{\pa\famtot}$ denote the inclusion
\begin{equation}
r_{\pa\famtot}:T^{*}_{\pa\famtot}(\famrel)\longhookrightarrow T^{*}(\famrel).
\label{efe.1}\end{equation}
If 
\begin{equation}
    d\Phi:T_{\pa\famtot}(\famrel)\longrightarrow T(\boubasrel)\times \bbR
\label{efe.3}\end{equation}
is the (extended) differential of $\boumap,$ then identifying the tangent
bundles with the cotangent bundles via some choice of metrics, there is a
well-defined families index 
\begin{equation*}
\ind_{\AS}: \K(T^{*}_{\pa\famtot}(\famrel))\longrightarrow
\Kc(T^{*}(\boubasrel)\times\bbR)\cong \Kco(T^{*}(\boubasrel)).
\end{equation*}

\begin{proposition}
An elliptic family, $P\in\fcP{\boumap}{m}{\famrel;\sE},$ can be perturbed
to be fully elliptic by the addition of some
$Q\in\fcP{\boumap}{-\infty}{\famrel;\sE}$ if and only if the index of the
class of its symbol at the boundary
\begin{equation}
\ind_{\AS} [r_{\pa\famtot}^{*}\sigma(P)]=0 \in \Kco(T^{*}(\boubasrel)) 
\label{faficu.222}\end{equation}
and then $Q$ can be chosen so that $P+Q$ is invertible.
\label{efe.2}\end{proposition}

\begin{proof}
By the families index theorem of Atiyah and Singer, there is a suspended
perturbation $I(Q)$ of order $-\infty$ such that $I(P)+I(Q)$ is invertible
if and only if \eqref{faficu.222} holds. With such a choice $P+Q$ is a
fully elliptic family. Since it follows from Proposition~\ref{faficu.133}
that the index map on perturbations of the identity is surjective, we may
compose on the left with a Fredholm family of the form
$\Id+\fcP{\boumap}{-\infty}{\famrel;E_{-}}$ to get an operator of index
zero and then perturb, as in Section~\ref{Analyticindex} to make it invertible.
\end{proof}

In the particular case of a single elliptic cusp operator, the obstruction
is in $\Kto(\pt)\cong\{0\},$ so there is no obstruction to such a
perturbation. In the case of a family of elliptic cusp operators, the
obstruction is in $\Kto(\boubas),$ which is not the trivial group in
general. However $\ind_{\AS}\circ r_{\pa\famtot}^{*}\kappa(P)$ is
always zero in $\Kto(\boubas)$ by the cobordism invariance of the index.

There is a parallel discussion in the odd case.  The fibration \eqref{efe.3} 
also induces an odd index
\begin{equation}
\ind_{\AS}: \Kco(T^{*}_{\pa\famtot}(\famrel))\longrightarrow 
\Kco(T^{*}(\boubasrel)\times \bbR)\cong \Kc(T^{*}(\boubasrel)).
\label{efe.3a}\end{equation}

\begin{proposition}\label{efe.4}
Suppose $P\in\fcP{\boumap}{0}{\famtot\times[0,1]/\boubas\times[0,1];E}$ 
is a family of elliptic operators which is the identity at $t=0$ and $t=1$
and let $[\sigma(P)]\in\Kco(T^{*}(\famrel))$ be the class of its symbol,
then $P$ can be perturbed by
$Q\in\fcP{\boumap}{-\infty}{\famtot\times[0,1]/\boubas\times[0,1];E}$ 
with $Q\big|_{\boubas\times\{0,1\}}=0$ to be invertible if and only if
\begin{equation*}
\ind_{\AS} r_{\pa\famtot}^{*}[\sigma(P)]=0\in\Kc(T^{*}(\boubasrel)).
\end{equation*}
\end{proposition}

\begin{proof}
One could proceed as in the proof of Proposition~\ref{efe.2}.  However, there
is an alternative proof which is more suggestive in this case.  As discussed
in \cite{HHH}, one can define the odd index \eqref{efe.3a} in the following
way. Notice that $r_{\pa\famtot}^{*}\sigma(P),$ which is the symbol of the indicial family, gives a class in
\begin{equation*}
\Kco(T^{*}_{\pa\famtot}(\famrel))\cong\Kco(T^{*}(\pa \famrel)\times \bbR).
\end{equation*}
If $t\in [0,1]$ denotes the suspension parameter, then by assumption 
$r_{\pa M}^{*}\sigma(P)$ is the identity at $t=0$ and $t=1.$ So there
is no obstruction to perturb $P$ by smoothing operators so that one gets
a 1-parameter family 
\begin{equation}
  t\mapsto N(P_{t})\in \Gsbn{1}(\pa \famrel;E)
\label{efe.5}\end{equation}
such that $N(P_{0})=\Id$ and $N(P_{1})\in \Gsfc(\pa \famrel;E),$
where $\Gsbn{1}(\pa \famrel;E)$ is the group of invertible elliptic 
fibred suspended operators of order 0. Via the identification (if $\dim
\boubas=\dim \pa \famtot,$ this is only true after stabilization)
\begin{equation}
  \pi_{0}(\Gsfc(\pa \famrel;E))\cong \Kco(T^{*}(\boubasrel)),
\label{efe.6}\end{equation}
the connected component in which $N(P_{1})$ lies gives a class in
$\Kco(T^{*}(\boubasrel))$ which is precisely $\ind_{\AS}[r^{*}_{\pa 
\famtot} \sigma(P)].$  In particular, it does not depend on the choice 
of the 1-parameter family \eqref{efe.5}.

Now, if the index is zero, this means $N(P_{1})$ is in the connected component
of the identity in $\Gsfc(\pa \famrel;E),$ so via some smooth deformation,
it can be arranged that $N(P_{1})=\Id$ as well, so $P$ can be perturbed by
$Q\in\fcP{\boumap}{-\infty}{\famrel;E}$ with
$Q\big|_{\boubas\times\{0,1\}}=0,$ to become fully elliptic. Conversely, if
such a perturbation $Q$ exists, then the index is given by the K-class
corresponding to $N(P_{1}+Q_{1})=\Id,$ which is necessarily zero.
\end{proof}

The obstruction result of Proposition~\ref{efe.2} and Proposition~\ref{efe.4}
indicates that the 
short exact sequence of Lemma~\ref{faficu.89} fails to be exact in the
more general setting of fibred cusp operators.  However, as in the case of
the K-theory of a pair of spaces, there is a 6-term exact sequence given by
\begin{equation}
\xymatrix{\Kc(T^*(\boubasrel))\ar[r]^{i_{0}}&
\fcK{\boumap}(\fammap)\ar[r]^(.4){\sigma_{0}}&
\Kc(T^*(\famrel))\ar[d]^{I_{0}}\\
\Kco(T^*(\famrel))\ar[u]_{I_{1}}&
\fcKo{\boumap}(\fammap)\ar[l]_(.35){\sigma_{1}}&
\Kco(T^*(\boubasrel)).\ar[l]_{i_{1}} }
\label{stes.2}\end{equation}

To define $i_{k}$ for $k\in\bbZ_{2},$ consider the group
\begin{equation*}
 \Gsfcn{1+k}(\pa\famrel)\cong\{\Id+Q;Q\in
\fsn{1+k}^{-\infty}(\pa\famrel),\Id+Q\mbox{ is invertible}\}.
\end{equation*}
Then using spectral sections techniques as in \cite{rochon} or the
projections in Section~\ref{Analyticindex}, one has an identification 
\begin{equation}
   \Kcn{-k}(T^*(\boubasrel))\cong \pi_{0}(\Gsfcn{1+k}(\pa \famrel)). 
\label{stes.3}\end{equation}
Strictly speaking, the result is only true provided $\dim X>0,$ but in the
case $\dim X=0,$ it is only necessary to allow some stabilization.
Any element of $\Gsfcn{1+k}(\pa \famrel)$ can be seen as the indicial 
family of a family of fully elliptic operators with symbol given by the 
identity. This gives a map
\begin{equation}
 \pi_{0}(\Gsfcn{1+k}(\pa \famrel))\longrightarrow \fcKn{-k}{\boumap}(\fammap)
\label{stes.4}\end{equation}
and we define $i_{k}$ by composing \eqref{stes.3} with \eqref{stes.4}.
The symbol maps have already been defined and the boundary map $I_{0}$ is
given by
\begin{equation*}
I_{0}= \ind_{\AS}\circ r_{\pa\famtot}^{*}: \Kc(T^*(\famrel))\longrightarrow 
\Kco(T^*(\boubasrel)),              
\end{equation*}
while $I_{1}$ has a similar definition 
\begin{equation*}
I_{1}= \ind_{\AS}\circ r_{\pa\famtot}^{*}: \Kco(T^*(\famrel)) \longrightarrow 
  \Kcn{-1}(T^{*}(\boubasrel)\times \bbR) \cong \Kc(T^*(\boubasrel)),
\end{equation*}
where the last identification is by Bott periodicity.

\begin{theorem}
The diagram \eqref{stes.2} is exact.
\label{stes.5}\end{theorem}

\begin{proof} The exactness at $\fcK{\boumap}(\fammap)$ and
$\fcKo{\boumap}(\fammap)$ follows rather directly from the definition since
homotopies of the symbol can be lifted to homotopies of the joint symbol,
see the discussion above in the cusp case. Exactness at $\Kc(T^*(\famrel))$
and $\Kco(T^*(\famrel))$ follows from Proposition~\ref{efe.2} and its
suspended version.

The exactness at $\Kc(T^{*}(\boubasrel))$ can be seen from the alternative
definition of the odd index used in the proof of Proposition~\ref{efe.4}.
Indeed, suppose that $\alpha\in \Kc(T^{*}(\boubasrel))$ is in the image of
$I_{1}.$  According to the proof of Proposition~\ref{efe.4}, this means
that there is a family of fully elliptic operators $t\mapsto P_{t}\in
\fcP{\boumap}{0}{\famrel;E}$ such that $P_{0}=\Id,$ $\sigma(P_{1})=\Id$ and 
\begin{equation*}
           [N(P_{1})]\in \pi_{0}(\Gsfc(\pa \famrel; E)\cong
\Kc^{0}(T^{*}(\boubasrel))
\end{equation*}
corresponds to $\alpha.$ The homotopy $t\mapsto (\sigma(P_{t}), N(P_{t}))$
between the identity and $(\Id,N(P_{1}))$ then shows that $i(\alpha)=0$ in 
$\fcK{\boumap}(\fammap).$  

Conversely, suppose that $i(\alpha)=0.$ If $N(P_{1})\in \Gsfc(\pa \famrel;
E)$ represents $\alpha,$ then after some stabilization, one can assume
that there is a homotopy of invertible joint symbols $t\mapsto
(\sigma(P_{t}),N(P_{t}))$ between $(\Id, \Id)$ and $(\Id, N(P_{1})).$ The
family symbol $t\mapsto \sigma(P_{t})$ then defines a class $\beta\in
\Kco(T^{*}(\famrel))$ such that $I_{1}(\beta)=\alpha,$
which establishes the exactness at $\Kc(T^{*}(\boubasrel)).$

To prove the exactness at $\Kco(T^{*}(\boubasrel)),$ 
one can proceed in a similar
way. Indeed, the diagram
\begin{equation}
\xymatrix{
\Kcn{0}(T^{*}(\famrel))\ar[r]^{I_{0}}\ar[d]^{b}& \Kcn{-1}(T^{*}(\boubasrel))
\ar[dr]^{i_{1}}\ar[d]^{b} &   \\
\Kcn{-2}(T^{*}(\famrel))\ar[r]^{I_{-2}}& \Kcn{-3}(T^{*}(\boubasrel))
\ar[r]^{j} & \fcKo{\boumap}(\fammap)}
\label{efe.6a}\end{equation}
commutes, where the vertical arrows are given by Bott periodicity and
$I_{-2}=\ind_{\AS} r_{\pa M}^{*}$ in terms of the families index map
\begin{equation*}
   \ind_{\AS}: \Kcn{-2}(T^{*}_{\pa\famtot}(\famrel))\to 
\Kcn{-2}(T^{*}(\boubasrel)\times \bbR)\cong \Kcn{-3}(T^{*}(\boubasrel))
\end{equation*}
and  $j$ is the map \eqref{stes.4} obtained using the identification
(when $\dim \boubas=\dim\pa \famtot,$ this is only true after stabilization)
\begin{equation}
          \Kcn{-3}(T^{*}(\boubasrel))\cong \pi_{0}(\Gsfcn{2}(\pa\famrel)).
\label{efe.30}\end{equation}
Notice that as opposed to \eqref{stes.3}, no Bott periodicity is involved
in \eqref{efe.30}, hence the right triangle in \eqref{efe.6a} is
commutative. The fact that the left square is commutative follows from the 
commutativity of the families index with Bott periodicity (and more generally
with the Thom isomorphism).

The exactness of the bottom row of \eqref{efe.6a} can be proved by applying
the proof of the exactness at $\Kc(T^{*}(\boubasrel))$ (really
$\Kcn{-2}(T^{*}(\boubasrel)))$ with one extra suspension parameter.  Since
the Bott periodicity maps in \eqref{efe.6a} are isomorphisms, 
this implies that 
\begin{equation*}
\Kc(T^{*}(\famrel))\overset{I_{0}}{\longrightarrow}
\Kcn{-1}(T^{*}(\boubasrel))\overset{i_{1}}{\longrightarrow}
\fcKo{\boumap}(\fammap) 
\end{equation*}
is also exact in the middle.
\end{proof}

In the scattering case (when stabilization is necessary in the arguments
above), using Lemma~\ref{faficu.128}, the 6-term
exact sequence \eqref{stes.2} reduces to the 
exact sequence in K-theory associated
to the inclusion of the boundary of $T^*(\famrel)$
\begin{equation}
\xymatrix{\Kco(T^{*}_{\pa\famtot}(\famrel))\ar[r]&
\Kc(T^{*}(\famrel),T^*_{\pa\famtot}(\famrel))\ar[r]&
\Kc(T^{*}(\famrel))\ar[d]\\
\Kco(T^*(\famrel))\ar[u]&
\Kco(T^*(\famrel),T^{*}_{\pa\famtot}(\famrel))\ar[l]&
\Kc(T^*_{\pa\famtot}(\famrel)).\ar[l] }
\label{stes.6}\end{equation} 

\section{Poincar\'e duality, the general case}\label{GeneralPoincare}

Consider the short exact sequence of $C^{*}$ algebras \eqref{faficu.37}
and the associated 6-term exact sequence
\begin{equation}
\xymatrix{\Kf{\fambas}{\CO(\boubas),\CO(\fambas)}\ar[r]^{\iota^{*}}&
\Kf{\fambas}{\CO_{\boumap}(\famtot),\CO(\fambas)}\ar[r]^{s^{*}}&
\Kf{\fambas}{\COo(\famtot),\CO(\fambas)}\ar[d]^{\delta}\\
\Kfo{\fambas}{\COo(\famtot),\CO(\fambas)}\ar[u]^{\delta}&
\Kfo{\fambas}{\CO_{\Phi}(\famtot,\CO(\fambas)}\ar[l]^{s^{*}}&
\Kfo{\fambas}{\CO(\boubas),\CO(\fambas)}\ar[l]^{\iota^{*}}.}
\label{gs.2}\end{equation}
Each term of this sequence can be related to the corresponding term in
$\eqref{stes.2}$ via a quantization map, namely the quantization maps of
Section~\ref{KK.0}
\begin{equation}
\quan:\fcKn{k}{\boumap}(\fammap)\longrightarrow
\Kfn{k}{\fambas}{\CO_{\boumap}(\famtot),\CO(\fambas)},\ k\in\bbZ_2,
\label{gs.3}\end{equation} 
and the quantization maps of \cite{Melrose-Piazza1} and of Atiyah and Singer,
\begin{gather}
\qr: \Kcn{k}(T^{*}(\famrel))\longrightarrow
\Kfn{k}{\fambas}{\COo(\famtot),\CO(\fambas)},
\\
\qb: \Kcn{k}(T^{*}(\boubasrel))\longrightarrow 
\Kfn{k}{\fambas}{\CO(\boubas), \CO(\fambas)}.
\label{gs.4}\end{gather}

\begin{theorem}
The quantization maps $\quan,$ $\qr,$ and $\qb$ are isomorphisms giving a
commutative diagram
\begin{equation*}
\xymatrix{
\cdots\Kc(T^*(\boubasrel))\ar[r]^{i}\ar[d]^{\qb}&
\fcK{\boumap}(\fammap)\ar[r]^(.4){\sigma}\ar[d]^{\quan}&
\Kc(T^*(\famrel))\cdots\ar[d]^{\qr}\\
\cdots
\Kf{\fambas}{\CO(\boubas),\CO(\fambas)}\ar[r]^{\iota^{*}}&
\Kf{\fambas}{\CO_{\boumap}(\famtot),\CO(\fambas)}\ar[r]^{s^{*}}&
\Kf{\fambas}{\COo(\famtot),\CO(\fambas)}\cdots 
}
\end{equation*}
between the 6-term exact sequences \eqref{stes.2} and \eqref{gs.2}.  In 
particular, this implies the Poincar\'{e} duality of Theorem~\ref{faficu.88}.
\label{gs.5}\end{theorem}

The fact that $\qr$ and $\qb$ are isomorphisms follows from the Poincar\'{e}
duality result of Kasparov \cite{kasparov} and its extensions by the first author and Piazza in
\cite{Melrose-Piazza1}.  Then, given the commutativity of the 
diagram between the two 6-term exact sequence, the fact that $\quan$ is
an isomorphism follows from the fives lemma. So it remains to check that
the diagram commutes.

For $k\in\bbZ_{2},$ the commutativity of 
\begin{equation}
\xymatrix{\fcKn{k}{\boumap}(\fammap)\ar[r]^(.4){\sigma}\ar[d]^{\quan}&
\Kcn{k}(T^*(\famrel))\ar[d]^{\qr}\\
\KK^{k}_{\fambas}(\CO_{\boumap}(\famtot),\CO(\fambas))\ar[r]^{s^{*}}&
\KK^{k}_{\fambas}(\COo(\famtot),\CO(\fambas))
}
\label{gs.6}\end{equation}
is clear from the definition of $\quan$ and $\qr.$  The proof of the 
commutativity of the four other squares of the diagram is more involved and 
will occupy the remainder of this section.

Let us first consider the commutativity of 
\begin{equation}
\xymatrix{
\Kcn{k}(T^*(\boubasrel))\ar[r]^{i}\ar[d]^{\qb}&
\fcKn{k}{\boumap}(\fammap)\ar[d]^{\quan} \\
\Kfn{k}{\fambas}{\CO(\boubas),\CO(\fambas)}\ar[r]^{\iota^{*}}&
\Kfn{k}{\fambas}{\CO_{\boumap}(\famtot),\CO(\fambas)}.}
\label{gs.7}\end{equation}
We will only provide a proof in the even case since the odd case is similar.
The first step is to describe the quantization map $\qb$ in terms of 
indicial families instead of symbols. Given a class $\alpha\in\Kc(T^{*}
(\boubasrel)),$ let $p\in \Gsfc(\pa \famrel; E_{+})$ be a 
representative of this class.  Consider the manifold 
\begin{equation*}
         \cn= \pa\famtot\times [0,1]
\end{equation*} 
which can be seen as a collar neighborhood $\pa\famtot$ in $\famtot.$
Consider the associated fibration
\begin{equation}
\cnfm= \pa\fammap\circ \pi_{1}: \cn\longrightarrow B,
\label{gs.8}\end{equation}
where $\pi_{1}:\pa\famtot\times [0,1]\longrightarrow \pa\famtot $ is the projection
on the first factor.
Since the boundary of $\cn$ has two parts $\pa\famtot_{0}=\pa\famtot\times\{0\}$
and $\pa\famtot_{1}=\pa\famtot\times\{1\},$ let 
\begin{equation*}
         \boumap_{i}: \pa\famtot_{i}\longrightarrow\boubas_{i},\ i\in\{0,1\},
\end{equation*}
denote the two boundary fibration maps and let
\begin{equation*}
   \cnbm:\pa \cn\longrightarrow\boubas_{0}\cup\boubas_{1}
\end{equation*} 
be the total boundary map. Let $P\in 
(\Id+\fcP{\cnbm}{-\infty}{\cn; E_{+}})$ be a Fredholm
operator with indicial family the identity at $\pa\famtot_{0}$ and
by $p$ at $\pa\famtot_{1}.$  Then Lemma~\ref{Kasp-mod} gives an associated
$\KK$-class
\begin{equation*}
          [P]\in \Kf{\fambas}{\CO_{\cnbm}(\cn),\CO(\fambas)}
\end{equation*}
and the pull-back map 
\begin{equation*}
   d^{*}:\Kf{\fambas}{\CO_{\cnbm}(\cn), \CO(\fambas)}\longrightarrow
\Kf{\fambas}{\CO(\boubas), \CO(\fambas)}
\end{equation*}
where $d:\CO(\boubas)\hookrightarrow \CO_{\cnbm}(\cn)$ is the obvious 
inclusion, gives a KK-class
\begin{equation*}
                  d^{*}[P]\in \Kf{\fambas}{\CO(\boubas), \CO(\fambas)}.
\end{equation*}
Thus, this procedure gives a well-defined quantization map
\begin{equation*}
\qb':\Kc(T^{*}(\bourel))\longrightarrow\Kf{\fambas}{\CO(\boubas),\CO(\fambas)}. 
\end{equation*}
\begin{lemma}
$\qb'=\qb.$
\label{gs.10}\end{lemma}

\begin{proof}
Let $\alpha\in\Kc(T^{*}(\bourel))$ be given. Following the discussion in
\S\ref{Minusinfinity} we are reduced to the scattering case, so 
\begin{equation*}
           \cn=\boubas\times [0,1]
\end{equation*}
in our construction. Let $p$ be an indicial family representing the 
K-class $\alpha$ and let $P$ be as above.  Let
$P_{t}\in \fcP{\cnbm}{0}{\cn/\fambas;E_{+}},$ $t\in [0,1]$ be a homotopy 
through families of fully elliptic operators such that $P_{0}=P$ and
$P_{1}$ has trivial indicial families both at $\pa\famtot_{0}$ and
$\pa\famtot_{1}$ and so with symbol $\sigma$ having K-class 
\begin{equation*}
  [\sigma]\in \Kc(T^{*}(\cn/\fambas),T^{*}_{\pa \cn}(\cn/\fambas))
\cong \Kc( T^{*}(\boubas\times (0,1)/\fambas))
\end{equation*}
identifying via Bott periodicity with $\alpha$ (see \cite{rochon} for 
details). That such a homotopy exists is discussed in \cite{MelroseGST} and
\S\ref{Scattering}. By considering a family of positive definite
approximate inverse square roots for $P_{t},$ we construct an operator
homotopy $(\Hv, \mu, \F_{t})$ of Kasparov modules, which means that $P_{0}$
and $P_{1}$ define the same element in $\Kf{\fambas}{\CO(\boubas),\CO(\fambas)}.$

On the other hand, let us quantize partially the symbol $\sigma$ of $P_{1}$ in
the $T^{*}[0,1]$ direction to get a family of elliptic operators
$\widehat{p}_{1}$ parameterized by $T^{*}(\boubasrel)$ and acting on the
Hilbert bundle 
\begin{equation*}
           \Ld(\cn/D;E_{+})\longrightarrow T^{*}(\boubasrel)
\end{equation*}
with typical fibre $\Ld([0,1];E_{+}).$  As discussed in \cite{Luke}, we can
interpret
\begin{equation*}
    \widehat{p}_{1}: \Ld(\cn/D;E_{+})\longrightarrow \Ld(\cn/D;E_{+})
\end{equation*}
as a symbol on a Hilbert bundle.  The notion of ellipticity leads to
very restrictive conditions in this context.  But according to example
1.7 in \cite{Luke}, $\widehat{p}_{1}$ is elliptic and once it is quantized
we get back $P_{1}$ modulo compact operators.  Assume without loss of
generality that $\sigma$ (and $\widehat{p}_{1}$) is homogeneous of degree
zero in the fibres of $T^{*}(D/B)\longrightarrow D$ outside a compact neighborhood of the 
zero section.  Assume also without loss of generality that 
$E_{+}$ is a trivial bundle on $\cn.$
Then as discussed in \cite{Luke}, it is possible to deform
$\widehat{p}_{1}$ through a homotopy $\widehat{p}_{t},$ $t\in [1,2]$ of 
elliptic symbols so that $\widehat{p}_{2}$ takes the form
\begin{equation*}
    \widehat{p}_{2}=\widehat{p}_{2}' \oplus \widehat{p}_{2}'':
V\oplus V^{\perp} \longrightarrow \widehat{p}_{1}V
\oplus (\widehat{p}_{1}V)^{\perp} 
\end{equation*}   
where $V$ is a sub-bundle of $\Ld(\cn/D;E)$ of finite corank on which 
$\widehat{p}_{1}$ is injective and where $\widehat{p}_{2}'$ is invertible
and constant along the fibres of $T^{*}(D/B)\longrightarrow D.$  Thus, when we quantize
$\widehat{p}_{2}$ we get an operator of the form
\begin{equation*}
   P_{2}=P_{2}'\oplus P_{2}'' : \Ld(D;V)\oplus \Ld(D;V^{\perp})\longrightarrow
                \Ld(D;\widehat{p}_{1}V)\oplus \Ld(D; 
(\widehat{p}_{1}V)^{\perp})
\end{equation*}
with $P_{2}'$ invertible.  Using positive definite approximate inverse square
roots of $(P_{2}')^{*}P_{2}'$ and $(P_{2}'')^{*}P_{2}'',$ we can associate a KK-class
$\beta\in \Kf{\fambas}{\CO(\boubas),\CO(\fambas)}.$  By homotopy invariance,
$P_{1}$ gives the same K-class so $\beta=\qb'(\alpha).$  On the other hand,
the Kasparov module coming from $P_{2}'$ is degenerate, so $\beta$ can be 
defined using $P_{2}''.$ The K-class associated to the symbol of $P_{2}''$
is the families index of $\widehat{p}_{1}$ (cf. \cite{Luke}) which by the
families index of Atiyah-Singer should be precisely $\alpha\in
\Kc(T^{*}(\bourel)).$  This means that 
\begin{equation*}
                      \qb'(\alpha)=\beta= \qb(\alpha).
\end{equation*}
\end{proof}

Thus we can use indicial families instead of symbols to define the quantization
map $\qb.$  Consider again the manifold $\cn.$  Since the boundary of 
$\cn$ has two disconnected parts (which can themselves be disconnected), we
can consider, instead of $\fcK{\cnbm}(\cnfm),$ the group 
$\fcK{\Phi_{1}}(\cnfm)$ of stabilized homotopy classes of invertible
joint symbols with indicial family given by the identity at the boundary
face $\pa\famtot_{0}.$  If we define $\CO_{\boumap_{1}}(\cn)$ to be the
$C^{*}$ algebra of continuous functions
\begin{equation}
 \CO_{\Phi_{1}}(\cn)=\left\{f\in \CO(\cn);f\big|_{\pa\famtot_{1}}
=\boumap_{1}^{*}g \text{ for some }g\in\CO(\boubas_{1})\right\},
\label{gs.11}\end{equation}
then a quantization map
\begin{equation*}
  \qcn: \fcK{\Phi_{1}}(\cnfm)\longrightarrow \Kf{\fambas}{\CO_{\Phi_{1}}(\cn),\CO(\fambas)}
\end{equation*}
can be defined as follows.  Given a joint symbol $(\sigma, N)$ representing
a class in $\alpha\in \fcK{\Phi_{1}}(\cnfm),$ one considers a family of fibred
cusp operators 
$P$ of order zero with joint symbol given by $(\sigma,N).$  Deforming
$\sigma$ if needed, we can assume $P$ acts as the identity in a small collar
neighborhood of $\pa M_{0}.$  Then in the 
usual fashion, one can construct a Kasparov module in 
$\Kf{\fambas}{\CO_{\Phi_{1}}(\cn),\CO(\fambas)}$ by considering a positive
definite approximate inverse square root 
to P which acts as the identity in a collar 
neighborhood of $\pa M_{0}.$  As in the definition of $\quan,$ one can check
that the associate KK-class only depends on $\alpha\in\fcK{\Phi_{1}}(\cnfm).$

There is also a natural map
\begin{equation}
   \icn: \Kc(T^{*}(\bourel))\longrightarrow \fcK{\Phi_{1}}(\cnfm)
\label{gs.12}\end{equation}
which to $\alpha\in \Kc(T^{*}(\bourel))$ associates the 
class $\icn(\alpha)$ represented by any joint symbol with trivial indicial 
family on $\pa M_{0},$ indicial family with K-class given by $\alpha$ on
$\pa M_{1}$ and with trivial symbol.
\begin{lemma}
The diagram
\begin{equation*}
\xymatrix{
\Kc(T^*(\boubasrel))\ar[r]^{\icn}\ar[d]^{\qb}&
\fcK{\Phi_{1}}(\cnfm)\ar[d]^{\qcn} \\
\Kf{\fambas}{\CO(\boubas),\CO(\fambas)}\ar[r]^{\rcn^{*}}&
\Kf{\fambas}{\CO_{\Phi_{1}}(\cn),\CO(\fambas)}
} 
\end{equation*}
is commutative, 
where $\rcn: \CO_{\Phi_{1}}(\cn)\longrightarrow \CO(\boubas)$ is the restriction
to $\pa M_{1}.$
\label{gs.13}\end{lemma}
\begin{proof}
Let $\alpha\in \Kc(T^{*}\boubasrel)$ be given.
Let $e: \CO(D)\hookrightarrow \CO_{\Phi_{1}}(\cn)$ be the obvious injective
map of $C^{*}$ algebras so that $  \rcn\circ e=\Id.$ Then
\begin{equation}
       \rcn^{*}\qb(\alpha)= \rcn^{*}e^{*}\qcn\circ \icn(\alpha).
\label{gs.16}\end{equation}
Consider the $C^{*}$ algebra
\begin{equation*}
\CO_{1}(\cn)=\left\{f\in\CO(\cn);f\big|_{\pa\famtot_{1}}=0\right\}.
\end{equation*} 
From the commutativity of the diagram
\begin{equation}
\xymatrix{\fcK{\Phi_{1}}(\cnfm)\ar[r]^(.4){\sigma}\ar[d]^{\qcn}&
\Kc(T^*(\cn/B))\ar[d]^{\qr}\\
\Kf{\fambas}{\CO_{\boumap_{1}}(\cn),\CO(\fambas)}\ar[r]&
\Kf{\fambas}{\CO_{1}(\cn),\CO(\fambas)}}
\label{gs.14}\end{equation}
and the exactness in the middle of 
\begin{gather}
 \Kc(T^*(\boubasrel))\overset{\icn}{\longrightarrow} 
 \fcK{\Phi_{1}}(\cnfm)\overset{\sigma}{\longrightarrow} \Kc(T^*(\cn/B)), \\
\Kf{\fambas}{\CO(\boubas),\CO(\fambas)}\overset{\rcn}{\longrightarrow}
\Kf{\fambas}{\CO_{\boumap_{1}}(\cn),\CO(\fambas)}\longrightarrow 
\Kf{\fambas}{\CO_{1}(\cn),\CO(\fambas)},
\label{gs.15}\end{gather}
we deduce that there exists $\beta\in\Kf{\fambas}{\CO(\boubas),\CO(\fambas)}$
such that $\rcn^{*}\beta=\qcn\circ\icn(\alpha).$ Since
$\rcn\circ e=\Id,$
\begin{equation*}
            \beta= e^{*}\rcn^{*}\beta= e^{*} \qcn\circ\icn (\alpha),  
\end{equation*}
which implies by \eqref{gs.16} that
\begin{equation*}
           \qcn\circ \icn(\alpha) = \rcn^{*} \beta= \rcn^{*} (e^{*}
\qcn\circ\icn 
(\alpha))=\rcn^{*}\qb(\alpha).
\end{equation*}
\end{proof}

\begin{lemma} For $k\in\bbZ_{2},$ the diagram \eqref{gs.7} is commutative.
\label{gs.17}\end{lemma}

\begin{proof}
As noted above, we will only provide a proof for the case $k=0,$
the case $k=1$ being similar. Think of $\cn$ as a collar neighborhood
of $\famtot$ where $\pa\famtot_{1}$ is identified with $\pa\famtot.$  Let 
\begin{equation*}
          j: \CO_{\boumap}\famtot)\longrightarrow \CO_{\boumap_{1}}(\cn)
\end{equation*}
be the restriction map on $\cn\subset\famtot.$ Then consider the diagram
\begin{equation}
\xymatrix{
\Kc(T^*(\boubasrel))\ar[r]^{\icn}\ar[d]^{\qb}&
\fcK{\boumap_{1}}(\cnfm)\ar[d]^{\qcn}\ar[r]^{\ecn}& \fcK{\boumap}(\fammap)
\ar[d]^{\quan} \\
\Kf{\fambas}{\CO(\boubas),\CO(\fambas)}\ar[r]^{\rcn^{*}}&
\Kf{\fambas}{\CO_{\Phi_{1}}(\cn),\CO(\fambas)}\ar[r]^{j^{*}}& 
\Kf{\fambas}{\CO_{\boumap}\famtot),\CO(\fambas)}.}
\label{gs.18}\end{equation}
where $\ecn$ is the extension of a representative in $\fcK{\boumap_{1}}(\cnfm)$
by the identity inside $\famtot.$  Clearly, $\ecn\icn= i$ and $j^{*}\rcn^{*}=
\iota^{*},$  so the lemma will be proven provided we show the diagram
\eqref{gs.18} is commutative.  By the previous lemma, we only need to show 
that the right square is commutative. Consider the two graded Hilbert
modules over $\CO(\fambas)$
\begin{gather}
\Hv_{1}=\Ld(\cn;\sE)= \Ld(\cn;E_{+})\oplus \Ld(\cn;E_{-}) \\
\Hv_{2}=\Ld(\overline{\famtot\setminus \cn};\sE)=
\Ld(\overline{\famtot\setminus \cn};E_{+})\oplus
\Ld(\overline{M\setminus \cn};E_{-}).
\label{gs.19}\end{gather}
Then given $\alpha\in \fcK{\boumap_{1}}(\cnfm),$ we can represent
$j^{*}\qcn(\alpha)$ by a Kasparov module of the form $(\Hv_{1},
j^{*}\mu,\F_{1})$ and at the same time $\quan \ecn(\alpha)$ can be
represented by a Kasparov module of the form $(\Hv_{1}\oplus\Hv_{2},
j^{*}\mu\oplus l^{*}\nu,\F_{1}\oplus \F_{2})$ where 
\begin{equation*}
           l:\CO_{\boumap}\famtot)\longrightarrow \CO(\overline{\famtot\setminus \cn}) 
\end{equation*}
is the restriction map,
\begin{gather}
\mu: \CO_{\Phi_{1}}(\cn)\longrightarrow \mathcal{B}(\Hv_{1}) \\
\nu: \CO(\overline{M\setminus \cn})\longleftrightarrow \mathcal{B}(\Hv_{2}) 
\label{gs.20}\end{gather}
are the obvious actions given by multiplication and
\begin{equation*}
             \F_{2}=\left( \begin{array}{cc}
                      0 & \Id \\
                    \Id & 0
                      \end{array} \right) .
\end{equation*}  
Since $(\Hv_{2},\nu, \F_{2})$ is a degenerate Kasparov module,
\begin{equation*}
j^{*}\qcn(\alpha)= \quan \ecn(\alpha).
\end{equation*}
\end{proof}

Finally we need to show that the diagram
\begin{equation}
\xymatrix{
\Kcn{k}(T^*(\famrel))\ar[d]^{\qr}\ar[r]^{I}& 
\Kcn{k-1}(T^*(\boubasrel))\ar[d]^{\qb} \\
\Kfn{k}{\fambas}{\COo(\famtot),\CO(\fambas)}\ar[r]^{\delta}& 
\Kfn{1-k}{\fambas}{\CO(\boubas),\CO(\fambas)}}
\label{gs.21}\end{equation}
is commutative for $k\in \bbZ_{2},$ which requires an understanding of
the boundary homomorphism $\delta.$  Let us again limit our attention to the
even case $k=0,$ the case $k=1$ being similar. Recall that in topological
K-theory (see for instance \cite{Atiyah1}), the
boundary homomorphism of the 6-term exact sequence associated to a pair of 
spaces $(X,Y)$ is defined via the cone space $X\cup CY$ which is obtained 
from $X$ and 
\begin{equation*}
CY= Y\times [0,1]/Y\times\{1\}
\end{equation*}
by identifying $Y\subset X$ with $Y\times\{0\}\subset CY.$ There is an
exact sequence
\begin{equation*}
\Kt(X\cup CY, X)\overset{m^{*}}{\longrightarrow}
\tKt(X\cup CY) \overset{k^{*}}\longrightarrow \tKt(X).
\end{equation*}
Under the identifications 
\begin{equation*}
   \tKt^{-1}(Y)\cong\Kt(X\cup CY, X),\
\Kt(X,Y)\cong \tKt(X\cup CY),
\end{equation*}
this becomes
\begin{equation*}
\tKt^{-1}(Y)
\overset{\delta}{\longrightarrow}
  \Kt(X,Y)   \longrightarrow \tKt(X)
\end{equation*}
where $\delta$ is the boundary homomorphism of the 6-term exact sequence 
associated to the pair $(X,Y).$

In $\KK$-theory, one can define the boundary homorphism in a 
similar way, introducing the mapping cone
\begin{equation}
  \mc=\left\{(x,f) \in \CO_{\boumap}(\famtot)\oplus
\COo([0,1)\times\boubas);\ \iota(x)=f(0)\right\},
\label{gs.22}\end{equation}    
where 
\begin{equation*}
  \COo([0,1)\times\boubas)=\left\{ f\in\CO([0,1]\times\boubas);
  f\big|_{\{1\}\times\boubas}=0\right\}
\end{equation*}
is the $C^{*}$ closure of $\CO_{c}([0,1)\times\boubas).$ There is a natural
  inclusion
\begin{equation}
e:\COo(\famtot)\ni x \longmapsto (x,0)\in \mc
\label{gs.23}\end{equation}
which induces a map
\begin{equation}
   e^{*}: \Kf{\fambas}{\mc, \CO(\fambas)} \longrightarrow 
\Kf{\fambas}{\COo(\famtot), \CO(\fambas)}.
\label{gs.24}\end{equation}
This map is an isomorphism as can be seen by interpreting the map $e^*$ as
multiplication (on the left) by $[e]\in \Kf{\fambas}{\COo(\famtot),\mc}$
associated to the $C^{*}$ homomorphism $e$ (see\cite{Blackadar1}). 
Then $[e]$ is a $\KK$-equivalence between $\COo(\famtot)$ and $\mc,$
which is to say that it has an inverse $u\in\Kf{\fambas}{\mc,\COo(\famtot)}$
with respect to the Kasparov product. 

Given this one can then define the boundary homomorphism as
\begin{equation}
  \delta= j^{*}(e^{*})^{-1}: \KK_{\fambas}(\COo(\famtot),\CO(\fambas))
\longrightarrow \KK_{\fambas}(\COS(\boubas),\CO(\fambas))\cong \KK^{1}_{\fambas}
(\CO(\boubas),\CO(\fambas))
\label{gs.25}\end{equation} 
where $j: \COS(D)\hookrightarrow \mc$ is the natural inclusion.  In terms of
Kasparov product, $\delta$ is multiplication on the left by $j^{*}u\in
\KK_{\fambas}(\COS(\boubas),\COo(\famtot)).$ 

In the present context, it is possible to give an alternative definition of the 
mapping cone which is especially useful. Consider as before a collar 
neighborhood $\cn=\pa\famtot\times[0,1]$ of $\pa\famtot$ with $\pa\famtot$
identifies with $\pa\famtot_{1}=\pa\famtot\times\{1\}.$ Then
\begin{equation}
   \mc \cong\left\{f\in \COo(\famtot);f\big|_{\cn}= (\boumap\times\Id)^{*}g
   \text{ for some }g\in \COo(D\times [0,1))\right\}
\label{gs.26}\end{equation}
by identifying $\bfamtot= \overline{\famtot/\cn}$ with $\famtot.$ The 
isomorphism \eqref{gs.26} gives
\begin{equation}
\varepsilon: \mc\hookrightarrow \COo(\famtot) .
\label{gs.27}\end{equation}
The map $e$ in \eqref{gs.23} is then described by
\begin{equation*}
 e: \COo(\famtot)\tilde{\longrightarrow}
\COo(\bfamtot)\hookrightarrow \mc.
\end{equation*}

\begin{lemma}
The boundary homomorphism $\delta$ of the 6-term exact sequence
\eqref{gs.2} is given by the pull-back map 
\begin{equation*}
 j^{*}\varepsilon^{*}: \Kf{\fambas}{\COo(\famtot),\CO(\fambas)}
\longrightarrow \Kf{\fambas}{\COS(\boubas),\CO(\fambas)}.
\end{equation*}
\label{gs.28}\end{lemma}

\begin{proof}
If $\sigma$ is a symbol on $T^{*}(\famrel)$ then $e^{*}\varepsilon^{*}
\qr(\sigma)$ can be obtained by quantizing the symbol 
\begin{equation*}
          \sigma'=\sigma\big|_{T^{*}(\bfamrel)}
\end{equation*}
on $T^{*}(\bfamrel)$ and making the identification $\COo(\bfamtot)\cong
\COo(\famtot).$  Since $\sigma'$ and $\sigma$ are homotopic when $\bfamtot$
and $\famtot$ are identified, they correspond to the same K-class.  Since
$\sigma$ was arbitrary, this shows $\qr= e^{*}\varepsilon^{*}\qr.$
Consequently,
\begin{equation*}
       (e^{*})^{-1}\qr= (e^{-*})^{-1}e^{*}\varepsilon^{*}\qr= \varepsilon^{*}
\qr.
\end{equation*}
Since $\qr$ is an isomorphism, $(e^{*})^{-1}= \varepsilon^{*}$ and it follows
that $\delta= j^{*}(e^{*})^{-1}=j^{*}\varepsilon^{*}.$
\end{proof}

After identifying the tangent bundle with the cotangent bundle via some
metric, The boundary fibration induces a 
fibration
\begin{equation}
  \cnf: T^{*}(\cn/\fambas)\longrightarrow T^{*}(\boubas\times [0,1]/\fambas).  
\label{gs.29}\end{equation}
Let 
\begin{equation}
  \ind_{\cnf}: \Kc(T^{*}(\cn/\fambas))\longrightarrow
  \Kc(T^{*}(\boubas\times[0,1]/\fambas))  
\label{gs.30}\end{equation}
be the topological index family map associated to this fibration.

\begin{lemma} The diagram
\begin{equation*}
\xymatrix{
  \Kc(T^{*}(\famrel))\ar[r]^{r^{*}}\ar[d]^{\qr}&
   \Kc(T^{*}(\cn/\fambas))\ar[r]^{\ind_{\cnf}}\ar[d]^{\qr}&
  \Kc(T^{*}(\boubas\times[0,1]/\fambas))\ar[d]^{\qr} \\
\KK_{\fambas}(\COo(\famtot),\CO(\fambas))\ar[r]^{\ccn^{*}}&
\KK_{\fambas}(\COo(\cn),\CO(\fambas))\ar[r]^{\dcn^{*}}&
\KK_{\fambas}(\COS(\boubas),\CO(\fambas))
}
\end{equation*}
is commutative, where
\begin{equation*}
  \ccn:\COo(\cn)\hookrightarrow \COo(\famtot), \quad \dcn:\COS(\boubas)
\hookrightarrow \COo(\cn), \quad r: T^{*}(\cn/\fambas))\hookrightarrow
T^{*}(\famrel),        
\end{equation*}
are the natural inclusions.
\label{gs.31}\end{lemma}

\begin{proof}
The commutativity of the first square is clear.  The proof of the 
commutativity of the second square is essentially a consequence of the 
Atiyah-Singer families index theorem.  Indeed, it allows us to define 
$\ind_{\cnf}$ as an analytical families index.  Then, using quantization of 
symbols over Hilbert bundles as in \cite{Luke}, one can represent
(cf. the proof of Lemma~\ref{gs.10}) $\dcn^{*}\qr(\alpha)$ by a 
Kasparov module of the form
\begin{equation*}
          (\Hv_{1}\oplus\Hv_{2}, \mu, \F_{1}\oplus \F_{2})
\end{equation*}  
where 
\begin{equation*}
\begin{gathered}
\Hv_{1}=\Ld(\boubas\times[0,1]; V_{+})\oplus\Ld(\boubas\times[0,1]; 
V_{-}), \\
 \Hv_{2}=\Ld(\boubas\times[0,1]; V_{+}^{\perp})\oplus\Ld(\boubas\times[0,1]; 
V_{-}^{\perp}), 
\end{gathered}
\end{equation*}
and $V_{\pm}$ is a sub-bundle of $\Ld(\cn/(D\times[0,1]);E_{\pm})$ with
$V_{+}$ and $V_{-}$ of same finite corank.  As usual, 
\begin{equation*}
    \mu: \COo(D\times [0,1]) \longrightarrow \mathcal{B}(\Hv_{i})
\end{equation*}
denotes multiplication.  One can do this in such a way that $(\Hv_{1},\mu,
\F_{1})$ is degenerate and $(\Hv_{2},\mu,\F_{2})$ represents 
$\qr\circ\ind_{\cnf}(\alpha),$ which establishes the commutativity of the 
left square. 
\end{proof}

\begin{lemma}
We have the identity $\delta \qr= \qr \kappa I_{0}$ where
\begin{equation*}
         \kappa: \Kcn{-1}(T^{*}(\boubasrel))\longrightarrow
	 \Kc(T^{*}(\boubas\times[0,1]/\fambas)) 
\end{equation*}
is the canonical isomorphism induced from the homotopy equivalence
\begin{equation*}
(T^{*}([0,1]),\infty)\sim (\bbR,\infty).
\end{equation*}
\label{gs.32}\end{lemma}
\begin{proof}
Since clearly $\dcn^{*}\ccn^{*}=j^{*}\varepsilon^{*}=\delta,$ we deduce
from Lemma~\ref{gs.31} that 
\begin{equation*}
    \delta \qr =\dcn^{*}\ccn^{*}\qr = \qr \ind_{\cnf}r^{*}, 
\end{equation*}
so we can conclude from the fact that $\ind_{\cnf}\circ r^{*}= \kappa I_{0},$
which is a consequence of the commutativity of the topological index map
with the Thom isomorphism. 
\end{proof}

Thus, Lemma~\ref{gs.32} reduces the proof of the commutativity of 
\eqref{gs.21} in the even case to the following statement.
\begin{lemma}
The diagram
\begin{equation*}
\xymatrix{
  \Kcn{-1}(T^{*}(\boubasrel))\ar[r]^{\kappa}\ar[d]^{\qb}&
 \Kc(T^{*}(\boubas\times[0,1]/\fambas))\ar[d]^{\qr} \\
\KK_{\fambas}(\CO(\boubas),\COS(\fambas))\ar[r] & \KK_{\fambas}(
\COS(\boubas),\CO(\fambas))}
\end{equation*}
is commutative, where the bottom arrow is Bott periodicity in KK-theory.
\label{gs.33}\end{lemma}
\begin{proof}
This essentially follows from the family version of proposition 11.2.5 in
\cite{Higson-Roe} and the alternative definition of odd quantization using
self-adjoint operators.
\end{proof}

In the odd case, a similar argument reduces the proof of the commutativity of 
\eqref{gs.21} to the following lemma.
\begin{lemma}
The diagram 
\begin{equation*}
\xymatrix{
  \Kc(T^{*}(\boubasrel))\ar[r]\ar[d]^{\qb}&
 \Kc(T^{*}(\boubas\times(0,1)\times[0,1])/[\fambas\times(0,1)] ))\ar[d]^{\qr}\\
\Kf{\fambas}{\CO(\boubas),\CO(\fambas)}\ar[r] &
\Kf{\fambas}{\COS(\boubas),\COS(\fambas)}}
\end{equation*}
is commutative, where the top and bottom arrows are Bott periodicity in 
K-theory and KK-theory respectively.
\label{gs.34}\end{lemma}
\begin{proof}
This essentially follows from the family version of proposition 11.2.5 in
\cite{Higson-Roe} and the alternative definition of odd quantization using
self-adjoint operators.
\end{proof}

\section{Adiabatic passage to the cusp case}\label{Adiabaticpassage}

As noted in the introduction, and confirmed by the properties of the 6-term
exact sequence above, it is the `cusp' case amongst the possible fibred
cusp structures on a give fibration which is universal.

\begin{proposition}\label{faficu.215} For any fibration with fibred cusp
  structure there is a natural map, given by an adiabatic limit,
  $q_{\ad}:\fcK{\boumap}(\fammap)\longrightarrow \cuK{\fammap}$ into the cusp
  K-group which commutes with the analytic index and gives a commutative
  diagram 
\begin{equation}
\xymatrix{
&
\fcK{\boumap}(\fammap)\ar[d]^{q_{\ad}}\ar[dl]_{\inda}\ar[dr]^{\sigma}&\\
\Kt(\fambas)&
\cuK{\fammap}\ar[r]_(.4){\sigma}\ar[l]^{\inda}&
\Kc(T^*(\famrel)).}
\label{faficu.216}\end{equation}
\end{proposition}

\begin{proof} Given a family $P\in\fcP{\boumap}0{\famtot;\sE}$ we
  construct an adiabatic family
  $P(\epsilon)\in\adcuP{\boumap}0{\famtot;\sE}$ which are cusp
  pseudodifferential operators for $\epsilon >0$ but degenerate, in the
  adiabatic limit, to the fibred cusp operator $P$ at $\epsilon =0.$ However, it
  is perhaps better to think of the construction as being in the opposite
  direction. The definition and properties of such adiabatic fibred cusp
  operators are discussed in Appendix~\ref{Adiabatic} where they are
  defined directly through their Schwartz kernels which are defined on the
  space given by blow-ups in \eqref{faficu.185}, in our case with the finer
  fibration of the boundary being the given one, $\boumap,$ and the coarser
  fibration being $\pa\fammap$ corresponding to the cusp structure.

Thus there are three boundary faces of the double space in
\eqref{faficu.185} which meet the diagonal (and the kernels are supposed to
vanish rapidly at all other boundary faces). At the `old' face $\epsilon
=0$ we wish to recover the given operator $P$ by the map
\eqref{faficu.191}; this fixes the kernel precisely on that face. The
adiabatic front face can be constructed, as in \eqref{faficu.185}, by the
last blow-up and it is fibred over the lifted variable 
\begin{equation}
\tau =\frac{x-\epsilon}{x+\epsilon}\in[-1,1]
\label{faficu.217}\end{equation}
representing the adiabatic passage, with fibre which is just
the front face for the fibred cusp calculus. On the fibre at $\tau=1$ the
kernel is already fixed to be that of $P.$ Thus, we may simply choose to
extend the kernel to the whole adiabatic face to be independent of the
`angle' $\tau.$ This fixes the kernel of the adiabatic normal
operator and also fixes the boundary value, corresponding to $\tau=-1,$ of
the kernel on the front face for the cusp blow-up as $\epsilon \downarrow0.$
Simply extend it to that face as a conormal distribution with respect to the
diagonal; in fact we also choose an extension which is conormal to the
diagonal in the interior and consistent with these boundary values (and
also the requirement of the vanishing at the `non-diagonal' boundaries).

This fixes an element $P(\epsilon)\in\adcuP{\boumap}0{\famtot;\sE}.$ The
choices ensure that in the sense of \eqref{faficu.191} 
\begin{equation}
A(P(\epsilon ))=P\Mand \ad(P)^{-1}\text{ exists,}
\label{faficu.218}\end{equation}
since the invertibility here comes from the invertibility of the normal
operator of $P.$ Now, for $\epsilon \le\delta,$ for some $\delta >0,$ the
hypotheses of Proposition~\ref{faficu.192} apply and show that $P(\delta
)\in\cuP0{\famtot;\sE}$ is fully elliptic. Thus we can define 
\begin{equation}
q_{\ad}:\fcK{\boumap}(\fammap)\ni[P]\longmapsto [P(\delta )]\in\cuK{\fammap}
\label{faficu.219}\end{equation}
provided the independence of choices is shown.

Clearly the image in \eqref{faficu.219} is independent of the choice of
$\delta$ small and the choices in the construction can all be related
by homotopies. Similarly an homotopy of $P$ through fully elliptic fibred
cusp operators lifts to a homotopy of $P(\delta)$ in fully elliptic cusp
operators and the behaviour under stabilization and composition with bundle
isomorphism is appropriate to guarantee that \eqref{faficu.219} is well
defined. The construction also ensures that this adiabatic limit commutes
with the passage to the symbol, giving commutativity in \eqref{faficu.216}.
\end{proof}

Under the Poincar\'{e} duality of Section~\ref{GeneralPoincare}, the 
corresponding adiabatic map in KK-theory is given by pull-back. Let
$\iota_{\ad}:\Ccu(\famtot)\hookrightarrow \Cfc(\famtot)$ be the natural
inclusion.  

\begin{proposition}\label{20.7.2005.3}
The diagram 
\begin{equation*}
  \xymatrix{
 \fcK{\boumap}(\fammap)\ar[r]^{q_{\ad}}\ar[d]^{\quan}&
\cuK{\fammap}\ar[d]^{\quan} \\
\Kf{B}{\Cfc(\famtot),\CO(\fambas)}\ar[r]^{\iota_{\ad}^{*}} &  
\Kf{B}{\Ccu(\famtot),\CO(\fambas)}     }
\label{faficu.223}\end{equation*}
is commutative.
\end{proposition}

\begin{proof}
If $s: \COo(\famtot)\hookrightarrow \Ccu(\famtot)$ and 
$i:\CO(B)\hookrightarrow \Ccu(\famtot)$ are the natural inclusions,
 then from the fact that 
$q_{\ad}$ preserves the index  and the homotopy class of the symbol, we see
that  
\begin{equation*}
            i^{*}\iota_{\ad}^{*}\quan = i^{*}\quan q_{\ad} \quad \mbox{and}
\quad s^{*}\iota_{\ad}^{*}\quan = s^{*}\quan q_{\ad}.  
\label{faficu.226}\end{equation*}
Thus, the result follows from the split short exact sequence of 
Lemma~\ref{aco.8}.
\end{proof}
\section{The extension of fibred cusp structures}\label{efc.0}

The definition of the topological index uses an embedding of a given family
of cusp structures into a product setting. In fact we show how to embed a
general fibred cusp structure (since this can also be used to define the
index without using Proposition~\ref{faficu.215}). We shall use as `model'
structure (for a given base manifold $\fambas)$ the products
\begin{equation}
\begin{gathered}
\efamtot=\mfamtot,\ p>0,\text{ for cusp structures and}\\
\efamtot=\Mfamtot,\ p,q>0
\text{ for fibred cusp structures.}
\end{gathered}
\label{faficu.116}\end{equation}
The model fibration, being just projection onto the first factor will be
written
\begin{equation}
\mfammap:\efamtot\longrightarrow \fambas
\label{faficu.117}\end{equation}
and in the case of fibred cusp structures the model boundary fibration is
the projection onto the first two factors (restricted to the boundary) 
\begin{equation}
\Mboumap:\pa\efamtot\longrightarrow\eboubas=\Mboubas.
\label{faficu.118}\end{equation}

\begin{definition}\label{faficu.149} A fibration with fibred cusp structure
  $\efammap:\efamtot\longrightarrow \fambas$ is an \emph{extension} of
  a given fibration with fibred cusp structure
$\fammap:\famtot\longrightarrow \fambas$ if
  $i: \famtot\hookrightarrow \efamtot$ and $j: \boubas\hookrightarrow \eboubas$
  have the following properties
\begin{enumerate}
\item\label{e.one}
$i$ embeds $\famtot$ as a `product type' submanifold in the sense that it
  is an embedding and 
\begin{equation}
i(\pa\famtot)=\pa\efamtot\cap i(\famtot)
\label{faficu.120}\end{equation}
and the pull-back under $i$ of a defining function for $\pa\efamtot$ is a
defining function for $\pa\famtot.$
\item\label{e.two} $\fammap=\Mfammap\circ i.$ 
\item\label{e.three} $j\circ\boumap=\Mboumap\circ i.$
\end{enumerate}
Notice that the inclusion $j$ is completly specified by the inclusion $i$.
\label{extension}\end{definition}

\begin{proposition}\label{faficu.114} For any family of fibred cusp
structures as in \eqref{faficu.2} there is an extension  
\begin{equation}
i:\famtot\longrightarrow \efamtot=\Mfamtot\label{faficu.115}\end{equation}
provided that $p$ and $q$ are large enough. 
In the case of cusp structures there is an extension
$i:\famtot\longrightarrow\mfamtot,$ for $p$ sufficiently large.
\end{proposition}

\begin{proof} First consider the cusp case. By Whitney's Embedding Theorem,
  we may embed $\famtot$ in Euclidean space of sufficiently large dimension 
\begin{equation}
i':\famtot\longrightarrow \bbR^k.
\label{faficu.122}\end{equation}
To replace this by an embedding of the desired product type into the
interior of a closed ball, consider the embedding
\begin{equation}
(x,i'):\famtot\longrightarrow [0,1)\times\bbR^k,
\label{faficu.123}\end{equation}
where $x\in\CI\famtot)$ is the boundary defining function (it may be assumed
without loss of generality that $0\le x\le 1.)$ But $[0,1)\times\bbR^k$ can
be identified, using a smooth map $e,$ with an open neighbourhood of a
piece of the boundary of $\bbB^{p+1},$ $p=k.$ It follows directly that this
is an embedding of product type in the sense of \eqref{e.one}. Then let
\begin{equation}
i=(\fammap,e(x,i')):\famtot\longrightarrow \mfamtot
\label{faficu.125}\end{equation}
be the product of this map with $\famtot$ as a map into
$\fambas\times\bbB^{p+1}.$ This is still an embedding with property
\eqref{e.one} and has the property \eqref{e.two} (and so \eqref{e.three}
as well).

In the case of a fibred cusp structure, we can start with the constructions
above for the embedding of the underlying cusp structure. Now, take an
additional embedding of $\boubas$ into a Euclidean space
\begin{equation*}
 d:\boubas\hookrightarrow \bbR^{k}
\end{equation*} 
for convenience again $\bbR^k$ by increasing $k$ if
necessary. Using a product structure near the boundary this
fibration can be extended inwards near the boundary to a fibration over
$[0,\epsilon)_{x}\times\boubas$ and we can consider the smooth map 
\begin{equation}
(\Id ,d\circ \boumap):[0,\epsilon)_{x}\times \pa M
\longrightarrow[0,\epsilon)\times\bbR^k  .
\label{faficu.121}\end{equation}
Taking $0<\epsilon<1$ small enough and using radial retraction on $\bbR^{k},$
this can be extended to a smooth map
\begin{equation*}
       (x,d'): M\longrightarrow [0,1)\times \bbR^{k}.
\end{equation*}
For the overall embedding we can then take 
\begin{equation}
i=(\fammap,e'(x,d'),e''\circ i'):\famtot\longrightarrow \Mfamtot, \quad p=q=k,
\label{faficu.124}\end{equation}
with
\begin{equation*}
\begin{gathered}
e':[0,1)\times\bbR^k\times\bbR^{k}\hookrightarrow\bbB^{p+1},\
p=k,\ q=k  \\
e": \bbR^{k}\hookrightarrow \bbS^{q},
\end{gathered}
\end{equation*}
being identifications of open sets.  We can also take 
$j:D\hookrightarrow \fambas\times \bbS^{p}$ to be 
\[
    j= (\fammap, e''(0,d')).
\]
Since $i'$ is an embedding, this is
also and satisfies \eqref{e.one} -- only the boundary is mapped into the
boundary, since $x$ is a boundary defining function of $\famtot.$ Now,
\eqref{e.two} holds for the same reason as before and \eqref{e.three}
follows from the definition of $j.$
\end{proof}

\section{Multiplicativity}\label{Multiplicativity}

Recall one form of the `multiplicativity property' for the index of
pseudodifferential operators as introduced by Atiyah and Singer. Consider an
iterated fibration, for the moment of compact manifolds without boundary
but later allowing $\famfib$ to have a boundary
\begin{equation}
\xymatrix{\sfamfib\ar@{-}[r]&\bfamtot\ar[d]^{\tfammap}&\pa\sfamfib=\emptyset\\
\famfib\ar@{-}[r]&\famtot\ar[d]^{\fammap}\\
&\fambas.}
\label{faficu.135}\end{equation}
Taking a connection on $\tfammap$ allows the fibre cotangent bundle
$T^*(\tfamrel)$ to be identified (naturally up to homotopy) as a subbundle
of $T^*(\bfamrel)$ which then splits
\begin{equation}
T^*(\bfamrel)=T^*(\tfamrel)\oplus(\tfammap)^*T^*(\famrel).
\label{29.4.2005.3}\end{equation}
Choose cutoff functions $\chi_1$ and $\chi_2$ which are homogeneous of
degree $0,$ smooth outside the zero section and are respectively supported
outside the two summands and equal to one in a conic neighbourhood of the
other with $\chi_1+\chi_2=1.$
 
Let $A\in\Ps0{\famrel;\sE}$ and $B\in\Ps0{\tfamrel;\sG}$ be families of
elliptic pseudodifferential, then the matrix
\begin{equation}
\begin{pmatrix}
\chi_1\sigma (B)& -\chi_2\sigma (A)^*\\
\chi_2\sigma (A)&\chi_1\sigma (B)
\end{pmatrix}
\label{29.4.2005.2}\end{equation}
acting from sections of $(\sE\otimes\sG)_+ = 
E_+\otimes G_+\oplus E_-\otimes G_-$ to
$(\sE\otimes\sG)_-=E_-\otimes G_+\oplus E_+\otimes G_-$ is a family of
elliptic symbols for the overall fibration $\bfamtot\longrightarrow\fambas.$

\begin{proposition}\label{29.4.2005.1} (Atiyah and Singer) If the family $B$
 has trivial one-dimensional index in $\Kt(\famtot)$ then any operator in
 $\Ps0{\bfamrel;\sE\otimes\sG}$ with symbol \eqref{29.4.2005.2} has
 index in $\Kt(\fambas)$ equal to that of $A.$
\end{proposition}

\begin{proof} 
To prove this we construct a natural `product type' calculus for the top
fibration which includes the fibrewise pseudodifferential operators and an
appropriate class of lifts of operators on the fibres of the lower
fibration and then do a deformation to the `true' pseudodifferential
calculus.

To ease the notational burden, at least initially, we consider the case of
the numerical index and suppose that the second fibration, $\phi$ has just
one fibre. The general case will be proved as part of the extension below
to the cusp calculus. The product-type algebra is discussed in
Appendix~\ref{Appendix-PT}. It contains the algebra of
pseudodifferential operators on the total space, $\bfamtot$ but has a
double order filtration and we denote the filtered spaces
$\ptP{\tfammap}{m,m'}{\bfamtot;\sH}$ for a $\bbZ_2$-graded vector bundle
$\sH=(H_+,H_-)$ over $\bfamtot;$   
\begin{equation}
\Ps m{\bfamtot;\sH}\subset\ptP{\tfammap}{m,m}{\bfamtot,\sH}\ \forall\ m \in\bbZ
\label{29.4.2005.4}\end{equation}
(there is no particular problem with real order, but since we do not need
them, we shall not bother with them here.)

The basic properties of these operator are similar to those of regular
pseudodifferential operators 
\begin{equation}
\begin{gathered}
\ptP{\tfammap}{m_2,m'_2}{\bfamtot,\sH_1}\circ
\ptP{\tfammap}{m_1,m'_1}{\bfamtot,\sH_2}
\subset\ptP{\tfammap}{m_1+m_2,m'_1+m'_2}{\bfamtot,\sH},\\
\text{provided }(\sH_1)_+=(\sH_2)_{-},\ \sH_+=(\sH_2)_+,\ \sH_-=(\sH_1)_-.
\end{gathered}
\label{29.4.2005.5}\end{equation}
The main difference between this product-type algebra and the usual one is
that there are two related symbol maps. The `usual' principal symbol map
is modified to a short exact sequence 
\begin{equation}
\begin{gathered}
\xymatrix@1{
\ptP{\tfammap}{m-1,m'}{\bfamtot,\sH}\ar[r]&
\ptP{\tfammap}{m,m'}{\bfamtot,\sH}\ar[r]^(.5){\sigma _{m,m'}}&
  \ptPS{\tfammap}{m,m'}{\bfamtot,\sH} } , \\
  \ptPS{\tfammap}{m,m'}{\bfamtot,\sH}=
\CI([S^*\bfamtot,\tfammap^*S^*(\famtot)];
\hom(\sH)\otimes_+N^{-m}\otimes_+N_{\ff}^{-m'}).
\end{gathered}
\label{29.4.2005.6}\end{equation}

Here, $S^*\bfamtot$ is the sphere bundle of $T^*\bfamtot$ (best thought of as the
boundary of the radial compactification of this bundle) which is then blown
up at the submanifold given by the corresponding sphere bundle at infinity
of the cotangent bundle of the base. Thus
$[S^*\bfamtot,\tfammap^*S^*(\famtot)]$ is the `old' boundary hypersurface
of the manifold with corners $[\com{T^*\bfamtot},\tfammap^*S^*(\famtot)].$
The bundle $N$ is the normal bundle in this sense and $N_{\ff}$ is the
normal bundle to the front face, thought of as a trivial bundle over
$[S^*\bfamtot,\tfammap^*S^*\famtot].$ Thus these factors just 
represent the growth order of symbols at the boundaries. This `standard'
symbol is multiplicative in the obvious sense that 
\begin{equation}
\sigma _{m_1+m_2,m'_1+m'_2}(BA)=\sigma _{m_2,m'_2}(B)\sigma _{m_1,m'_1}(A)
\label{29.4.2005.7}\end{equation}

More interestingly, there is a non-commutative symbol, which has much in
common with the indicial family for the fibred cusp calculus; it is called here
the \emph{base family.} Since $S^*\famtot\longrightarrow\famtot$ is
a fibration we can lift the fibres of
$\tfammap:\bfamtot\longrightarrow\famtot$ to give a fibration
$S^*_{\tfammap}\famtot\longrightarrow S^*\famtot$ which has the same
typical fibre as $\tfammap.$ This second symbol gives a short exact sequence 
\begin{equation}
\ptP{\tfammap}{m,m'-1}{\bfamtot,\sH}\hookrightarrow
\ptP{\tfammap}{m,m'}{\bfamtot,\sH}\overset{L_{m,m'}}\longrightarrow
\Ps{m'}{S^*_{\tfammap}\famtot/S^*\famtot;\sH\otimes_+ N^{-m}}.
\label{29.4.2005.8}\end{equation}
Here $N$ is the inverse homogeneity bundle on $T^*\famtot,$ which is to say the
normal bundle to the boundary of the radial compactification 
$\com{T^*\famtot}.$ Again this sequence is multiplicative in the sense that 
\begin{equation}
L_{m_1+m_2,m'_1+m'_2}(BA)=L_{m_2,m'_2}(B)L_{m_1,m'_1}(A)
\label{29.4.2005.9}\end{equation}
under \eqref{29.4.2005.5}.

Although these two maps are separately surjective they are related through
the symbol map on the image of \eqref{29.4.2005.8}. Namely the combination
of the two maps gives a short exact sequence
\begin{multline}
\xymatrix@1{\ptP{\tfammap}{m-1,m'-1}{\bfamtot,\sH} \, \ar@{^(->}[r]&
\ptP{\tfammap}{m,m'}{\bfamtot,\sH}\ar[r]^{(\sigma,L)}&
\ptjS{\tfammap}{m,m'}{\bfamtot,\sH}}\\
\ptjS{\tfammap}{m,m'}{\bfamtot,\sH}=\big\{(a,A)\in
\Ps{m'}{S^*_\psi X/S^*X;\sH\otimes_+ N^{-m}}\oplus 
\ptPS{\tfammap}{m,m'}{\bfamtot,\sH}; \\
\sigma _{m'}(A)=a\big|_{\pa [S^*\bfamtot,\psi^*S^*(\famtot)]}\big\}.
\label{29.4.2005.11}\end{multline}

There are appropriate Sobolev spaces on which these operators are
bounded. Since we are interested principally in the case of operators of
order $0$ it suffices to note that 
\begin{multline}
\ptP{\tfammap}{0,0}{\bfamtot;\sH}\ni A:L^2(\bfamtot;H_+)\longrightarrow L^2(\bfamtot;H_-)\\
\text{ is bounded and is compact iff }\\
\sigma _{0,0}(A)=0,\ N_{0,0}(A)=0.
\label{29.4.2005.10}\end{multline}
From this and standard constructions it follows that $A$ in
\eqref{29.4.2005.10} is Fredholm iff
\begin{equation}
\exists\
\sigma _{0,0}(A)^{-1}\in \ptPS{\tfammap}{0,0}{\bfamtot,\sH^{-}}
\Mand N_{0,0}(A)^{-1}\in\Ps{m'}{S^*_{\tfammap}\bfamtot/S^*\famtot;\sH^-}.
\label{29.4.2005.12}\end{equation}

As well as \eqref{29.4.2005.4} the fibrewise pseudodifferential operators
may also be considered as product-type operators
\begin{equation}
\begin{gathered}
\Ps m{\bfamtot/\famtot;\sH}
\subset\ptP{\tfammap}{0,m}{\bfamtot;\sH},\\ 
\sigma _{0,m}(A)=\sigma _m(A), \
N_0(A)=A\ \forall\ A\in\Ps{m}{\bfamrel;\sH}.
\end{gathered}
\label{29.4.2005.13}\end{equation}
In particular such an operator of order $0$ is Fredholm on $L^2$ if and
only it is invertible, at least if the fibration has non-trivial base.

It is also important to see that pseudodifferential operators on the base
are in a sense included in the product type pseudodifferential operators on
the total space. Suppose that $E$ is a $\bbZ_2$-graded vector bundle over
$\famtot$ which is embedded in $\CI(\bfamtot; \sF)$ for some
$\bbZ_2$-bundle $\sF.$ Thus there is a pair of families of finite rank,
smoothing, projections $\pi_\pm\in\Ps{-\infty}{\tfamrel;F_\pm}$ projecting
onto the fibre of $E_\pm$ over each point of $\famtot.$ Then suppose that
$A\in\Ps{m}{\famtot;\sE}$ is some pseudodifferential operator. With the
identification of bundles, we may identify $A$ as an operator from
$\CI(\bfamtot;F_+)$ to $\CI(\bfamtot;F_-)$ and then
\begin{equation}
A=\pi_-A\pi_+\in\ptP{\tfammap}{m,-\infty}{\bfamtot;\sF}.
\label{29.4.2005.14}\end{equation}
In this case  
\begin{equation}
N_{m}(A)=\pi_-\sigma _m(A)\pi_+\in
\Ps{-\infty}{S^*_{\tfammap}\bfamtot/S^*\famtot;\sF\otimes_+N_{-m}}.
\label{29.4.2005.15}\end{equation}

Under the hypothesis of the Proposition that we are trying to prove -- for
the moment only in the case of a single operator -- we have a fibre family of
trivial index of rank one. So, by smoothing perturbation we may assume that
this $B\in\Ps0{\tfamrel;\sG}$ is surjective and has a null bundle which is
a trivial line bundle. Taking the $\bbZ_2$-bundle from the base we may
form the extended operator which we denote 
\begin{multline}
\begin{pmatrix}B&0\\0&B^*
\end{pmatrix}\in\Ps0{\tfamrel;\sE\otimes\sG}
\subset\ptP{\tfammap}{0,0}{\bfamtot;\sE\otimes\sG},\\
\sE\otimes\sH=(E_+\otimes
G_+\oplus E_-\otimes G_-,E_+\otimes G_-\oplus E_-\otimes G_+).
\label{29.4.2005.16}\end{multline}
Since the lifted bundles $E_\pm$ are trivial on each fibre, this has null
bundle precisely $E_+$ as a subbundle of $\CI\famtot;E_+\otimes H_+)$ and cokernel
bundle given by $E_-$ as a subbundle of $\CI(E_+\otimes H_-).$ Thus we can
interpret $A$ as mapping this null bundle into $E_-\otimes G_+$ and
so form the operator 
\begin{equation}
\begin{pmatrix}B&0\\A&B^*
\end{pmatrix} \in\ptP{\tfammap}{0,0}{\bfamtot;\sE\otimes\sG}.
\label{29.4.2005.17}\end{equation}
Directly from the definition, this is a Fredholm operator on $L^2.$
Since it has null space just the null space of $A$ as a subspace of the
null space of $B$ and cokernel the complement of the range of $A$ as a
subspace of the complement of the range of $B^*$ we conclude directly that 
\begin{equation}
\inda\begin{pmatrix}B&0\\A&B^*\end{pmatrix}=\inda(A).
\label{faficu.3}\end{equation}

Now we proceed to deform this operator as an elliptic family in
$\ptP{\tfammap}{0,0}{\bfamtot;\sE\otimes\sG}.$ First choose an
element $\tilde A\in\Ps0{\bfamtot;\sE\otimes G_+}$ with symbol
$\chi _2\sigma (A)\otimes\Id_{G_+}.$ As an element of
$\ptP{\tfammap}{0,0}{\bfamtot;\sE\otimes G_+}$ its symbol is the
lift of $\chi_2\sigma (A)$ and its base family is the lift of this symbol
to the front face, interpreted as a bundle map. Consider the 1-parameter family 
\begin{equation}
\begin{pmatrix}B&-s\tilde A^*\\ s\tilde A +(1-s)A&B^*
\end{pmatrix}\in \ptP{\tfammap}{0,0}{\bfamtot;\sE\otimes\sG}, s\in[0,1].
\label{faficu.4}\end{equation}
This is fully elliptic, \ie both symbol and base family are invertible for
all $s\in[0,1].$ For the symbol itself this is rather clear since it is the
family  
\begin{equation}
\begin{pmatrix}\sigma (b)&-s\chi _2\sigma(a)^*\\ s\chi_2\sigma(a)&\sigma(b)^*
\end{pmatrix}.
\label{faficu.146}\end{equation}
On the other hand the base family is 
\begin{equation}
\begin{pmatrix}B&-s\sigma (A)^*\\ s\sigma (A) +(1-s)\pi_-\sigma(A)\pi_+&B^*
\end{pmatrix}
\label{faficu.147}\end{equation}
as an operator on each fibre of $\tfammap,$ lifted to $S^*\famtot.$ Thus,
$\sigma (A)$ is simply being extended from the null space of $B$ to the
whole bundle and this family is invertible for all $s\in[0,1].$

Then we may choose $\tilde B\in\Ps0{\bfamtot;E_+\otimes\sG}$ with symbol
$\chi _1\sigma (B)\otimes\Id_{E_+}$ and similarly extend the
operator\eqref{faficu.4} at $s=1$ to a 1-parameter
elliptic family
\begin{equation}
\begin{pmatrix}(1-r)B+r\tilde B&-\tilde A^*\\ \tilde A &(1-r)B^*+r\tilde B^*
\end{pmatrix}\in\ptP{\tfammap}{0,0}{\bfamtot;\sE\otimes\sG}, r\in[0,1].
\label{faficu.5}\end{equation}
Again this is a fully elliptic family and at $r=1$ reduces to an element of
$\Ps0{\bfamtot;\sE\otimes\sG}$ which has index equal to that of $A$ and symbol
given by \eqref{29.4.2005.2}.
\end{proof}

\begin{remark}\label{faficu.204} In this construction it is already clear
that the K-class of the symbol of the resulting operator only depends on
the K-classes of the symbols of $A$ and $B$ and so this defines a map
\begin{equation}
M_b:\Kc(T^*(\famrel))\longrightarrow \Kc(T^*(\bfamrel)).
\label{faficu.210}\end{equation}

It should also be noted that if the operator (or family of operators) $A$
is actually invertible then the lifted family is invertible and this
invertibility can be maintained through the homotopy back to the standard
algebra. This is used below for full ellipticity in the cusp setting.

Despite this, in general a fully elliptic family of product-type cannot be
deformed through elliptics into the `standard' subspace since this
involves deforming the base family, through invertibles, to bundle maps and
there may be an obstruction to this in K-theory.
\end{remark}

Next we consider the corresponding construction in the cusp case; the main
obstacle to this extension is notational! Thus we have an iterated fibration. 
Here $\fammap:\famtot\longrightarrow \fambas$ is our usual fibration of compact
manifolds with fibre a compact manifold with boundary. The `top' fibration
$\tfammap:\bfamtot\longrightarrow \famtot$ is assumed to have compact fibre
a manifold without boundary (in the application to lifting it is a
sphere.) Again we are given an elliptic family $B\in\Psi^0(\tfamrel;\sG)$
which is surjective and has a trivial 1-dimensional null bundle where $\sG$
is some $\bbZ_2$-graded bundle over $\bfamtot.$

We refer to Appendix~\ref{PTFCO} for a discussion of product-type
cusp operators in this setting. These are quite analogous to the
product-type pseudodifferential operators discussed above. Such an operator
$P\in\ptcuP{\tfammap}{m,m'}{\bfamrel;\sF}$ for any $\bbZ_2$-graded
bundle $\sF$ over $\bfamtot$ acts on (weighted) smooth sections as in
\eqref{faficu.68} but there are now \emph{three} `symbol' maps. The
`commutative' symbol corresponds to the usual cusp symbol modified in a way
analogous to that of the symbol in the product-type pseudodifferential
calculus in \eqref{29.4.2005.6}
\begin{multline}
\ptcuP{\tfammap}{m-1,m'}{\bfamrel;\sF}\hookrightarrow
\ptcuP{\tfammap}{m,m'}{\bfamrel;\sF}\overset{\sigma _{m,m'}}\longrightarrow
 \ptcuPS{\tfammap}{m,m'}{\bfamrel;\sF},    \\
\ptcuPS{\tfammap}{m,m'}{\bfamrel;\sF}=
\CI([\cuS^*\bfamtot,(\tfammap)^*\cuS^*(\famrel)];
N^{-m}\otimes N_{\ff}^{-m'}\otimes\hom(\sF)).
\label{faficu.137}\end{multline}
Secondly is the indicial family, corresponding to the map \eqref{faficu.70}
for the usual cusp calculus, but now taking values in the suspended version
of the product-type calculus for the restriction of the fibration to the
preimage of the boundary of $\famtot$
\begin{equation}
\xymatrix@1{x\ptcuP{\tfammap}{m,m'}{\bfamrel;\sF}\ar[r]&
\ptcuP{\tfammap}{m,m'}{\bfamrel;\sF}\ar[r]^(.45){N}&
\sptP{\pa\tfammap}{m,m'}{\pa\bfamrel;\sF}}.
\label{faficu.138}\end{equation}
Finally there is an analogue of the non-commutative symbol in
\eqref{29.4.2005.8}, namely a short exact sequence 
\begin{equation}
\xymatrix@1{\ptcuP{\tfammap}{m,m'-1}{\bfamrel;\sF}\ar[r]&
\ptcuP{\tfammap}{m,m'}{\bfamrel;\sF}\ar[r]^(.37){L}&
\Psi^{m'}(\cuS^{*}_{\tfammap}M/\cuS^{*}M;\sF\otimes_{+}N^{-m}) }
\label{faficu.139}\end{equation}
taking values in the pseudodifferential operators on the fibres of
$\tfammap$ but lifted to the fibrewise cusp cosphere bundle $\cuS^*(\famrel).$
Whilst separately surjective these three maps are related and combine to
give a `full symbol algebra' consisting of the elements with the
compatibility conditions
\begin{multline}
\ptcujS{\tfammap}{m,m'}(\bfamrel)=\big\{(a,I,\beta)\in \\
 \ptcuPS{\tfammap}{m,m'}{\bfamrel;\sF}\times 
\sptP{\pa\tfammap}{m,m'}{\pa\bfamrel;\sF}\times
\Psi^{m'}(\cuS^{*}_{\tfammap}M/\cuS^{*}M;\sF\otimes_{+}N^{-m});\\
a\big|_{\pa}=\sigma(I),\ a\big|_{\ff}=\sigma (\beta),\ N(\beta)=L(I)
\big\}.
\label{faficu.140}\end{multline}
This space itself corresponds to the short exact sequence 
\begin{equation}
\xymatrix@1{x\ptcuP{\tfammap}{m-1,m'-1}{\bfamrel;\sF}\ar[r]&
\ptcuP{\tfammap}{m,m'}{\bfamrel;\sF}\ar[r]^{(\sigma,N,L)}&
\ptcujS{\tfammap}{m,m'}{(\bfamrel,\sF)}}
\label{faficu.142}\end{equation}
which captures compactness and Fredholm properties on the appropriate
Sobolev spaces. For us it suffices to note that 
\begin{multline}
A\in\ptcuP{\tfammap}{0,0}{\bfamrel;\sF}\Longrightarrow
A:L^2(\bfamrel;F_+)\longrightarrow L^2(\bfamrel;F_-)\text{ is bounded and}\\
A\text{ is compact}\Longleftrightarrow \sigma (A)=0,\ N(A)=0,\ L(A)=0\\
A\text{ is Fredholm}\Longleftrightarrow\ \exists\
(\sigma (A),\ N(A),\ L(A))^{-1}\in
\ptcujS{\tfammap}{-m,-m'}(\bfamrel,\sF^-).
\label{faficu.141}\end{multline}
In the latter case, we say $A$ is \emph{fully elliptic}.
As with the product-type pseudodifferential operators in the boundaryless
case, it is important that the three different spaces of pseudodifferential
operators, namely the ordinary pseudodifferential operators on the fibres
of the smaller fibration, the ordinary cusp pseudodifferential operators
on the fibres of the larger fibration and the cusp pseudodifferential
operators on the base fibration all lift into this larger space: 
\begin{equation}
\begin{gathered}
\Ps m{\tfamrel;\sF}\subset\ptcuP{\tfammap}{0,m}{\bfamrel;\sF}\\
\cuP m{\bfamrel;\sF}\subset\ptcuP{\tfammap}{m,m}{\bfamrel;\sF}\\
\cuP m{\famrel;\sE}\subset \ptcuP{\tfammap}{m,-\infty}{\bfamrel;\sF}
\end{gathered} 
\label{faficu.144}\end{equation}
where $\sE$ is a $\bbZ_2$ graded vector bundle embedded as a subbundle of
$\CI(\tfamrel;\sF)$ as a bundle over $\famtot.$

\begin{proposition}\label{faficu.143} For an iterated fibration with cusp
structure as in \eqref{faficu.135}, if an elliptic family $B\in\Ps
0{\tfamrel;\sF}$ is surjective and has trivial one-dimensional null bundle
$H$ then for any fully elliptic element $A\in\cuP 0{\famrel;\sE}$ the operator 
\begin{equation}
\begin{pmatrix}B&0\\A&B^*
\end{pmatrix} \in\ptcuP{\tfammap}{0,0}{\bfamrel;\sF\otimes\sE},
\label{faficu.145}\end{equation}
where the inclusions are through \eqref{faficu.144}, is a fully elliptic
element with the same index as $A$ in $\Kt(\fambas)$ and furthermore
this element is deformable, through fully elliptic elements, into
$\cuP0{\bfamrel;\sF\otimes\sE}.$ This results again in a map 
\begin{equation}
M_b:\cuK{\fammap}\longrightarrow \cuK{\fammap\circ\tfammap}
\label{faficu.211}\end{equation}
which only depends on $b=[\sigma(B)]\in\Kc(T^*(\tfamrel)).$
\end{proposition}

\begin{proof} As noted above, the first part of the result, that
\eqref{faficu.145} gives a fully elliptic element of the product-type
calculus on the total fibration, and that \eqref{faficu.3} again holds,
follows as in the proof of Proposition~\ref{29.4.2005.1} with only minor
notational changes.

For the second, deformation, part of the statement we need to proceed
with more care. First, the arguments of Proposition~\ref{29.4.2005.1} apply
unchanged in the sense that \eqref{faficu.4} and \eqref{faficu.5} together
give a symbolically elliptic deformation, indeed the computation of the
symbol (although now for a family in the base) is the same as there, as is
the computation of the base family -- the latter is still given by
\eqref{faficu.147} provided the symbol of $A$ is interpreted as the cusp
symbol. Thus we have constructed a family $P_t,$ in the product-type cusp
calculus, where we can relabel and change parameterization to make the
family smooth in $t\in[0,1].$ This family is elliptic for all $t,$ has
invertible base family for all $t,$ is Fredholm at $t=0$ and is in the
ordinary cusp calculus at $t=1.$ So, consider the indicial family of $P_t.$
This is a suspended family of product-type pseudodifferential operators
(acting on the fibres of a fibration) which is fully elliptic, \ie is
elliptic and has invertible base family, and which is invertible at $t=0.$
Since it is also invertible for large values of the suspending variable, it
follow by standard arguments that it can be perturbed by a family of smoothing
operators (see Remark~\ref{faficu.204}) which vanishes near infinity in the
suspension variable and vanishes near $t=0$ to be invertible for all values
of $t.$ Modifying the family in this way results in a deformation as desired.

That this construction is symbolic as far as $B$ is concerned and `fully
symbolic' as far as $A$ is concerned, in the sense of \eqref{faficu.211},
follows readily from the construction. Namely it certainly behaves
well under stabilization and homotopy, with the homotopy parameter simply
being added to the base variables. So it remains to show that it is stable
under different regularizations of the null bundle of $B;$ in fact any two
such stabilizations are homotopic.
\end{proof}

Finally we comment on the generalization of this result to the fibred cusp
case although this is not used below. Thus we again consider an iterated
fibration \eqref{faficu.135} where now the second fibration has a fibred
cusp structure, $\boumap:\pa\famtot\longrightarrow \boubas.$ There are two
extreme possibilities for a fibred cusp structure for the overall
fibration. First, it is always possible to `add' the fibres of
$\pa\bfamtot\longrightarrow\pa\famtot$ to the fibres of the boundary, \ie
to take the boundary fibration for the overall fibration to be
\begin{equation}
\bboumap=\boumap\circ\tfammap.
\label{faficu.207}\end{equation}

\begin{proposition}\label{faficu.209} For an iterated fibration as in
\eqref{faficu.135} where the second fibration has a fibred cusp structure
and the top fibration has the fibred cusp structure
\eqref{faficu.207}, if an elliptic family $B\in\Ps 0{\tfamrel;\sF}$ is
surjective and has trivial one-dimensional null bundle $H$ then for any
fully elliptic element $A\in\fcP{\boumap}0{\famrel;\sE}$ the operator
\begin{equation}
\begin{pmatrix}B&0\\A&B^*
\end{pmatrix} \in\ptfcP{\tfammap}{\bboumap}{0,0}{\bfamrel;\sF\otimes\sE},
\label{faficu.145a}\end{equation}
where the inclusions are through \eqref{faficu.144}, is a fully elliptic
element with the same index as $A$ in $\Kt(\fambas)$ and furthermore
this element is deformable, through fully elliptic elements, into
$A_B\in\fcP{\bboumap}0{\bfamrel;\sF\otimes\sE}.$ This results again in a map 
\begin{equation}
M_b:\fcK{\boumap}{(\fammap)}\longrightarrow
\fcK{\boumap\circ\tfammap}{(\fammap\circ\tfammap)}
\label{faficu.211a}\end{equation}
which only depends on $b=[\sigma(B)]\in\Kc(T^*(\tfamrel)).$
\end{proposition}

\begin{proof} The arguments in the proof of Proposition~\ref{faficu.143} go
through essentially without change, the only difference being the
appearance of more suspension parameters in the normal operator.
\end{proof}

There is a second extreme possibility, other than \eqref{faficu.207}, for
the fibred cusp structure on the top fibration in \eqref{faficu.135},
corresponding to the fibres for $\bboumap:\bfamtot\longrightarrow \bboubas$
being of the same dimension as for $\boumap$ in which case all the `new'
boundary variables are in the base of the fibration. Then there 
is a commutative diagram of fibrations
\begin{equation}
\xymatrix{\sfamfib\ar@{-}[r]&\pa\bfamtot\ar[d]^{\tfammap}
\ar[r]^{\bboumap}&\bboubas\ar[d]&
\ar@{-}[l]\sfamfib\\
&\pa\famtot\ar[r]^{\boumap}&\boubas}
\label{faficu.208}\end{equation}

In this case there is a difficulty in extending the construction above, in
that it is not quite clear what the family $B$ should be. Approached
directly $B$ would need to be a family which is adiabatic in the boundary
variable. Since the index theory for such families has not been properly
developed we have chosen to proceed more indirectly by using
Proposition~\ref{faficu.215}. However, it is possible instead to use
Proposition~\ref{faficu.209} above and then make an adiabatic limit back to
the case where the extra variables are in the base of the fibration. Since
this involves a rather delicate homotopy, and is not used below, the
details are omitted.

\section{Lifting and excision}\label{Lifting}

Next we consider the lifting construction for cusp K-theory with respect
to an extension of a given fibration, see Definition~\ref{faficu.149}.

\begin{proposition}\label{faficu.150} If $\efammap:\efamtot\longrightarrow
  \fambas$ is an extension of a fibration $\fammap:\famtot\longrightarrow
  \fambas$  with cusp structure in the sense of Definition~\ref{faficu.149}
  then there is a well-defined lifting homomorphism
\begin{equation}
(\efammap/\fammap)^!:\cuK{\fammap}\longrightarrow \cuK{\efammap}
\label{faficu.151}\end{equation}
which induces a commutative diagram \eqref{faficu.35} where the map on
 the right is an absolute form of the lifting map of Atiyah and Singer.  When
 $\efamtot=\bbB^{2N+1}\times \fambas$, 
  this reduces to pull-back under Poincar\'e duality 
\begin{equation}
\xymatrix{\cuK{\efammap}\ar[r]^(0.35){\quan}&\Kf{\boubas}
{\CO_{\cusp}(\efamtot),\CO(\fambas)}\\
\cuK{\fammap}\ar[u]^{(\efammap/\fammap)^!}\ar[r]^(0.35){\quan}&
\Kf{\boubas}{\CO_{\cusp}(\famtot),\CO(\fambas)}\ar[u]_{R^*}}
\label{faficu.206}\end{equation}
where $R:\CO_{\cusp}(\efamtot)\to \CO_{\cusp}(\famtot)$ is the restriction
map.
\end{proposition}

\begin{proof} The normal bundle to $\famtot$ in $\efamtot$ is a real vector
  bundle which models the extended fibration near $\famtot.$ Thus, if we
  take the one-point compactification of the fibres we obtain a bundle of
  spheres $\pi:S_V\longrightarrow \famtot$ which gives an extension of the
  original fibration but is also an iterated fibration in the sense of
  Proposition~\ref{faficu.143}. Moreover, following the approach of Atiyah
  and Singer, there is a `Bott element', which we can realize as a family 
\begin{equation}
B\in\Ps0{S_V/\famtot;\bbL}
\label{faficu.152}\end{equation}
for a $\bbZ_2$-graded bundle $\bbL$ over $S_V$ with $L_+$ and $L_-$
  identified near the section at infinity and such that $B-\Id$ has kernel
  with support disjoint from the section at infinity (in either
  factor). Now, we can apply Proposition~\ref{faficu.143} to `lift' each
  fully elliptic element of $\cuP0{\famrel;\sE}$ to a fully elliptic element of
  $\cuP0{S_V/\fambas;\sF}.$ A brief review of the construction shows that
  the condition that all operators be equal to the identity near the
  section at infinity can be maintained throughout. In view of this
  property of the kernel, including an identification of $F_+$ and $F_-$
  near the section at infinity, these operators may be trivially extended
  into $\cuP0{\efamrel;\sF}$ to be fully elliptic. The nature of this
  construction shows that this does indeed construct a map
  \eqref{faficu.151} and that it leads to a commutative diagram
  \eqref{faficu.35}.  When $\efamtot= \bbB^{2N+1}\times \fambas$, the
  commutativity of diagram~\eqref{faficu.206} follows from proposition~\ref{faficu.20}
  and the fact that $\efammap= \fammap\circ R$.
\end{proof}

Note that the Poincar\'e duality isomorphism and \eqref{faficu.206} show
that the lifting map is independent of choices and is well-behaved under
composition. This is also relatively easy to see directly.

\section{Product with a ball and a sphere}\label{Products}

We have shown in Proposition~\ref{faficu.114} that a fibred cusp structure
over $\fambas$ can always be `trivialized' by embedding it in one of the
product cases with both fibrations being projections, in which case the
diagram \eqref{faficu.2} becomes
\begin{equation}
\xymatrix{&&\bbB^{M+1}\times\bbS^{N}\ar@{-}[r]&
\fambas\times\bbB^{M+1}\times\bbS^{N}\ar[d]^{\fammap=\pi}\\
\bbS^N\ar@{-}[r]&\bbS^M\times\bbS^N\ar@{^(->}[ur]^{\pa}
\ar[d]^{\boufibmap=\pi'}
\ar@{-}[r]&\fambas\times\bbS^M\times\bbS^N
\ar[r]^{\pa\fammap=\pi}\ar[d]^{\boumap=\pi'}
\ar@{^(->}[ur]^{\pa}&\fambas\\
&\bbS^M\ar@{-}[r]&\ar[ur]\fambas\times\bbS^M}
\label{faficu.98}\end{equation}
with $\pi$ being projection onto the leftmost factor and $\pi'$ off the
rightmost.

\begin{proposition}\label{faficu.20} If $M>0$ and $N\ge 0$ in the product spaces in
\eqref{faficu.98}, then
\begin{equation}
   \fcK{\pi'}(\pi)\cong
\begin{cases}
\Kt(\fambas)& M+N\text{ even}\\
\Kt(\fambas)\oplus\Kt(\fambas)&M+N\text{ odd}
\end{cases}.
\label{faficu.99}\end{equation}
In the cusp case with $\pi: \fambas\times \bbB^{M+1}\to B$, we have 
\[
\cuK{\pi}\cong
\begin{cases}
\Kt(\fambas)& M\text{ even}\\
\Kt(\fambas)\oplus\Kt(\fambas)&M\text{ odd}
\end{cases}.
\]
\end{proposition}

\begin{proof} From the Poincar\'e duality isomorphism,
Theorem~\ref{faficu.88}, we can use KK-theory to perform the
computation.  Now, for any product
$\famtot=\fambas\times\famfib$ in which the overall fibration is the projection
onto the left factor and the boundary fibration is the product with some
fibration of the boundary of $\boufibmap:\pa\famfib\longrightarrow\boubas,$ 
\begin{equation}
\CO_{\Id\times\boufibmap}(\fambas\times\famfib)=
\CO(\fambas)\hat\otimes\CO_{\boufibmap}(\famfib)
\label{faficu.100}\end{equation}
is the completed tensor product. From the general properties of KK-theory
it follows that 
\begin{equation}
\begin{aligned}
\KK_{\fambas}(\CO_{\Id\times\boufibmap}(\fambas\times\famfib),\CO(\fambas))
&=\KK_{\fambas}(\CO(\fambas)\hat\otimes\CO_{\boufibmap}(\famfib) ,\CO(\fambas))
\\
&\cong
\KK(\CO_{\boufibmap}(\famfib) ,\CO(\fambas)) \\
&\cong \Kt(\fambas)\otimes\KK(\CO_{\boufibmap}(\famfib),\bbC),
\end{aligned}
\label{faficu.101}\end{equation}
where in the last step we used the K\"unneth formula for KK-theory
(see for instance Theorem 23.1.3 in \cite{Blackadar1}) and we used the
fact (proved below) that $\KK(\CO_{\boufibmap}(\famfib),\bbC)$ is a free 
$\bbZ$-module.
Thus it suffices to consider the case in which the base in
\eqref{faficu.98} is reduced to a point and to show then that
$\fcK{\pi'}(\pi)$ reduces to one or two copies of $\bbZ$ according
to parity. In fact the general argument, with the base factor retained, is
not much more complicated.

Thus we need only to consider
\begin{equation}
\xymatrix{\bbS^N\ar@{-}[r]&\bbS^M\times\bbS^N\ar@{^(->}[r]^{\pa}
\ar[d]^{\boufibmap=\pi'}&\bbB^{M+1}\times\bbS^{N}\\
&\bbS^M}
\label{faficu.102}\end{equation}
and compute the fibred-cusp K-theory in this case.

First consider the cusp case, where there is no factor of $\bbS^N.$ The
same argument applies to reduce the computation to the case of
$\bbB^{M+1}.$ So consider the (split) short exact sequence \eqref{faficu.23},
proved in Section~\ref{aco.0}. Since $K(\{\text{pt}\})=\bbZ,$ essentially
by definition, it suffices to check that 
\begin{equation}
\Kc(\bbB^{M+1}\times\bbR^{M+1})=
\begin{cases}\{0\}&M\text{ even}\\
\bbZ&M\text{ odd.}
\end{cases}
\label{faficu.104}\end{equation}
Here, recall that the ball is taken to be closed, so the K-theory is
absolute in that factor and compactly supported in the Euclidean
factor. Since the ball is contractible, the only source for K-classes is 
\begin{equation*}
\Kc(\bbR^{M+1})\equiv\Kc(\bbB^{M+1}\times\bbR^{M+1})
\label{faficu.105}\end{equation*}
so \eqref{faficu.104} follows.

In the fibred cusp case, set
$Z_{M,N}=\bbB^{M+1}\times\bbS^N.$   The obstruction to perturbing an elliptic
operator so that it becomes fully elliptic lies in
\begin{equation*}
\Kt(T^{*}\bbS^{M}\times\bbR)\cong\Kt(\bbS^{M}\times \bbR^{M+1}).
\label{faficu.106}\end{equation*}
Consider the 6-term exact sequence coming from the inclusion
\begin{equation*}
    r_{\pa \bbB^{M+1}}: T^{*}\bbS^{M}\times\bbR\cong
                          T^{*}(\bbB^{M+1})\big|_{\pa \bbB^{M+1}}
\hookrightarrow  T^{*}(\bbB^{M+1}),
\end{equation*}
namely
\begin{equation}
\xymatrix{
\Kcn{0}(\bbR^{2M+2})\ar[r]&
\Kcn{0}(T^{*}(\bbB^{M+1}))
\ar[r]^{r_{\pa \bbB^{M+1}}^{*}}& \Kcn{0}(T^{*}(\bbS^{M})\times\bbR ) \ar[d]\\
\Kcn{1}(T^{*}(\bbS^{M})\times\bbR )\ar[u]&\ar[l]_{r_{\pa \bbB^{M+1}}^{*}}
 \Kcn{1}(T^{*}(\bbB^{M+1}))&
\ar[l] \Kcn{1}(\bbR^{2M+2}) }  
\label{faficu.108}\end{equation}
using the identification
\begin{equation}
  \Kcn{k}(T^{*}(\bbB^{M+1}), T^{*}(\bbS^{M})\times\bbR )\cong 
   \Kcn{k}(\bbR^{2M+2})\cong
\begin{cases}
     \bbZ & k=0,\\
    \{0\} & k=1.
\end{cases} 
\label{iden}\end{equation}
From \eqref{faficu.104} and \eqref{iden} and using the fact that 
\begin{equation}
  \Kcn{k}(\bbR^{2M+2})\longrightarrow
\Kcn{k}(T^{*}(\bbB^{M+1}) )  
\label{iden.1}\end{equation}
maps to zero, we get that
\begin{equation}
    \K(T^{*}\bbS^{M})\cong
\begin{cases}
  \bbZ\oplus\bbZ & M\text{ even}\\
    \bbZ & M\text{ odd}
\end{cases},\quad 
\Kt^{1}_{c}(T^{*}\bbS^{M})\cong
\begin{cases}
  \{0\} & M\text{ even} \\
    \bbZ  & M\;\text{ odd.}
\end{cases}
\label{efe.7}\end{equation}

Thus, in particular, when $M$ is even, there is no obstruction to the
existence of a smoothing perturbation that makes an elliptic fibred cusp
operator fully elliptic. So consider the image of
\begin{equation*}
   I_{1}: \Kcn{1}(T^{*}(\bbB^{M+1}\times\bbS^{N}))\longrightarrow \Kcn{0}(T^{*}(\bbS^{M})).
\end{equation*}
Since $\pi': \bbS^{M}\times\bbS^{N}\longrightarrow \bbS^{M}$ extends to
\begin{equation*}
P': \bbB^{M+1}\times\bbS^{N}\longrightarrow \bbB^{M+1}
\end{equation*}
by projecting on the right factor
\begin{equation*}
                I_{1}=\ind_{\pi'}r^{*}_{\pa(\bbB^{M+1}\times\bbS^{N})}=
              r^{*}_{\pa(\bbB^{M+1})} \ind_{P'}.
\end{equation*}
But $\ind_{P'}$ is clearly surjective, so the image of $I_{1}$ is the same
as the image of $r^{*}_{\pa(\bbB^{M+1})}.$  Thus, we conclude from 
the 6-term exact sequence \eqref{faficu.108} that
\begin{equation*}
       I_{1}({\Kcn{1}(T^{*}(\bbB^{M+1}\times\bbS^{N}))})\cong 
         \begin{cases}
  \bbZ & M\text{ even} \\
    \{0\}  & M\;\text{ odd.}
\end{cases}
\end{equation*}
Thus, when $M$ is even, there is a short exact sequence
\begin{equation*}
\xymatrix@1{ \bbZ \ar[r]&
\fcK{\pi'}(\bbB^{M+1}\times \bbS^{N})\ar[r]&
\Kc(T^{*}(\bbB^{M+1}\times\bbS^{N})).}
\end{equation*}
Since
\begin{equation*}
   \K(T^{*}(\bbB^{M+1}\times\bbS^{N})) \cong \begin{cases}
  \{0\} & N\text{ even}\\
    \bbZ & N\text{ odd}
\end{cases}
\label{faficu.112}\end{equation*}
is a free $\bbZ$-module, this sequence splits and hence
\begin{equation*}
  \fcK{\pi'}(\bbB^{M+1}\times \bbS^{N}) 
 \cong \bbZ \oplus \Kt(T^{*}(\bbB^{M+1}\times\bbS^{N}))\cong   
\begin{cases}
  \bbZ & N\text{ even}\\
    \bbZ\oplus\bbZ  & N\text{ odd.}
\end{cases}
\label{faficu.113}\end{equation*}
When $M$ is odd, we have instead the short exact sequence 
\begin{equation*}
\K(T^{*}\bbS^{M})\longrightarrow \fcK{\pi'}(\bbB^{M+1}\times\bbS^{N})
\longrightarrow\ker(I_{0})
\label{faficu.109}\end{equation*}
which splits by sending an element of $\ker(I_{0})$ to a full symbol with
null index. Thus we conclude that
\begin{equation*}
\fcK{\pi'}(\bbB^{M+1}\times \bbS^{N})\cong \ker(I_{0})
\oplus\bbZ \cong   \K(T^{*}(\bbB^{M+1}\times \bbS^{N}))     
\label{faficu.110}\end{equation*}
since the space of obstruction is $\Kco(T^{*}\bbS^{M})\cong \bbZ$ and 
$I_{0}$ is surjective in this case, as one can see from the 6-term 
exact sequence \eqref{faficu.108}. Using again \eqref{efe.7}
\begin{equation}
\begin{aligned}
  \fcK{\pi'}(\bbB^{M+1}\times \bbS^{N}) &\cong
 \Kt_{c}(T^{*}(\bbB^{M+1}\times\bbS^{N}))\\
  &\cong \Kt_{c}( T^{*}\bbS^{N}\times\bbR^{M+1}) \\
 &\cong \Kt_{c}(T^{*}\bbS^{N})\cong   
\begin{cases}
  \bbZ\oplus\bbZ & N\text{ even}, \\
    \bbZ & N\text{ odd.}
\end{cases}
\end{aligned}
\label{efe.10}\end{equation}
\end{proof}

\section{The topological index}\label{ti.0}

Consider a fibration as in \eqref{faficu.2}. Let  
\begin{equation*}
    i:\famtot \hookrightarrow \fambas\times \bbB^{p+1},
\end{equation*}
with $p$ even, be an embedding as in Proposition~\ref{faficu.114}, where the
fibration structure on $\fambas\times\bbB^{p+1}$ is given by the projection
on the right factor
\begin{equation*}
\pi:\fambas\times \bbB^{p+1}\longrightarrow \fambas   .
\end{equation*}
Proposition~\ref{faficu.150} gives a well-defined lifting
homomorphism
\begin{equation}
    (\pi/\fammap)^!:\cuK{\fammap}\longrightarrow \cuK{\pi}.\label{ti.1}
\end{equation}
By the proof of Proposition~\ref{faficu.20}, 
\begin{equation*}
     \cuK{\pi}\overset{\quan}{\longrightarrow} \KK_{\fambas}(\CO_{\pa\pi}(
 \fambas\times \bbB^{p+1}), \CO(\fambas))\overset{\ind}{\longrightarrow}
   \Kt(B)
\end{equation*}
is an isomorphism.
\begin{definition}
the topological index map $\indt: \fcK{\boumap}\longrightarrow\Kt(B)$ 
is defined by
\begin{equation*}
\indt=\ind \circ \quan\circ  (\pi/\fammap)^! \circ q_{\ad}\; .
\end{equation*}
The fact that it does not depend on the choice of the embedding follows
from the stability of the lifting map under repeated embedding, but in any
case will follow from Theorem~\ref{faficu.18}!
\label{ti.2}\end{definition}

\begin{proof}[Proof of Theorem~\ref{faficu.18}]
The theorem follows from the commutativity of the diagram \eqref{faficu.35}
stated in Proposition~\ref{faficu.150}, and the commutativity of diagram
\eqref{faficu.220} stated in Proposition~\ref{faficu.215}.
\end{proof}

\section{Families of Atiyah-Patodi-Singer type}\label{APS}

The Atiyah-Patodi-Singer index theorem of \cite{APS1} was originally proved
with the idea of obtaining a generalization of Hirzebruch's signature
theorem for the case of manifolds with boundary. Since then, variants of
the proof of the theorem were obtained, for instance in \cite{Cheeger1} and
\cite{MelroseAPS}. An extension to the family case was discussed in
\cite{Bismut-Cheeger1}, \cite{Bismut-Cheeger2} and in \cite{MR99a:58144},
\cite{MR99a:58145}.  In all these proofs, including the original one, heat
kernel techniques for Dirac operators play an important r\^{o}le. As
opposed to the Atiyah-Singer index theorem, the Atiyah-Patodi-Singer index
theorem was originally restricted to Dirac operators, but using trace
functional techniques, a pseudodifferential generalization was obtained in
\cite{Melrose-Nistor} in the setting of cusp operators. Such a
generalization was also discussed by Piazza in \cite{Piazza} for
b-pseudodifferential operators.

In \cite{Dai-Zhang1}, Dai and Zhang provided an interesting proof of the
Atiyah-Patodi-Singer index theorem by embedding the manifold with boundary
into a large ball.  Relating the Dirac operator of interest with one
defined on the large ball via a careful analysis, they were able to take
advantage of the simple topology of the ball to get the
Atiyah-Patodi-Singer index theorem for the original operator. This is
certainly closely related to our constructions above, however, the methods
of \cite{Dai-Zhang1} are essentially analytical and no K-theory is
involved.

We proceed to briefly recall the setting of Atiyah-Patodi-Singer boundary
problem and its reformulation in terms of cusp pseudodifferential
operators, leading to a K-theory index. Let $Z$ be an even dimensional 
Riemannian manifold with nonempty boundary $\pa  Z=X$ with a Riemannian
metric which is of product type near the boundary, so there is a
neighborhood $X\times [0,1)\subset Z$ of the boundary in which the metric 
takes the form
\begin{equation}
g= du^{2} + h_{X}
\label{dirac.1}\end{equation}
where $u\in [0,1)$ is the coordinate normal to the boundary and $h_{X}$ is the 
pull-back of a metric on $X$ via the projection
$X\times[0,1)\longrightarrow X.$ Let $\sE$ be a Hermitian vector bundle
over $Z$ with a Clifford module structure for the metric \eqref{dirac.1}
and with a unitary Clifford connection which is constant in the normal
direction near the boundary under a product trivialization. This defines a
generalized Dirac operator
\begin{equation*}
\eth:\mathcal{C}^{\infty}(Z,\sE)\longrightarrow \mathcal{C}^{\infty}(Z,\sE).
\end{equation*} 
In the neighborhood
$X\times [0,1)\subset Z$ of the boundary described above, it takes the form
\begin{equation}
\eth^{+} =\gamma \left(\frac{\pa}{\pa u} + A \right)
\label{dirac.2}\end{equation}
where $\gamma =cl(u):\sE \big|_{X}\longrightarrow\sE \big|_{X}$ is given
by Clifford multiplication by the normal differential and 
$A: \mathcal{C}^{\infty}(X,\left. E_{+} \right|_{X})\longrightarrow 
\mathcal{C}^{\infty}(X,\left. E_{+}  
\right|_{X})$ is a Dirac operator on $X$ such that
\begin{equation}
  \gamma^{2}=-\Id,\ \gamma^{*}=-\gamma,\ A\gamma=-\gamma A,\ A^{*}=A .
\label{dirac.3}\end{equation}   
Consider the spectral boundary 
condition
\begin{equation}
\varphi \in \mathcal{C}^{\infty}(Z,\sE),\ P(\varphi\big|_{X})=0,
\label{dirac.4}\end{equation}
where $P$ is the projection onto the nonnegative spectrum of $A.$ Then
\begin{equation}
       \eth^{+}: W^{1}_{P}\longrightarrow \Ld(Z;E_{-})
\label{APS.1}\end{equation}
is a Fredholm operator,
where 
\begin{equation*}
    W^{1}_{P}=\left\{f\in \operatorname{H}^{1}(X;E_{+});P(f\big|_{X})=0\right\}
\end{equation*}
is a subspace of $\operatorname{H}^{1}(X;E_{+}),$ the Sobolev space
of order 1.

Atiyah, Patodi and Singer show in \cite{APS1} that the index of $\eth^{+}$
is given by
\begin{equation*}
\ind(\eth^{+})=\int_{Z}\Ahat(Z)\Ch'(\sE) - \frac{h+\eta}{2}
\end{equation*}
where $\Ch'(\sE)$ is the twisting Chern of $\sE,$ $\Ahat$ is the
$\Ahat$-genus, $h=\dim \ker A$ and $\eta$ is the eta invariant of $A.$

As discussed in \cite{APS1}, one can alternatively describe the index problem
by adding a cylindrical end to the manifold with boundary $Z.$  More precisely,
$Z$ may be enlarged by attaching the half-cylinder $(-\infty,0)\times X$ to
the boundary $X$ of $Z.$ Call the resulting manifold $Z'.$ The metric,
being a product near the boundary, can be naturally extended to this
half-cylinder, which makes the resulting manifold a complete Riemannian
manifold. The bundle, Clifford structure, connection and hence the Dirac
operator also have natural translation-invariant extensions to $Z'$
using the product structure near the boundary. On $Z',$ it is possible to
think of $\eth$ as a cusp operator by compactifying to a compact manifold
with boundary (diffeomorphic to the original $Z)$ by replacing $u$ by the
variable
\begin{equation}
  x=-\frac{1}{u} \in[0,1)
\label{APS.3}\end{equation}
for $u\in (-\infty,-1).$ The extension down to $x=0,$ gives the manifold
with boundary $\overline{Z'},$ and $x$ is a boundary defining function for
$\pa \overline{Z'}\cong X$ which defines a cusp structure.  Let us denote
by $\dc$ the natural extension of $\eth$ to $\overline{Z'}.$  Near the
boundary of $\overline{Z'},$ $\dc$ takes the form
\begin{equation}
\dc=\gamma\left( x^{2}\frac{\pa}{\pa x}+A\right) 
\label{APS.4}\end{equation}
and so is clearly an elliptic cusp differential operator. 

\begin{lemma}
If $A$ is invertible, then
\begin{equation*}
         \dc^{+}: \operatorname{H}^{1}(\overline{Z'};E_{+})
\longrightarrow \Ld(\overline{Z'};E_{-})
\end{equation*}
is Fredholm and has the same index has the operator \eqref{APS.1}.
\label{APS.5}\end{lemma}

\begin{proof} Recall (\cite{Mazzeo-Melrose4}) that a cusp operator is
Fredholm if and only if it is elliptic and its indicial family is invertible.
The indicial family of $\dc^{+}$ is given by
\begin{equation*}
 e^{-i\frac{\tau}{x}}\dc^{+}e^{i\frac{\tau}{x}}
\big|_{x=0} = \gamma( A-i\tau),\ \tau\in\bbR.
\end{equation*}
Since $A$ is self-adjoint, this is invertible for all $\tau$ if and only if
$A$ is invertible. Thus, $\dc$ is Fredholm if and only if $A$ is
invertible. That $\dc^{+}$ has the same index as \eqref{APS.1} then follows
from Proposition 3.11 in \cite{APS1} and Proposition 9 in \cite{Mazzeo-Melrose4}.
\end{proof}

When $A$ is not invertible, it is still possible to relate \eqref{APS.1} with
a Fredholm cusp operator.  In fact, since this is basically the problem we 
encounter when we consider the family version, let us immediately generalize to
this context.  For the family case consider a fibration, \eqref{faficu.6},
of a manifold with boundary and we assume now that $\eth$ is a family of
Dirac operators parameterized by $B$ which as before are of product near
the boundary of $\pa M.$ Thus, there is an associated family $A$ of
self-adjoint Dirac operators on the boundary. The main difficulty in the
family case is that interpreted directly, the spectral boundary condition
need not be smooth. In \cite{MR99a:58144}, this difficulty was
overcome by introducing the notion of a spectral section.
\begin{definition}
A spectral section for a family of elliptic self-adjoint operators 
\begin{equation*}  A \in  
\operatorname{Diff}^{1}(\pa \famrel;\left. \sE\right|_{\pa \famtot})
\end{equation*} 
is a family of 
self-adjoint projections 
$P\in \Psi^{0}(\pa \famrel; \left.\sE\right|_{\pa \famtot})$ such that for
some smooth function $R:B\longrightarrow [0,\infty)$ (depending on $P$) and every $b\in B,$
\begin{equation*}
               A_{b}f=\lambda f \Longrightarrow
  \begin{cases}
P_{b}f=f & \Mif\lambda > R(b), \\
P_{b}f=0  & \Mif\lambda < -R(b).
\end{cases}
\end{equation*}
\label{APS.6}\end{definition}
Such a spectral section always exists for the boundary family and any such
choice gives a smooth family of boundary problems
problem 
\begin{equation}
\eth^{+}: W^{1}_{P}\longrightarrow \Ld(\famrel;E_{-})
\label{APS.7}\end{equation}
where
\begin{equation*}
    W^{1}_{P}=\left\{ f\in
    \operatorname{H}^{1}(\famrel;E_{+});P(f\big|_{\pa M})=0\right\}.
\end{equation*}
The family $\eth^{+}$ in \eqref{APS.7} is Fredholm so has a well defined
families index. As before, one can attach a cylindrical end and
get a new fibration $\phi:\famtot'\longrightarrow\fambas$ where the family
of operators $\eth$ naturally extends.  By making the  change of variable
$x=-\frac{1}{u},$ one get the a family of cusp operators 
$\dc$ which takes the form \eqref{APS.4} near the boundary.

\begin{lemma}
There exists $Q\in \cuP{-\infty}{\famrel;E_{+},E_{-}}$ such that 
$\dc^{+} + Q$ is a fully elliptic (hence Fredholm) family with the same family 
index as \eqref{APS.7}.
\label{APS.8}\end{lemma}

\begin{proof}
This is carried out in section 8 of \cite{MR99a:58144} in the context of
b-pseudodifferential operators instead of cusp operators. Since the
relationship corresponds to the introduction of the transcendental variable
$1/u$ (for cusp) instead of $e^u$ (for b-) one can check that the argument
continues to hold for cusp operators with only minor modifications. 
\end{proof}

\begin{proof} [Proof of Theorem~\ref{faficu.153}]
Let $(\eth,P)$ be as in the proposition.  Define
$[(\eth,P)]\in\cuK{\fammap}$ to be the K-class associated to the operator
$\dc^{+}+Q$ of Lemma~\ref{APS.8}. Then by the lemma 
\begin{equation*}
             \ind(\eth,P)=\inda([(\eth,P)])=\indt([(\eth,P)]).
\end{equation*}
That $[(\eth,P)]$ is canonically defined, that this, does not depend on the
choice of $Q$ in Lemma~\ref{APS.8}, is a consequence of the 
relative families index theorem of \cite{fipomb}.  Thus, given the Poincar\'{e}
duality of Theorem~\ref{faficu.88} and the commutativity of diagram
\eqref{faficu.30}, Theorem~\ref{faficu.153} follows.
\end{proof}

\appendix

\section{Product-type pseudodifferential operators}\label{Appendix-PT}

There are many variants of pseudodifferential operators `of product type'
to be found in the literature, see Shubin \cite{shubin3}, Rodino
\cite{MR0397800}, Melrose-Uhlmann \cite{MR81d:58052}, H\"ormander
\cite{hormander3}. Here we describe, succinctly, a particularly natural
algebra of such operators associated to a fibration (and for families to an
iterated fibration). First consider the local setting.

On $\bbR^{d_1}_y\times\bbR^{d_2}_z$ we consider operators corresponding to
this product thought of as a fibration over the first factor. The class of
symbols admitted is determined by the \ci\ structure on the manifold with
corners 
\begin{equation}
X_{d_1,d_2}=\bbR^{d_1}\times\bbR^{d_2}\times
[\com{\bbR^{d_1+d_2}};\pa\com{\bbR^{d_1}\times\{0\}}].
\label{16.6.2005.1}\end{equation}
That is, take the radial compactification of $\bbR^{d_1+d_2}$ and blow up
the boundary (at infinity) of the radial compactification of the subspace
$\bbR^{d_1}\times\{0\}.$ Let $\rho,$ $\rho _{\ff}\in\CI(X_{d_1.d_2})$ be
defining functions for the two boundary hypersurfaces, the first being the
`old hypersurface' at infinity and the second that produced by the blow up.

Now, if $a\in\rho ^{-m}\rho _{\ff}^{-m'}\CIc(X_{d_1,d_2})$ then it satisfies
the estimates 
\begin{equation}
|\pa_y^\alpha \pa_z^\beta \pa_{\eta}^\gamma \pa_{\zeta }^\delta a|
\le C_{\alpha ,\beta ,\gamma ,\delta }
(1+|\zeta |)^{m-m'-|\delta|}(1+|\eta |+|\zeta|)^{m'-|\gamma|},
\label{16.6.2005.2}\end{equation}
as follows by noting that one can take
$\rho=(1+|\zeta|^2)^{-\frac12}$ and $\rho _{\ff}=(1+|\zeta
|^2)^{\frac12}(1+|\eta |^2+|\zeta |^2)^{-\frac12}.$ This gives the overall
weight in \eqref{16.6.2005.2} with no differentiation. The vector fields 
\begin{equation*}
\pa_{y_j},\ \pa_{z_k},\ \eta _i\pa_{\eta _j},\ \zeta _l\pa_{\zeta _k},\
\zeta _k\pa_{\eta _j} 
\label{16.6.2005.3}\end{equation*}
all lift to be smooth on $\bbR^{d_1+d_2}\times\com{\bbR^{d_1+d_2}}$ and
tangent to the boundary and within the boundary to the submanifold blown up
in \eqref{16.6.2005.1} so
\begin{equation}
\pa_y^\alpha \pa_z^\beta \pa_{\eta}^\gamma \pa_{\zeta }^\delta a \in
\rho ^{-m+\gamma+\delta}\rho _{\ff}^{-m'+\gamma}\CIc(X_{d_1,d_2});
\label{20.7.2005.1}\end{equation}
this leads directly to the estimates \eqref{16.6.2005.2}. Consider the
kernels on $\bbR^{2d_1+2d_2}$ defined by Weyl quantization of symbols 
\begin{multline}
A(y,z,y',z')=(2\pi)^{-d_1-d_2}\int e^{i(y-y')\cdot\eta+(z-z')\cdot\zeta}
a(\frac{y+y'}2,\frac{z+z'}2,\eta ,\zeta )d\eta d\zeta \Mwhere  \\
a=a_1+a_2+a_3,\ a_1\in \rho ^{-m}\rho _{\ff}^{-m'}\CIc(X_{d_1,d_2}),\\
a_2\in\rho ^{\infty}\cS(\bbR^{2d_2};
\rho_\eta^{-m'}\CIc(\bbR^{d_1}\times\com{\bbR^{d_1}})),\ 
a_3\in\cS(\bbR^{2d_1+2d_2}).
\label{16.6.2005.4}\end{multline}
Note that the three terms in the amplitude in \eqref{16.6.2005.4} are
really of the same type and the second and third can be included in the
first, except that the support conditions are relaxed to rapid decay at
infinity. Thus, the third class of symbols corresponds to Schwartz
kernels. The second class corresponds 
to Schwartz functions in $z,z'$ with values in
the the classical pseudodifferential operators of order $m'$ on
$\bbR^{d_1}$ and with kernels having bounded support in $y+y'.$ Since these
kernels are actually Schwartz if the singularity at $y=y'$ is cut out, the
effect of the second two terms is simply to admit the kernels which are
Schwartz functions of $z,z'$ with values in the pseudodifferential kernels
of order $m'$ on $\bbR^{d_1}$ with bounded singular support (in $y)$ and
Schwartz tails. Similarly addition of these two terms `completes' the first
term in admitting appropriate tails at infinity to ensure that

\begin{proposition}\label{16.6.2005.5} The operators with kernels as in
  \eqref{16.6.2005.4} act on $\cS(\bbR^{d_1+d_2})$ and form a bifiltered
  algebra with the orders $m,m'\in\bbZ;$ omitting the first term in
  \eqref{16.6.2005.4} gives an ideal, as does omitting the first two
  terms. The filtration is delineated by two symbol maps 
\begin{equation}
\begin{gathered}
\sigma _{m,m'}(A)=\rho ^m\rho _{\ff}^{m'}a_1\big|_{\rho
  =0}\in\CIc(\bbR^{d_1+d_2}\times[\bbS^{d_1+d_2-1};\bbS^{d_1-1}];N_{m,m'})\\
\begin{aligned}
L(A)=(2\pi)^{-d_2}\int e^{i(z-z')\cdot\zeta}
(b_1&+b_2)(y,\frac{z+z'}2,\hat\eta,\zeta )d\zeta,\\
b_1&=(\rho _{\ff}^{m}a_1)\big|_{\rho _{\ff}=0},\ b_2=a_2\big|_{\bbS^{d_1-1}}
\end{aligned}
\end{gathered}
\label{16.6.2005.6}\end{equation}
which are homomorphisms into the algebras of functions and parameterized
 pseudodifferential operators on $\bbR^{d_2}$ respectively.
\end{proposition}

\begin{proof} All these conclusions follow from the standard methods for
  proving the composition formula for pseudodifferential operators on
  Euclidean space, \ie some form of stationary phase. The fact that the two
  symbol maps are homomorphism follows by oscillatory testing.
\end{proof}

Now, these objects can be transferred to a general fibration of compact
manifolds by localization. Thus, the kernels in \eqref{16.6.2005.4} are
smooth outside the \emph{two} submanifolds 
\begin{equation}
\{y=y',z=z'\}\cup\{y=y'\}
\label{16.6.2005.7}\end{equation}
and the singularity is determined by the Taylor series of $a_1$ at the
boundary of $X_{d_1,d_2}$ and the Taylor series of $a_2$ at the boundary of
the ball. Furthermore these singularities are locally determined in the
sense that the singularity on the diagonal near a point $(\bar y,\bar z)$ is
determined by $a_1$ near $(\bar y,\bar z,\eta ,\zeta )$ and by $a_2$ near
$(\bar z,\bar z,\bar y,\eta).$ The singularity near a point on $y=y'$ away
from the diagonal is determined by $a_1$ and $a_2$ near that point $y=\bar
y.$ So, for a fibration we can associate an algebra
of operators using such localizations.

\begin{definition}\label{16.6.2005.8} If
  $\efammap:\efamtot\longleftrightarrow \famtot$ is a fibration of compact
  manifolds without boundary then $\ptP{\efammap}{m',m}{\efamtot;\sE)$ for
  a $\bbZ_2$-graded vector bundle over $\efamtot$ consists of those
  operators $A:\CI(\efamtot;E_+)\longrightarrow \CI(\efamtot;E_-)$ which
  have Schwartz kernels on $\efamtot^2$ which under local trivializations
  of the bundles and densities are matrices of distributions which are
\begin{enumerate}
\item Smooth away from the fibre diagonal.
\item Near a point of the complement of the diagonal in the fibre diagonal
  are given by smooth functions of the fibre variables with values in the
  classical pseudodifferential kernels of order $m'$ on the base.
\item Near a point of the diagonal are of the form \eqref{16.6.2005.4}.
\end{enumerate}
\end{definition}

As a consequence of Proposition~\ref{16.6.2005.5} we then conclude that the
scalar operators in $\ptP{\efammap}{m',m}{\efamtot}$ form a bigraded
algebra for $m,m'\in\bbZ$ and that there are corresponding global symbol
homomorphisms giving surjective, and multiplicative, maps
\begin{equation}
\begin{gathered}
\sigma :\ptP{\efammap}{m',m}{\efamtot}\longrightarrow
\CI([\com{S^*\efamtot};\efammap^*S^*\famtot];N_{m',m})\\
L:\ptP{\efammap}{m',m}{\efamtot}\longrightarrow
\Ps{m}{\pi^*(\efamtot/\famtot};\sE\otimes_+ N_{m'}}.
\end{gathered}
\label{16.6.2005.9}\end{equation}
In the first map, $N_{m',m}$ is a trivial real bundle corresponding to the
coefficients $\rho ^{-m'}$ and $\rho^{-m}$ in the symbols and in the second
sequence, the pseudodifferential operators act on the fibres of the lift of
the fibration from $\famtot$ to $S^*\famtot$ using the  projection
$\pi:S^*\famtot\to\famtot$  and there is an additional
trivial line bundle corresponding to the factor $\rho _{\ff}^{-m'}.$ As
noted, both maps are surjective, but together they are constrained
precisely by the fact that the symbol of $L(A),$ as a family, is given by
$\sigma (A)$ evaluated at the front face of the blow up of $S^*\efamtot.$

Apart from these composition properties, and their natural generalizations
to the case of families over a second fibration, it is important to note
certain inclusions.

First,
\begin{equation}
\Ps{m}{\efamtot;\sE}\subset\ptP{\efammap}{m,m}{\efamtot;\sE}.
\label{16.6.2005.10}\end{equation}
Indeed, the definition of $\ptP{\efammap}{m,m}{\efamtot;\sE}$ is modelled
on one of the standard definitions of the usual pseudodifferential
operators, so it is enough to refer back to the model case. Directly from
the definition of the symbols in \eqref{16.6.2005.4}, the classical
symbols of order $m,$ which are just elements of
$\rho ^{-m}\CIc(\bbR^{d_1+d_2}\times\com{\bbR^{d_1+d_2}}),$ lift to product
type symbols of order $(m,m)$ under the blow up. This leads immediately to
\eqref{16.6.2005.10}. Note that in this case, the symbol in the
product-type sense is just the lift of the symbol in the usual sense
(actually loosing no information by continuity) and the base family is
simply again the symbol, evaluated on the lift of the cosphere bundle from
the base and interpreted as acting as bundle isomorphisms (so local
operators) on the fibres.

The second inclusion is of pseudodifferential operators acting on the
fibres of the fibration. By definition the kernels of these operators are
smooth families in the base variables, with values in the classical
pseudodifferential operators on the fibres; as operators on smooth sections over
the total space they are therefore of the form of a product of a classical
conormal distribution with a delta section on the fibre diagonal. Working locally
this reduces to \eqref{16.6.2005.4} with $a_1$ or $a_2$ actually independent
of $\eta;$ it follows directly that 
\begin{equation}
\Ps m{\efamtot/\famtot;\sE}\subset\ptP{\efammap}{0,m}{\efamtot;\sE}.
\label{faficu.156}\end{equation}
The symbol map reduces to the lift of the symbol on the fibres and the base
family of such an operator is the operator itself.

The third inclusion is simply of pseudodifferential operators on the base
but acting on a bundle $\sE$ which is embedded as a subbundle of
$\CI(\efamtot/\famtot;\sF)$ for some bundle $\sF$ over $\efamtot.$ This
embedding corresponds to a family of smoothing projections of finite rank
$\pi_{\pm}\subset\Ps{-\infty}{\efamtot/\famtot;\sF}$ and the kernel can then
be written, somewhat formally, as $\pi_-\cdot K(A)\cdot\pi_+.$ Again this is
everywhere locally of the form \eqref{16.6.2005.4}, with the fibre part of
order $-\infty$ and it follows that 
\begin{equation}
\Ps m{\famtot;\sE}\subset \ptP{\efammap}{m,-\infty}{\efamtot;\sF}.
\label{faficu.158}\end{equation}
In this case, the operator being of principal order $-\infty,$ the symbol
is zero at any order, but the base family is simply $\pi_-\sigma (A)\pi_+$
acting on sections of $\sF.$

We pass over without extensive comment the extension of this construction
to define families of operators, with respect to an overall fibration, and
more generally suspended families in which the cotangent variables of the
base are symbolic variables in the operators. These latter parameters can
always be incorporated into the operators as the duals of additional
Euclidean (base) variables in which the operators are translation invariant.

\section{Product-type fibred cusp operators}\label{PTFCO}

The discussion above of operators of product-type can be extended to fibred
cusp pseudodifferential operators. In such an extension, as in the subsequent
one to an adiabatic limit, the extension is based on the principle that
the product-type operators above are defined through a geometric class of
distributions, the product-type conormal distributions for the pair of
embedded submanifolds  
\begin{equation*}
\Diag\subset\Diag_{\efammap}.
\label{faficu.159}\end{equation*}
So, to generalize these operators to another setting it is only necessary
to start with a space of operators associated with the conormal
distributions on an embedded submanifold, replacing $\Diag,$ and to
replace these in turn by an appropriate class of product-type distributions.

This is precisely the case with the fibred cusp operators defined and
discussed in \cite{Mazzeo-Melrose4} for a compact manifold with boundary
$\famtot$ with a given fibration of the boundary
$\boumap:\pa\famtot\longrightarrow\boubas$ and a choice of normal
trivialization of the fibres. The latter choice can be represented by a choice of
boundary defining function $x\in\CI(\famtot).$ Then the product $\famtot^2$
on which kernels are normally defined, is replace by a blown-up version of
it. Namely first the corner is blown up
\begin{equation}
\famtot^2_{\bo}=[\famtot^2;(\pa\famtot)^2]
\label{faficu.162}\end{equation}
(when the boundary is not connected all products of pairs of boundary
components should be blown up.) Within the new, or front, face of this
manifold with corners the fibration and choice of normal trivialization
combine to define a submanifold 
\begin{equation}
\Gamma_{\fC{\boumap}}\subset\ff(\famtot^2_{\bo}).
\label{faficu.163}\end{equation}
Namely $\Gamma_{\fC{\boumap}}$ is the fibre diagonal in the boundary
variables intersected with the submanifold $\frac{x-x'}{x+x'}=0$ where it
should be observed that this \emph{is} a smooth function on
$\famtot^2_{\bo}$ and that the resulting submanifold does only depend on
the data giving the fibred cusp structure. Then the kernels for fibred cusp
pseudodifferential operators are simply the standard conormal
distributional sections with respect to the lifted diagonal of an
appropriate smooth bundle over
\begin{equation}
\Diag_{\fC{\boumap}}\subset\famtot_{\fC{\boumap}}^2=
[\famtot^2_{\bo};\Gamma_{\fC{\boumap}}].
\label{faficu.164}\end{equation}
The lifted diagonal is an embedded p-submanifold, \ie has a product type
decomposition at the boundaries. Thus the definition of the kernels as the
conormal distributions (which are also required to vanish to infinite order
at boundary faces not meeting the diagonal) is meaningful.

Here we consider the case of a fibration over the manifold with boundary,
$\efammap:\efamtot\longrightarrow \famtot$ which has compact fibres without
boundary. We suppose that $\famtot$ has a fibred cusp structure as above
and that the fibres of $\efammap$ are to be treated as part of the boundary
fibres, \ie we take the fibred cusp structure on $\efamtot$ given by the
fibration 
\begin{equation}
\eboumap:\pa\famtot\longrightarrow \boubas,\
\eboumap=(\pa\efammap)\circ\boumap 
\label{faficu.165}\end{equation}
and with the normal trivialization given by lifting an admissible defining
function on $\famtot.$

In the construction of $\famtot^2_{\fC{\boumap}},$ each blow up is near the
corners of $\famtot^2$ and the procedure is local with respect to the
fibration of the boundary. That is, if $\boumap$ is locally trivialized to a
product $O\times\boufib,$ consistent with the fibred cusp structure, where
$O$ is an open set in $[0,\infty)\times\bbR^l$ then over the preimage of
this set the stretched product is just 
\begin{equation}
\famtot^2_{\fC{\boumap}}\simeq O^2_{\scat}\times\boufib^2.
\label{faficu.166}\end{equation}
Here, $O^2_{\scat}=O^2_{\fC{\Id}}$ is the stretched product in the 
scattering case, that is, the case where the boundary map is the identity.
Thus, when the extra fibration is added it follows that, again locally near
boundary points and in small open sets $U,U'\subset\boufib$ (so that
$\efammap$ is also locally trivialized) 
\begin{equation}
\efamtot^2_{\fC{\eboumap}}\simeq O^2_{\scat}\times U\times U'\times\efamfib^2.
\label{faficu.167}\end{equation}
In fact $\efamtot_{\fC{\eboumap}}^2$ fibres over $\famtot^2_{\fC{\boumap}}$
with fibre which is $\efamfib^2.$ This shows that the geometric situation
of the lifted diagonal and the lifted fibre diagonal, which is just the
diagonal in $\famtot_{\fC{\boumap}},$ is completely analogous to the
product-type setting discussed above.

\begin{definition}\label{faficu.169} If $\efammap:\efamtot\longrightarrow
  \famtot$ is a fibration with fibres compact manifolds without boundary
  over a compact manifold with fibred cusp structure then the space
  $\ptfcP{\efammap}{\boumap}{m',m}{\efamtot;\sE}$ is given by the space of
  product-type conormal distributions, as in Definition~\ref{16.6.2005.8}, with
  respect to the lifted diagonal and fibred diagonal, and vanishing to
  infinite order at all boundary faces which do not meet these p-submanifolds.
\end{definition}

As shown in \cite{Mazzeo-Melrose4}, the composition properties of fibred
cusp operators follow from geometric considerations. Namely if the product
is interpreted as a push-forward for an appropriately defined triple space
(as in \cite{Mazzeo-Melrose4}) then locally for the fibred cusp calculus
the problem is the same uniformly up to the boundary, and hence follows
from the discussion above. Thus we conclude that 
\begin{equation}
\ptfcP{\efammap}{\boumap}{m',m}{\efamtot;\sE}
\circ\ptfcP{\efammap}{\boumap}{p',p}{\efamtot;\sF}\subset
\ptfcP{\efammap}{\boumap}{m'+p',m+p}{\efamtot;\sG}
\label{faficu.171}\end{equation}
provided the product makes sense, \ie $E_+=F_-,$ $G_+=E_+,$ $G_-=F_-.$
Furthermore there are now three `symbol homomorphisms'. Two of these are
modified versions of the corresponding homomorphism for the fibred cusp
calculus. Thus, the symbol map now takes values, as in \eqref{16.6.2005.9}, in
sections of the appropriate bundle over a blown-up version of the
fibred-cusp cosphere bundle 
\begin{equation}
\sigma:\ptfcP{\efammap}{\boumap}{m',m}{\efamtot}\longrightarrow
\CI([\com{\cuS^*\efamtot};\efammap^*\cuS^*\famtot];N_{m',m}).
\label{faficu.172}\end{equation}
Similarly the indicial operator, which corresponds geometrically to the
restriction of the kernel to the final front face in the blown-up space,
now takes values in product-type and suspended operators on the boundary 
\begin{equation}
N:\ptfcP{\efammap}{\boumap}{m',m}{\efamtot;\sE}\longrightarrow
\sfptP{\boumap}{\efammap}{m,m'}{\pa\efamtot;\sE}.
\label{faficu.173}\end{equation}
The base map can also be defined by oscillatory testing 
\begin{equation}
L:\ptfcP{\efammap}{\boumap}{m',m}{\efamtot;\sE}\longrightarrow
\fcP{\boumap}{m}{\pi^*(\efamtot/\famtot);\sE\otimes_+ N_{m'}}.
\label{faficu.174}\end{equation}
Each of these maps delineates a filtration of the algebra, corresponding to
the order $m,$ the degree of boundary vanishing $x$ and the order $m'.$ They
are separately surjective and jointly, in pairs or all together, subject
only to the natural compatibility conditions 
\begin{equation}
\sigma (N)=\sigma _{\pa},\ \sigma (L)=\sigma_{\ff},\ N(L)=L(N).
\label{faficu.175}\end{equation}

\section{Adiabatic limit in the fibres}\label{Adiabatic}

A notion of the adiabatic limit of pseudodifferential operators was
introduced in \cite{MR90m:58004}. Namely, for a fibration of compact
manifolds, for the moment without boundary,
$\efammap:\efamtot\longrightarrow \famtot,$ one can consider a sense in
which pseudodifferential operators on $\efamtot$ degenerate, by localizing
in the base variables, to become families of pseudodifferential
operators on the fibres of $\efammap.$

To do so, let $\epsilon \in[0,1]$ be the `adiabatic parameter.' If we
consider kernels given by conormal distributions with respect to the
submanifold
\begin{equation}
\Diag\times[0,1]\subset \efamtot^2\times[0,1]
\label{faficu.176}\end{equation}
we simply arrive at the $\epsilon$-parameterized pseudodifferential
operators on $\efamtot.$ On the other hand, if we blow up the fibre
diagonal at $\epsilon =0,$ introducing 
\begin{equation}
\efamtot_{\aD{\efammap}}^2=[\efamtot^2\times[0,1];\Diag_{\efammap}\times\{0\}]
\overset{\beta }\longrightarrow \efamtot^2\times[0,1]
\label{faficu.177}\end{equation}
we obtain a manifold with corners with two important boundary faces (we
ignore $\epsilon =1$ as being `regular'), the `old boundary' being the
proper lift, $\beta^{\#}\{\epsilon =0\}$ and the new `front face' produced
by the blow up. The diagonal has proper lift to a smooth p-submanifold 
\begin{equation}
\beta ^{\#}(\Diag\times[0,1])\subset\efamtot^2_{\aD{\efammap}}
\label{faficu.178}\end{equation}
which meets the boundary only in the front face. The space of adiabatic
pseudodifferential operators is then 
\begin{multline}
\adP{\efammap}m{\efamtot;\sE}=\{K\in
I^{m-\frac14}(\efamtot^2_{\aD{\efammap}},\beta ^{\#}(\Diag\times[0,1]); 
\hom(\sE)\otimes\Omega
_{\aD{\efammap}}); \\ 
K\equiv0\Mat\beta ^{\#}\{\epsilon =0\}\}.
\label{faficu.179}\end{multline}
In fact a neighbourhood of $\beta ^{\#}(\Diag\times[0,1])$ in
$\efamtot^2_{\aD{\efammap}}$ is diffeomorphic to a neighbourhood of
$\Diag\times[0,1]$ in $\efamtot^2\times[0,1]$ so we may legitimately think
of these kernels as having the `same singularities' as the ordinary
pseudodifferential families but with a different action. The (trivial)
density line bundle in \eqref{faficu.179} takes care of factors that arise
even for the identity.

These operators compose in the expected way and have two symbols. The first
is the usual symbol, now taking values in sections of the appropriate
bundle over the `adiabatic cosphere bundle' (which is a bundle over
$\efamtot\times[0,1])$ and giving a short exact sequence
\begin{equation}
\xymatrix@1{\adP{\efammap}{m-1}{\efamtot;\sE}\ar[r]&
\adP{\efammap}m{\efamtot;\sE}\ar[r]^(.35){\sigma }&
\syM{\aD{\efammap}}{m}{\efamtot\times[0,1];\hom(\sE)}.}
\label{faficu.180}\end{equation}
Secondly there is a symbol representing the limit at $\epsilon =0.$ It is a
suspended family of pseudodifferential operators on the fibres of the
fibration with (symbolic) parameters in the rescaled cotangent bundle of
the base 
\begin{equation}
\xymatrix@1{\epsilon\adP{\efammap}m{\efamtot;\sE}\ar[r]&
\adP{\efammap}m{\efamtot;\sE}\ar[r]^(.35){\aD{\efammap}}&
\fsP{\efammap}{m}{\sfamrel;\hom(\sE)}.}
\label{faficu.181}\end{equation}

Now, in the multiplicativity construction in Section~\ref{Multiplicativity}
we use a quite analogous construction to pass from fibred cusp operators
with respect to a fibration and a boundary fibration to fibred cusp
operators for the same fibration and a finer boundary fibration, that is
with smaller fibres. Thus we are `converting' some fibre variables in the
boundary to base variables. The main step is to carry this out on one
fibre, so we can consider the model case of a compact manifold with
boundary $\famfib$ with an iterated fibration of its boundary giving a
commutative diagram
\begin{equation}
\xymatrix{\pa\efamtot \ar[rd]_{\eboumap}\ar[rr]^{\boumap}&&\boubas,\\
&\eboubas\ar[ru]}
\label{faficu.182}\end{equation}
so the map from $\eboubas$ to $\boubas$ is also a fibration.

The construction of the fibred cusp calculus on $\boufib$ is briefly
discussed above. Adding an adiabatic parameter $\epsilon \in[0,1]$ we can
consider the $\epsilon$-parameterized fibred cusp calculus with respect to
$\eboumap$ (and some choice of fibred cusp structure) by replacing
\eqref{faficu.164} by 
\begin{equation}
\famfib_{\fC{\boumap}}^2\times[0,1]=
[\famfib^2_{\bo}\times[0,1];\Gamma_{\fC{\boumap}}\times[0,1]].
\label{faficu.183}\end{equation}
Now, at $\epsilon=0$ we consider the finer boundary fibration given by
$\eboumap,$ with a consistent fibred cusp structure (meaning such that there
is a global boundary defining function which induces both). Note that the
corresponding `lifted' fibre diagonal in the boundary is then a
p-submanifold 
\begin{equation}
\Gamma_{\fC{\eboumap}}\subset \Gamma_{\fC{\boumap}}\Longrightarrow
\Gamma_{\fC{\eboumap}}\cap\{\epsilon =0\}\subset
\Gamma_{\fC{\boumap}}\times[0,1].
\label{faficu.184}\end{equation}
In particular this means that the proper lift of this manifold under the
blow up in \eqref{faficu.183} is again a p-submanifold which can be blown
up, giving the new compact manifold with corners 
\begin{multline}
\famfib_{\fC{\boumap},\aD{\eboumap}}^2=[\famfib_{\fC{\boumap}}^2\times[0,1];
\Gamma_{\fC{\eboumap}}\cap\{\epsilon =0\}]\\
=[\famfib^2_{\bo}\times[0,1];\Gamma_{\fC{\boumap}}\times[0,1],
\Gamma_{\fC{\eboumap}}\cap\{\epsilon =0\}].
\label{faficu.185}\end{multline}

The kernels of the operators we consider will be defined on this
manifold. Ignoring $\epsilon =1$ as we shall, there are three boundary
faces which meet the proper lift of the diagonal (with a factor of
$[0,1]),$ which as usual is a p-submanifold. Two of these are the boundary
hypersurfaces produced by the blow-ups in \eqref{faficu.185}, the first is
essentially the front face of the $\Phi$-fibred cusp calculus with an
extra factor of $[0,1]$ but also blown up at $\epsilon =0$ corresponding to
the $\tPhi$-fibred cusp calculus. The second is the front face produced by
the last blow-up over $\epsilon =0$ and the third is the proper lift
of $\epsilon =0,$ this is simply the manifold with corners
$\boufib^2_{\fC{\eboumap}}.$

Now, the space of $\tPhi$-adiabatic, $\Phi$-fibred cusp pseudodifferential
operators on $\famfib,$ $\adfcP{\boumap}{\eboumap}{m}{\famfib}$ is identified
with the conormal distributions, of an appropriate density bundle, on
$\famfib^2_{\fC{\boumap},\aD{\eboumap}}$ with respect to the proper lift of
the diagonal and vanishing to infinite order at all boundary faces which do
not meet this lift. These operators map $\CI(\efamtot_{\ad};E_+),$ to
$\CI(\efamtot_{\ad};E_-)$ where 
\begin{equation}
\efamtot_{\ad}=[\efamtot\times[0,1];\pa\efamtot\times\{0\}]
\label{faficu.187}\end{equation}
and compose in the usual way.

Corresponding to this definition and the discussion above of boundary
hypersurfaces, there are four `symbol maps', the symbol, the indicial
operator, the adiabatic symbol and the adiabatic operator. 

The first two of these correspond to the symbol for the $\boumap$-fibred cusp
calculus and its indicial operator, with dependence on the parameter $\epsilon$
with a change of uniformity at $\epsilon =0.$ The third symbol is the real
transition between the $\boumap$- and $\eboumap$-fibred cusp calculi, and the
last is simply the limiting $\eboumap$-fibred cusp operator itself. More
precisely, the `usual' symbol becomes 
\begin{equation}
\xymatrix@1{
\adfcP{\boumap}{\eboumap}{m-1}{\efamtot;\sE}\ar[r]&
\adfcP{\boumap}{\eboumap}m{\efamtot;\sE}\ar[r]^(.45){\sigma }&
\syM{\aD{\eboumap},\fC{\boumap}}{m}{\efamtot_{\ad};\sE}.}
\label{faficu.186}\end{equation}
Note that the symbols here are section of the tensor product of $\hom(\sE)$
with a density bundle over the modified cosphere bundle to
$\efamtot_{\ad}.$ This coshere bundle is associated to the cotangent bundle
dual to the tangent
bundle of $\efamtot,$ rescaled near the boundary, with locally generating
sections over $\efamtot\times[0,1]$ near a boundary point of the corner given
by the
vector fields (which lift to be smooth on $\efamtot_{\ad})$ 
\begin{equation*}
x^2\pa_x,\ x\pa_y,\ (x^2+\epsilon ^2)^{\frac12}\pa_{y'},\ \pa_z
\label{faficu.188}\end{equation*}
where $x$ is the normal variable, the $y$'s are base variables for
$\boumap,$ the $y'$'s are the extra base variables for $\eboumap$ (so fibre
variables for $\boumap)$ and the $z$'s are fibre variables for $\eboumap.$

The normal, or indicial operator corresponds to the restriction of the
kernel to the front face of $\efamtot^2_{\fC{\boumap}}\times[0,1]$ after
part of its boundary is blown up in the last step in the construction,
\eqref{faficu.185}. Since this can be related directly to the adiabatic
construction for the fibre calculus of $\boumap$ with respect to $\eboumap$
it becomes a map into the corresponding adiabatic (and suspended) calculus 
\begin{equation}
\xymatrix@1{
\tilde x\adfcP{\boumap}{\eboumap}{m}{\efamtot;\sE}\ar[r]&
\adfcP{\boumap}{\eboumap}m{\efamtot;\sE}\ar[r]^(.45){N}&
\sadP{\boumap}{\eboumap}{m}{\pa \efamtot;\sE}.}
\label{faficu.189}\end{equation}
Note that the factor $\tilde x\in\CI(\efamtot_{\ad})$ is a defining function for
the boundary after the last blow-up, it can be taken to be $x(x^2+\epsilon
^2)^{-\frac12}.$

The transitional, adiabatic normal operator corresponds to the restriction
of the kernel to the front face produced in the final blow-up in
\eqref{faficu.185} and takes values in a suspended space of
pseudodifferential operators on the fibres of $\eboumap$ giving a short
exact sequence 
\begin{equation}
\xymatrix@1{
(x^2+\epsilon ^2)^{\frac12}\adfcP{\boumap}{\eboumap}{m}{\efamtot;\sE}\ar[r]&
\adfcP{\boumap}{\eboumap}m{\efamtot;\sE}\ar[r]^(.45){\ad}&
\Psi^{m}_{\eboumap-\sus(V)}(\pa\efamtot;\sE)},
\label{faficu.190}\end{equation}
where $V\to [-1,1]_{\tau}\times\eboubas$ is 
is the null bundle of the restriction 
${}^{\eboumap}T_{\pa\efamtot}{\famtot}\to T_{\pa\efamtot}\efamtot$ and 
\begin{equation*}
   \tau= \frac{x-\epsilon}{x+\epsilon} \in [-1,1]
\label{faficu.224}\end{equation*}
is a variable on the adiabatic front face.
Finally, the fourth map, into the finer fibred cusp calculus gives a short
exact sequence 
\begin{equation}
\xymatrix@1{
\tilde\epsilon\adfcP{\boumap}{\eboumap}{m}{\efamtot;\sE}\ar[r]&
\adfcP{\boumap}{\eboumap}m{\efamrel;\sE}\ar[r]^(.55){A}&
\fcP{\eboumap}{m}{\efamtot;\sE}.}
\label{faficu.191}\end{equation}
Here $\tilde\epsilon=\epsilon (x^2+\epsilon ^2)^{-\frac12}$ is a defining
function for the limiting boundary in $\efamtot_{\ad}.$

Notice that
\begin{equation}
P\in\adfcP{\boumap}{\eboumap}m{\efamrel;\sE},\ N(P)=0,\
\ad(P)=0\Longleftrightarrow  P\in x\adfcP{\boumap}{\eboumap}m{\efamrel;\sE}.
\label{faficu.195}\end{equation}
The various compatibility conditions are then given by
\begin{equation}
\begin{gathered}
A\circ\sigma(P)= \sigma\circ A(P), \quad \sigma\circ N(P)= N\circ \sigma (P),
\quad \ad\circ \sigma(P)= \sigma\circ \ad(P), \\
\ad\circ N(P)= \left. \ad(P)\right|_{\tau=-1}\quad 
\mbox{and} \quad N\circ A(P)= \left. \ad(P)
\right|_{\tau=1} .
\end{gathered} 
\label{20.7.2005.2}\end{equation}
In particular there is no compatibility condition between $A(P)$ and $N(P).$

As noted above there is in fact a fifth map corresponding to restriction to
$\epsilon =1$ (or really any positive value of $\epsilon)$ 
\begin{equation}
\xymatrix@1{
(1-\epsilon)\adfcP{\boumap}{\eboumap}{m}{\efamtot;\sE}\ar[r]&
\adfcP{\boumap}{\eboumap}m{\efamrel;\sE}\ar[r]^(.55){\big|_{\epsilon =1}}&
\fcP{\boumap}{m}{\efamtot;\sE}.}
\label{faficu.193}\end{equation}
In some sense, $P\in \adfcP{\boumap}{\eboumap}m{\efamtot;\sE}$ should be 
interpreted as a homotopy between $A(P)\in \fcP{\eboumap}{m}{\efamtot;\sE}$
and $P\big|_{\epsilon =1}\in \fcP{\boumap}{m}{\efamtot;\sE}.$  When both
$A(P)$ and $P\big|_{\epsilon =1}$ are Fredholm operators, one would expect them
to have the same index provided $P$ is a homotopy through Fredholm operators.
The next proposition makes this statement more precise.

\begin{proposition}\label{faficu.192} If
$P\in\adfcP{\boumap}{\eboumap}m{\efamtot;\sE}$ exists with $\sigma (P),$
$N(P)$ and $\ad(P)$ invertible in their pseudodifferential (or symbol)
calculi then $P\big|_{\epsilon =1}$ and $A(P)$ are both fully elliptic
and have the same (families) index.
\end{proposition}

\begin{proof} The compatibility conditions between $\sigma(P),$ $N(P)$ and
  $\ad(P)$ are such that if each is invertible within the calculus of
  pseudodifferential operators then the inverses are compatible. Thus,
  under this condition a parametrix can be constructed for $P.$ As usual,
  a symbolic parametrix can be improved to a full parametrix in
  the sense that $Q\in\adfcP{\boumap}{\eboumap}{-m}{\efamtot;\sE^-}$ satisfies
\begin{equation}
QP-\Id\in x^\infty\adfcP{\boumap}{\eboumap}{-\infty}{\efamtot;E_+},
\ PQ-\Id\in x^\infty\adfcP{\boumap}{\eboumap}{-\infty}{\efamtot;E_-}.
\label{faficu.194}\end{equation}
Note that we do \emph{not} achieve vanishing at $\epsilon =0$ since we have
  not assumed that $A(P)$ is invertible. In fact the errors in
  \eqref{faficu.194} are just smoothing operators, smooth in $\epsilon$ and
  with kernels vanishing to infinite order at the boundary. Following the
  discussion of the index in Section~\ref{Analyticindex} the index bundle
  (for a family of such operators) can be stabilized to a bundle over
  $[0,1]\times\fambas.$ Indeed, this follows from the fact that  
\begin{equation*}
     x^\infty\adfcP{\boumap}{\eboumap}{-\infty}{\efamtot;E_+}=
\CI([0,1], \dot{\Psi}^{-\infty}(\efamtot;E_{+}))
\label{faficu.225}\end{equation*}
where $\dot{\Psi}^{-\infty}(\efamtot;E_{+})=
x^{\infty}\Psi_{\fC{\boumap}}^{-\infty}(\efamtot;E_{+})$ does not depend on 
the choice of the boundary fibration $\boumap.$ 
It follows from this that the families index at
  $\epsilon =0$ and $\epsilon=1$ 
are the same and the former is the index bundle
  for $A(P),$ the latter for $P\big|_{\epsilon =1}.$
 \end{proof}

\providecommand{\bysame}{\leavevmode\hbox to3em{\hrulefill}\thinspace}
\providecommand{\MR}{\relax\ifhmode\unskip\space\fi MR }
\providecommand{\MRhref}[2]{%
  \href{http://www.ams.org/mathscinet-getitem?mr=#1}{#2}
}
\providecommand{\href}[2]{#2}

\end{document}